%% file: algorithmic.tex
\documentclass{amsart}

\usepackage[english]{babel}
\usepackage{csquotes}
\usepackage[style=alphabetic,
            isbn=false,
            doi=false,
            url=false,
            backend=bibtex,
            maxnames=5
           ]{biblatex}
\bibliography{sources}
\usepackage{amsmath,amsfonts,amsthm,mathtools, amssymb}
\usepackage{thmtools}
\usepackage{xcolor}
\usepackage{tikz-cd}
\usepackage{hyperref} 
\usepackage{aliascnt}
\usepackage{xparse}
\usepackage{tensor}
\usepackage{caption}

\usepackage{tikz}
\usetikzlibrary{matrix,arrows,calc}
\usetikzlibrary{positioning}
\usetikzlibrary{decorations.pathreplacing,decorations.markings,decorations.pathmorphing}
\usetikzlibrary{positioning,arrows,patterns}
\usetikzlibrary{cd}
\usetikzlibrary{intersections}

\tikzset{
	every loop/.style={very thick},
	comp/.style={circle,fill,black,,inner sep=0pt,minimum size=5pt},
	order bottom left/.style={pos=.05,left,font=\tiny},
	order top left/.style={pos=.9,left,font=\tiny},
	order bottom right/.style={pos=.05,right,font=\tiny},
	order top right/.style={pos=.9,right,font=\tiny},
	order node dis/.style={text width=.75cm},
	circled number/.style={circle, draw, inner sep=0pt, minimum size=12pt},
	below left with distance/.style={below left,text height=10pt},
    below right with distance/.style={below right,text height=10pt}
	}

\usepackage{listings}
\usepackage{textcomp}
\usepackage{color}

\usepackage{microtype}
\makeatletter\g@addto@macro\@verbatim{\microtypesetup{activate=false}}\makeatother%

\definecolor{mygreen}{rgb}{0,0.6,0}
\definecolor{mygray}{rgb}{0.5,0.5,0.5}
\definecolor{mymauve}{rgb}{0.58,0,0.82}
\definecolor{backcolour}{rgb}{0.95,0.95,0.92}

\newcommand{\BIC}{\operatorname{BIC}}

\makeatletter
\ifcsname phantomsection\endcsname
    \newcommand*{\@gobblenexttocentry}[9]{}
\else
    \newcommand*{\@gobblenexttocentry}[4]{}
\fi
\newcommand*{\addsubsection}{%
    \addtocontents{toc}{\protect\@gobblenexttocentry}%
    \subsection*}
\makeatother

\include{macros}

\begin{document}

\def\subsectionautorefname{Section}
\def\subsubsectionautorefname{Section}
\def\sectionautorefname{Section}
\def\equationautorefname~#1\null{(#1)\null}

\def\optmark{*}

\title[\MakeLowercase{\texttt{diffstrata}} -- a Sage package]
{\MakeLowercase{\texttt{diffstrata}} -- a Sage package for calculations in the tautological ring of the moduli space of Abelian differentials}

\author{Matteo Costantini}
\email{costanti@math.uni-bonn.de}
\address{Institut f\"ur Mathematik, Universit\"at Bonn,
Endenicher Allee 60,
53115 Bonn, Germany}
\author{Martin M\"oller}
\author{Jonathan Zachhuber}
\email{zachhuber@math.uni-frankfurt.de}
\address{Institut f\"ur Mathematik, Goethe-Universit\"at Frankfurt,
Robert-Mayer-Str. 6--8,
60325 Frankfurt am Main, Germany}
\email{moeller@math.uni-frankfurt.de}
\thanks{Research of the second and third author is supported
by the DFG-project MO 1884/2-1 and by the LOEWE-Schwerpunkt
``Uniformisierte Strukturen in Arithmetik und Geometrie''}

\begin{abstract}
The boundary of the multi-scale differential compactification of strata of abelian differentials admits an explicit combinatorial description.
However, even for low-dimensional strata, the complexity of the boundary requires use of a computer.
We give a description of the algorithms implemented in the SageMath package \texttt{diffstrata} to enumerate the boundary components and perform intersection theory on this space. In particular, the package can compute the Euler characteristic of strata using the methods developed by the same authors in \cite{strataEC}.
\end{abstract}

\maketitle
\tableofcontents

\input{sec_intro}

\input{sec_interface}

\input{sec_genstrata}

\input{sec_levelgraphs}

\input{sec_generating}

\input{sec_AGs}

\input{sec_mult}

\input{sec_gluing}

\input{sec_caching}

\input{sec_tests}

\printbibliography
\end{document}

%% file: macros.tex
\newcommand{\mynewtheorem}[4]{
  \if\relax\detokenize{#3}\relax %
    \if\relax\detokenize{#4}\relax %
	\declaretheorem[name=#2]{#1}
    \else
	\declaretheorem[numberwithin=#4,name=#2]{#1}
    \fi
  \else
	\declaretheorem[sibling=#3,name=#2]{#1}
  \fi
}
\makeatletter
 \providecommand{\@dotsep}{5}
\makeatother
\mynewtheorem{theorem}{Theorem}{}{section}
\theoremstyle{definition}
\mynewtheorem{lemma}{Lemma}{theorem}{}
\mynewtheorem{rem}{Remark}{lemma}{}
\mynewtheorem{prop}{Proposition}{lemma}{}
\mynewtheorem{cor}{Corollary}{lemma}{}
\mynewtheorem{definition}{Definition}{lemma}{}
\mynewtheorem{question}{Question}{lemma}{}
\mynewtheorem{assumption}{Assumption}{lemma}{}
\mynewtheorem{example}{Example}{lemma}{}
\mynewtheorem{algorithm}{Algorithm}{lemma}{}

\def\defbb#1{\expandafter\def\csname b#1\endcsname{\mathbb{#1}}}
\def\defcal#1{\expandafter\def\csname c#1\endcsname{\mathcal{#1}}}
\def\deffrak#1{\expandafter\def\csname frak#1\endcsname{\mathfrak{#1}}}
\def\defop#1{\expandafter\def\csname#1\endcsname{\operatorname{#1}}}
\def\defbf#1{\expandafter\def\csname b#1\endcsname{\mathbf{#1}}}
\def\deftt#1{\expandafter\def\csname#1\endcsname{\texttt{#1}}}

\makeatletter
\def\defcals#1{\@defcals#1\@nil}
\def\@defcals#1{\ifx#1\@nil\else\defcal{#1}\expandafter\@defcals\fi}
\def\deffraks#1{\@deffraks#1\@nil}
\def\@deffraks#1{\ifx#1\@nil\else\deffrak{#1}\expandafter\@deffraks\fi}
\def\defbbs#1{\@defbbs#1\@nil}
\def\@defbbs#1{\ifx#1\@nil\else\defbb{#1}\expandafter\@defbbs\fi}
\def\defbfs#1{\@defbfs#1\@nil}
\def\@defbfs#1{\ifx#1\@nil\else\defbf{#1}\expandafter\@defbfs\fi}
\def\defops#1{\@defops#1,\@nil}
\def\@defops#1,#2\@nil{\if\relax#1\relax\else\defop{#1}\fi\if\relax#2\relax\else\expandafter\@defops#2\@nil\fi}
\def\deftts#1{\@deftts#1,\@nil}
\def\@deftts#1,#2\@nil{\if\relax#1\relax\else\deftt{#1}\fi\if\relax#2\relax\else\expandafter\@deftts#2\@nil\fi}
\makeatother

\defbbs{ZHQCNPALRVW}
\defcals{OPQMNXYLTRAEHZKCFI}
\deffraks{apijklmnopqueR}
\defops{PGL,SL,Sp,mod,Spec,Re,Gal,Tr,End,GL,Hom,PSL,H,div,Aut,rk,Mod,R,T,Tr,Mat,Vol,MV,Res,vol,Z,diag,Hyp,ord,Im,ev,U,dev,c,CH,fin,pr,Pic,lcm,ch,td,LG,id,Sym,Aut}
\defbfs{uvzwp} %
\deftts{LevelGraph,EmbeddedLevelGraph,sage}

\def\i{\mathrm{i}}
\def\e{\mathrm{e}}

\def\ol{\overline}

\def\be{\begin{equation}}   \def\ee{\end{equation}}     \def\bes{\begin{equation*}}    \def\ees{\end{equation*}}
\def\ba{\be\begin{aligned}} \def\ea{\end{aligned}\ee}   \def\bas{\bes\begin{aligned}}  \def\eas{\end{aligned}\ees}
\def\={\;=\;}  \def\+{\,+\,} \def\m{\,-\,}

\newcommand*{\proj}{\mathbb{P}}
\newcommand{\barmoduli}[1][g]{{\overline{\mathcal M}}_{#1}}

\newcommand{\omoduli}[1][g]{{\Omega\mathcal M}_{#1}}

\DeclareDocumentCommand{\LMS}{ O{\mu} O{g,n} O{}}{\Xi\overline{\mathcal{M}}^{#3}_{#2}(#1)}
\DeclareDocumentCommand{\Romod}{ O{\mu} O{g,n} O{}}{\Omega\mathcal{M}^{#3}_{#2}(#1)}

\newcommand{\bfg}{{\bf g}}

\newcommand{\bfn}{{\bf n}}

\newcommand{\bfmu}{{\boldsymbol{\mu}}}

\newcommand{\wh}{\widehat}

\newlength{\halfbls}

%% file: sec_intro.tex
\section{Introduction}

An important tool in studying the moduli space of (meromorphic) abelian differentials $\bP\Omega\cM_g(\mu)$ is its modular compactification, the space of multi-scale differentials $\LMS$ for $\mu$ an integer partition of $2g-2$.
The boundary is of a combinatorial nature, parametrised, for any $\mu$, by finitely many labeled \emph{level graphs} \cite{LMS}.
However, already listing isomorphism classes of these graphs is a non-trivial task, and already for $g=3$ the number of components becomes so large that even listing them is unfeasible by hand.

Moreover, in \cite{strataEC} the tautological ring of $\LMS$ is described and calculations therein may be expressed purely in terms of the combinatorics of the boundary. 
Again, even the simplest calculations are extremely cumbersome to perform without the assistance of a computer.

The \verb|diffstrata| package provides a framework for calculations in the tautological ring of $\LMS$.
It is implemented in \verb|sage| \cite{sage} and is inspired by the package \verb|admcycles| \cite{admcycles} for calculations in the tautological ring of $\overline{\cM}_{g,n}$.
However, due to the differences in the structure of the boundary, the implementation and interface of the two packages have only little in common.
{A key step in the evaluation process is performed by \texttt{admcycles}, though.}

\verb|diffstrata| may be used naively for basic inquiries about strata. The following example contains the complete code required to calculate the (orbifold) Euler characteristic of $\bP\Omega\cM_{2}(2)$ (using \cite{strataEC}):
\begin{lstlisting}
sage: from admcycles.diffstrata import *
sage: X=Stratum((2,))
sage: X.euler_characteristic()
-1/40
\end{lstlisting}
We can also easily display information on the compactification of a stratum:
\begin{lstlisting}
sage: from admcycles.diffstrata import *
sage: X=Stratum((2,2))
sage: X.info()
Stratum: (2, 2)
with residue conditions: []

Genus: [3]
Dimension: 6
Boundary Graphs (without horizontal edges):
Codimension 0: 1 graph
Codimension 1: 20 graphs
Codimension 2: 86 graphs
Codimension 3: 147 graphs
Codimension 4: 110 graphs
Codimension 5: 30 graphs
Total graphs: 394
\end{lstlisting}
Moreover, we can compute the Masur--Veech volume as an intersection number of $\xi$-powers and $\psi$-classes, cf. \cite{CMSZ}:
\begin{lstlisting}
sage: from admcycles.diffstrata import *
sage: X=Stratum((1,1))
sage: (X.xi^2 * X.psi(1) * X.psi(2)).evaluate()
-1/720
sage: (X.xi^3 * X.psi(1)).evaluate()
-1/360
\end{lstlisting}
However, more extensive calculations necessarily require a deeper understanding of the syntax and the underlying objects of the package.

The most fundamental notion is that of an \emph{(enhanced) profile}, encoding isomorphism classes of level graphs using tuples of integers (essentially writing them as products of divisors).
It is explained in detail in \autoref{sec:AllGraphs} together with the algorithms used to list all boundary components.
This notion is used by \verb|diffstrata|, together with an encoding of $\psi$-polynomials, to encode all additive generators of the tautological ring and all calculations performed involve (formal) sums of these.

Implementing the recursive structure of the boundary, even when considering holomorphic strata, disconnected meromorphic strata with residue conditions will appear as levels of level graphs. 
To encode these, we use \emph{Generalised Strata}, as introduced in \cite{strataEC}. 
These may be thought of as container objects, where the graphs are stored and the calculations performed.

However, implementing any formula in the tautological ring will require some understanding of the subtleties surrounding enhanced profiles and the degeneration of level graphs, as well as extracting levels from graphs and working with these.
In \autoref{subsec:ec} we discuss, as an example, the implementation of the formula for the Euler characteristic using level-wise evaluations of top powers of the tautological class $\xi$, cf. \cite{strataEC}.

\addsubsection{Algorithmic aspects}

Implementing the compactification $\LMS$ poses a series of challenges not encountered in the boundary of $\overline{\cM}_{g,n}$.
First, constructing all level graphs satisfying precisely the conditions of \cite[Def.~1.1]{LMS} is a non-trivial problem.
While generating a codimension-one degeneration of a stable graph inside $\overline{\cM}_{g,n}$ is straight-forward (add an edge either as a loop or by splitting a vertex subject to the stability condition), for a \emph{level graph} this is the more subtle problem of \emph{splitting a level}.

To solve it, we must construct all divisors in any \emph{generalised} stratum, i.e. meromorphic, disconnected and with residue conditions, as these appear as levels already inside low-genus holomorphic strata. Since all the calculations
of the Chern classes of the logarithmic cotangent bundle and in particular
of the Euler characteristic (\cite{strataEC}) and also the computation of
Masur-Veech volumes (\cite{CMSZ}) happen in the tautological ring defined
by clutching
of \emph{non-horizontal divisors}, the \texttt{diffstrata} package treats
exclusively graphs without horizontal edges. See \cite[Section~8]{strataEC}
for other candidates of tautological rings that are potentially
larger but might actually agree with the tautological ring used here.
\par
We usually refer to two-level graphs for brevity as BICs (\emph{bicoloured
  graphs}, as in \cite{FP}).

The combinatorics arising from the distribution of point orders, level structures, genera and edges leads to the simplest approach, namely considering all enhanced level structures on stable graphs in $\overline{\cM}_{g,n}$ (e.g. using \verb|admcycles|), being too slow even for low-dimensional holomorphic strata.
Instead, we describe an algorithm for directly generating all BICs in a generalised stratum in detail in \autoref{subsec:bicgeneration}.

Moreover, this requires checking the Global Residue Condition (GRC), discovered in \cite{BCGGM}, for generalised strata, i.e. implementing the $\frakR$-GRC \cite[\S4]{strataEC}.
For this, we refine the combinatorial criterion of \cite{MUW} to work for the $\frakR$-GRC in \autoref{sec:GenStrata}.

Having generated all BICs, this allows us to construct any graph by recursively  clutching BICs (one for each level-crossing).
As a by-product, this yields discrete coordinates for the enhanced level graphs, as any graph splits uniquely into a product of distinct BICs: we therefore number the BICs of a stratum and associate to each level graph its \emph{profile}, the tuple of indices of BICs appearing as levels of this graph.
Unfortunately, this is not always injective: profiles may be reducible and we need an \emph{enhanced profile} to refer to a graph uniquely, see \autoref{subsec:degenerations} for details.

However, when multiplying two tautological classes, we require an implementation of the excess intersection formula for $\LMS$ \cite[\S 8]{strataEC} and this requires a good understanding of the degeneration graph of the boundary components.
Here again the notion of profile is key: it allows us to efficiently determine which graphs appear as degenerations and is essential in the computation of normal bundles and general intersections, see \autoref{sec:mult}.

Finally, even in low-genus holomorphic strata, calculations involving top-classes of the tautological ring are only feasible because we use extensive caching.
Again, our ability to replace graphs by enhanced profiles is essential, but already the higher-dimensional strata in genus three require more techniques to become manageable, see \autoref{sec:caching}.
Also, several computed values are written to files so that they can be easily reused between \verb|sage| sessions and may be easily precomputed on another machine and imported into the current session, see \autoref{subsec:files}.

\addsubsection{Open questions}

For the moduli space of (pointed) curves~$\barmoduli[g,n]$ enumerating the boundary
strata, also known as tropical curves, and providing a tight estimate for their
growth rate has been discussed at various places (\cite{ChanCombTT},
\cite{MPstablegraph}), but we are aware of a complete solution only in genus
zero (\cite{McMgenuszero}).
It seems interesting to address the analogous problem in $\LMS$, to count the
number of boundary strata of at least provide tight estimates for
their number. Some tables are given in Section~\ref{sec:AllGraphs} for the
number of graphs without horizontal nodes. Such estimates
would also be the basis for serious runtime analysis of the algorithms in
\verb|diffstrata|.\footnote{In practice, one runs out of memory before speed becomes a serious issue: Calculations in genus four take several hours but need several TB of RAM.}
\par
Related to this problem is the question of determining the number of components
of enhanced profiles for a given profile. There are obvious coarse upper bounds
for this number, but in practice the number is relatively small in average
(for fixed~$\mu$ and codimension). Can this be proven? On the algorithmic
side, it would be useful
to design an algorithm that implements directly the degeneration of enhanced
profiles. This could avoid the generation of the graphs altogether an lead
to significant speed-up.
\par
The alternating sum of the number of enhanced profiles appear to be zero
for all holomorphic strata (we thank A.~Neitzke for observing this during
a talk!). Can this be proven? Moreover, the same appears to hold for the number of
profiles, although the degeneration process of this is a main source of
the complexity of the algorithm, see \autoref{subsec:degenerations}.
\addsubsection{Installation}

The package \verb|diffstrata| is included with \verb|admcycles| version~$1.1$ or greater. See \cite{admcycles} for a detailed guide to installing it.

From now on, all examples will assume that the line
\begin{lstlisting}
sage: from admcycles.diffstrata import *
\end{lstlisting}
has been executed!

\addsubsection{Structure}

We begin by giving an overview of the package interface for the ``casual user'' in \autoref{sec:interface}, before explaining the implementation in more detail.
In \autoref{sec:GenStrata} we start by reviewing some of the mathematical background. 
In \autoref{sec:LG}, we explain the fundamental objects of \verb|diffstrata| and how they relate to their mathematical counterparts. The algorithms to construct all non-horizontal level graphs as well as determining the degeneration graph of a stratum and computing isomorphisms of graphs are described in \autoref{sec:AllGraphs}. \autoref{sec:AGs} explains how \verb|diffstrata| encodes and evaluates (using \verb|admcycles|) tautological classes on strata and in \autoref{sec:mult} the implementation of the excess intersection formula is discussed. In \autoref{clutching} the subtleties around splitting graphs around a level and clutching are discussed and the recursive evaluation of $\xi$ on levels is explained. Finally, in \autoref{sec:caching} we explain how and which computations are cached by \verb|diffstrata| and how to import pre-computed values and end in \autoref{sec:tests} by illustrating how \verb|diffstrata| calculates Euler characteristics and giving a few examples of cross-checks and tests.

\addsubsection{Acknowledgements}

We thank Vincent Delecroix and Johannes Schmitt for many helpful suggestions for implementing this package. We are also grateful to the
Mathematical Sciences Research Institute (MSRI, Berkeley)
and the Hausdorff Institute for Mathematics (HIM, Bonn),
where significant progress on this project was made during their programs
and workshops. The authors thank the MPIM, Bonn, for hospitality
and support for \cite{sage} computations. 
\par
\medskip

\def\listtheoremname{List of Algorithms}
\listoftheorems[ignoreall,show={algorithm}]

%% file: sec_interface.tex
\section{Basic Interface}
\label{sec:interface}

Before discussing the implementation details of \verb|diffstrata|, we briefly revisit the examples of the introduction. The first step is always generating a stratum:
\begin{lstlisting}
sage: X=Stratum((2,))
sage: print(X)
Stratum: (2,)
with residue conditions: []
\end{lstlisting}
Here we have defined a \verb|Stratum| object. The argument is a Python \verb|tuple| and may contain integers that sum to $2g-2$ to define any meromorphic stratum (note the trailing \verb|,| if there is only one entry!).
The \verb|print| statement displays information about the object and hints that this is in fact an instance of a \verb|GeneralisedStratum| that can be disconnected and have residue conditions at the poles, see \autoref{sec:LG} for a more detailed discussion.

Creating a \verb|Stratum| automatically performs a series of calculations. For example, all non-horizontal divisors (BICs) are generated and can now be accessed through~\verb|X|:
\begin{lstlisting}
sage: X.bics
[EmbeddedLevelGraph(LG=LevelGraph([1, 1],[[1], [2, 3]],[(1, 3)],{1: 0, 2: 2, 3: -2},[0, -1],True),dmp={2: (0, 0)},dlevels={0: 0, -1: -1}),
 EmbeddedLevelGraph(LG=LevelGraph([1, 0],[[1, 2], [3, 4, 5]],[(1, 4), (2, 5)],{1: 0, 2: 0, 3: 2, 4: -2, 5: -2},[0, -1],True),dmp={3: (0, 0)},dlevels={0: 0, -1: -1})]
\end{lstlisting}
This illustrates how \verb|diffstrata| represents level graphs internally.
The multitude of decorations makes the classes \verb|EmbeddedLevelGraph| and \verb|LevelGraph| a bit unwieldy and there should be little reason to enter them by hand. But they do store the essential information, they are a backbone of \verb|diffstrata|, and they appear frequently in the output. Details are described in \autoref{sec:LG}.

For a single \verb|EmbeddedLevelGraph|, e.g. an element of \verb|X.bics|, we may use its \verb|explain| method to produce a human-readable description of the graph:
\begin{lstlisting}
sage: X.bics[0].explain()
LevelGraph embedded into stratum Stratum: (2,)
with residue conditions: []
 with:
On level 0:
* A vertex (number 0) of genus 1
On level 1:
* A vertex (number 1) of genus 1
The marked points are on level 1.
More precisely, we have:
* Marked point (0, 0) of order 2 on vertex 1 on level 1
Finally, we have one edge. More precisely:
* one edge between vertex 0 (on level 0) and vertex 1 (on level 1) with prong 1.
\end{lstlisting}
Instead of entering graphs by hand, we should always use (enhanced) profiles to refer to them inside \verb|X|. For BICs, this is simply their index in \verb|X.bics|. We can list all profiles of a given length:
\begin{lstlisting}
sage: X.enhanced_profiles_of_length(2)
(((0, 1), 0),)
sage: X.enhanced_profiles_of_length(3)
()
\end{lstlisting}
Note that there are no profiles of length $3$ even though
\begin{lstlisting}
sage: X.dim()
3
\end{lstlisting}
The reason is that \verb|diffstrata| ignores all graphs with horizontal edges in the boundary. We can also retrieve the \verb|EmbeddedLevelGraph| from an (enhanced) profile:
\begin{lstlisting}
sage: X.lookup_graph((0,1))
EmbeddedLevelGraph(LG=LevelGraph([1, 0, 0],[[1], [2, 3, 4], [5, 6, 7]],[(1, 4), (2, 6), (3, 7)],{1: 0, 2: 0, 3: 0, 4: -2, 5: 2, 6: -2, 7: -2},[0, -1, -2],True),dmp={5: (0, 0)},dlevels={0: 0, -1: -1, -2: -2})
\end{lstlisting}
See \autoref{sec:AllGraphs} for details on profiles and graph generation.

The examples of the introduction also illustrated working in the tautological ring of \verb|X|. We may inspect the individual classes. Using \verb|print| gives a more readable output:
\begin{lstlisting}
sage: print(X.psi(1))
Tautological class on Stratum: (2,)
with residue conditions: []

1 * Psi class 1 with exponent 1 on level 0 * Graph ((), 0) +

sage: X.psi(1)
ELGTautClass(X=GeneralisedStratum(sig_list=[Signature((2,))],res_cond=[]),psi_list=[(1, AdditiveGenerator(X=GeneralisedStratum(sig_list=[Signature((2,))],res_cond=[]),enh_profile=((), 0),leg_dict={1: 1}))])
\end{lstlisting}
This illustrates how \verb|diffstrata| encodes elements of the tautological ring: 
a tautological class is represented by an \verb|ELGTautClass|, which is in turn a sum of \verb|AdditiveGenerator|s. Each \verb|AdditiveGenerator| corresponds to a $\psi$-monomial on a graph and thus carries the information of an enhanced profile and a \verb|leg_dict| encoding the $\psi$-powers:
every $\psi$-class is associated to a marked point of a level \cite[Thm. 1.5]{strataEC}, i.e. a leg of the graph. A $\psi$-monomial is thus encoded by a Python \verb|dict| with entries of the form \verb|l : n| where \verb|l| is the number of a leg of the graph and 
\verb|n| is the exponent of the $\psi$-class associated to this leg.
For example, we saw above that \verb|X.bics[0]| is the compact-type graph in the boundary of $\Omega\cM_2(2)$. We see from the \verb|LevelGraph| that the marked point is at leg \verb|2| (cf. \autoref{sec:LG}), the $\psi$-class at this point is therefore represented by the \verb|leg_dict| \verb|{2 : 1}|. We can enter this into \verb|diffstrata| as follows:
\begin{lstlisting}
sage: A = X.additive_generator(((0,), 0), {2 : 1})
sage: print(A)
Psi class 2 with exponent 1 on level 1 * Graph ((0,), 0)
\end{lstlisting}
Note that we had to use the \emph{enhanced} profile \verb|((0,), 0)| to refer to the graph.
Details and more examples may be found in \autoref{sec:AGs}.

Tautological classes may be added and multiplied. We can check that the
class~\verb|A| we defined agrees with the product of the $\psi$-class on
the stratum with the class of the graph:
\begin{lstlisting}
sage: A == X.psi(1) * X.additive_generator(((0,), 0))
True
\end{lstlisting}
Moreover, when squaring, e.g., the class of a graph, a normal bundle contribution appears:
\begin{lstlisting}
sage: print(X.additive_generator(((0,), 0))^2)
Tautological class on Stratum: (2,)
with residue conditions: []

-1 * Psi class 1 with exponent 1 on level 0 * Graph ((0,), 0) +
-1 * Psi class 3 with exponent 1 on level 1 * Graph ((0,), 0) +
\end{lstlisting}
The multiplication process is described in detail in \autoref{sec:mult}.

In the formulas for the Euler characteristic \cite{strataEC}, the class $\xi=c_1(\cO(-1))$ of the tautological bundle and its restriction $\xi_{B^{[i]}_\Gamma}$ to a level $i$ of a graph $\Gamma$ were key.
For a stratum, the class $\xi$ is easily accessible:
\begin{lstlisting}
sage: print(X.xi)
Tautological class on Stratum: (2,)
with residue conditions: []

3 * Psi class 1 with exponent 1 on level 0 * Graph ((), 0) +
-1 * Graph ((0,), 0) +
-1 * Graph ((1,), 0) +
\end{lstlisting}
Moreover, it is not difficult to compute $\xi_{B^{[i]}_\Gamma}$ (here for the top-level of the compact-type graph):
\begin{lstlisting}
sage: print(X.xi_at_level(0, ((0,),0)))
Tautological class on Stratum: (2,)
with residue conditions: []

1 * Psi class 1 with exponent 1 on level 0 * Graph ((0,), 0) +
\end{lstlisting}
More details and examples may be found in \autoref{topxieval} and \autoref{xiatlevel}.

We now describe these objects and the implementation in more detail.

%% file: sec_genstrata.tex
\section{Generalised Strata}
\label{sec:GenStrata}

We begin by briefly recalling the notions from \cite{LMS} and \cite{strataEC} in the generality that we here require.

\subsection{Strata with residue conditions}
\label{subsec:rescond}
To obtain a recursive structure on the boundary of $\LMS$, recall the definition of \emph{generalised stratum}, introduced in \cite[\S 4]{strataEC} to cover the case of a level of an enhanced level graph.
More precisely, we allow differentials on disconnected surfaces:
denote by $\mu_i=(m_{i,1},\dots,m_{i,n_i})\in\bZ^{n_i}$ the type of a differential,
i.e., we require that $\sum_{j=1}^{n_i} m_{i,j}=2g_i-2$ for some $g_i \in \bZ$
and $i=1,\ldots,k$. 
Then we define, for a tuple $\bfg = (g_1,\ldots,g_k)$
of genera and a tuple $\bfn = (n_1,\ldots,n_k)$ together with
$\bfmu = (\mu_1,\ldots, \mu_k)$, the disconnected stratum
\[\omoduli[\bfg,\bfn](\bfmu) \= \prod_{i=1}^k \omoduli[g_i,n_i](\mu_i)\,.\]
Note that the projectivized stratum $\bP\omoduli[\bfg,\bfn](\bfmu)$ is the quotient
by the diagonal action of $\bC^*$, not the quotient by the action of $(\bC^*)^k$.

Moreover, we consider subspaces of these \emph{cut out by residue conditions}.
More precisely, denote by $H_p \subseteq \cup_{i=1}^k \{(i, 1),\cdots (i, n_i)\}$ the subset of the marked points such that $m_{i,j} < -1$. 
Now consider vector spaces
$\frakR$ of the following special shape, modelled on the global residue condition
from \cite{BCGGM}:
for $\lambda$ a partition of $H_p$, with parts denoted
by $\lambda^{(k)}$, and a subset $\lambda_{\frakR}$ of the parts of $\lambda$, we define the $\bC$-vector space
\[\frakR\coloneqq \Bigl\{ r= (r_{i,j})_{(i,j) \in H_p} \in \bC^{H_p}
\quad \text{and} \sum_{(i,j) \in \lambda^{(k)}} r_{i,j}=0 \quad
\text{for all} \quad \lambda^{(k)} \in \lambda_{\frakR} \Bigr\}\,.\]
We denote the subspace of surfaces with residues in~$\frakR$ by
$\Romod[\bfmu][\bfg,\bfn][\frakR]$.

In \cite[Prop. 4.2]{strataEC} a modular compactification $\proj\LMS[\bfmu][\bfg,\bfn][\frakR]$ of $\proj \Romod[\bfmu][\bfg,\bfn][\frakR]$ is constructed in analogy to \cite{LMS}.
Consequently, the boundary components are parametrised by \emph{enhanced level graphs}.
More precisely, a \emph{level graph} is defined to be a stable graph together with a level function.
Recall that a \emph{stable graph} is a tuple $\Gamma = (V_i, H_i, E_i, g_i, v_i, \iota_i)_{i=1,\dotsc,k}$ consisting of \emph{vertices} $V_i$, a \emph{genus map} $g_i\colon V_i\to\bZ_{\geq 0}$,
\emph{legs} $H_i$ that are associated to the vertices by a \emph{vertex map} $v_i\colon H_i\to V_i$ and come with an involution $\iota_i\colon H_i\to H_i$, the two-cycles of which form the \emph{edges} $E_i\subseteq H_i\times H_i$ while the fixed points (denoted $H^m_i$) are in bijection with the $n_i$ marked points. 
The $k$ graphs $\Gamma_i=(V_i, H_i, E_i, g_i, v_i,\iota_i)$ are required to be connected and satisfy the usual stability conditions.
Moreover, we set $g = \sqcup g_i$, $v = \sqcup v_i$, $E = \sqcup E_i$, $ H=\sqcup H_i$, and $V=\sqcup V_i$.
Note that this data induces a unique bijection $o\colon\bigsqcup H^m_i\to\bfmu$ associating to each marked point $(i,j)$ the order $m_{i,j}$ of the differential.

A \emph{level function} on the vertices is a map
$\ell\colon V\to\bZ$, which we normalise to take values in $\{0,-1,\dotsc,-L\}$ and require that, for all edges $e\in E$, 
$\ell(v(e^+)) \geq \ell(v(e^-))$,
where we write $e\eqqcolon(e^+,e^-)\in H\times H$. 

Moreover, an \emph{enhancement} is a map $\kappa\colon E\to\bZ_{\geq 0}$ such that $\kappa(e)=0$ if and only if $e$ is horizontal (i.e. $\ell(v(e^+))=\ell(v(e^-))$), subject to the following stability condition:
define the \emph{degree} of a vertex~$v$ in $V$ to be
\[\deg(v) = \sum_{h \in H^m, v(h)=v} o(h) + \sum_{e \in E, v(e^+)=v}  (\kappa(e)-1)
- \sum_{e \in E, v(e^-)=v} (\kappa(e)+1) \,.\]
Then the enhancement is admissible, if
$\deg(v) = 2g(v)-2$ holds for every vertex~$v$ of~$\Gamma$.

The \emph{enhanced level graph} then consists of the triple $(\Gamma, \ell, \kappa)$. 
We denote the corresponding boundary component of $\proj\LMS[\bfmu][\bfg,\bfn][\frakR]$
by $D_\Gamma$.

\begin{rem}
\label{poleorders}
Note that $\kappa(e)$ corresponds to the number of \emph{prongs} at $e$ and, for a non-horizontal edge, the associated differential has a zero of order $\kappa(e)-1$ on the top component and a pole of order $-\kappa(e)-1$ on the bottom component (for horizontal edges there is a simple pole on each component). This gives an extension of $o$ to $H$.

Using this identification, the stability condition of \cite[\S 2]{LMS} is simply the requirement that the orders of zeros and poles sum to $2g_{i,j}-2$ on each vertex.
See \cite[\S 2]{LMS} and \cite[\S 3.2]{strataEC} for details.
\end{rem}

To determine, for a generalised stratum $\proj\LMS[\bfmu][\bfg,\bfn][\frakR]$, which enhanced level graphs give non-empty boundary components, we recall the $\frakR$-GRC from \cite[\S 4]{strataEC}:
starting with an enhanced level graph $\Gamma$,
we construct a new \emph{auxiliary level graph} $\widetilde{\Gamma}$
by adding, 
for each $\lambda^{(k)}\in \lambda_\frakR$, 
a new vertex $v_{\lambda^{(k)}}$ to $\Gamma$ at level $\infty$  and converting a tuple $(i,j)\in \lambda^{(k)}$
into an edge from the marked point $(i,j)$ to  the vertex $v_{\lambda^{(k)}}$. 
We then say that $\Gamma$ satisfies the
\emph{$\frakR$-global residue condition ($\frakR$-GRC)} if 
the tuple of residues at the legs in $H_p$ belongs to~$\frakR$ and for every
level $L<\infty$ of $\widetilde{\Gamma}$ and every connected component~$Y$
of the subgraph~$\widetilde{\Gamma}_{>L}$ one of
the following conditions holds.
\begin{enumerate}
\item The component~$Y$ contains a marked point with a prescribed
pole that is \emph{not} in $\lambda_\frakR$.
\item The component~$Y$ contains a marked point with a prescribed
pole $(i,j) \in H_p$ and
there is an $r \in \frakR$ with $r_{(i,j)} \neq 0$.
\item Let $e_1,\ldots,e_b$ denote the set of edges where~$Y$ 
intersects $\widetilde{\Gamma}_{=L}$. Then
\[ \sum_{j=1}^b\Res_{e_j^-}\eta_{v(e^-_j)}\=0\,,\]
where  $v(e^-_j)\in\widetilde{\Gamma}_{=L}$.
\end{enumerate}
By \cite[Prop. 4.2]{strataEC}, the boundary components $D_\Gamma$ of $\proj\LMS[\bfmu][\bfg,\bfn][\frakR]$ are parametri\-sed by enhanced level graphs $(\Gamma, \ell, \kappa)$ satisfying the $\frakR$-GRC.

For applications such as listing all graphs in the boundary of a stratum, it
is convenient to have a purely graph-theoretic criterion in analogy to the one shown for the classical GRC in \cite{MUW}.
The following
proposition strips the tropical language away in the criterion 
\cite[Theorem~1]{MUW} and generalises it to meromorphic strata
with residue conditions.

Let~$(\Gamma,\ell,\kappa)$ be an enhanced level graph. 

We call a vertex~$v$ of~$\Gamma$ \emph{inconvenient} if $g(v) = 0$,
if it is not adjacent to any edge~$e$ with enhancement $\kappa(e) = 0$
and if it is adjacent to a leg with a \emph{very high enhancement} in
the following precise sense: denote by $\frakp(v)$ the set of half-edges on the vertex $v$ that are poles (in the sense of \autoref{poleorders}).
Then the condition is that there is a $p\in\frakp(v)$ such that
\[o(p) > \sum_{p'\in\frakp(v)} (o(p') - 1) - 1.\]
\begin{prop} \label{prop:realizewithR}
The boundary stratum $D_\Gamma$ of $\Romod[\bfmu][\bfg,\bfn][\frakR]$
associated with the enhanced level graph~$(\Gamma,\ell,\kappa)$ is non-empty if and only if
both of the following conditions are satisfied
\begin{enumerate}
\item For every inconvenient vertex~$v$ of~$\Gamma$ there is
\begin{enumerate}
\item a simple cycle based at~$v$ that does not pass through
any vertex of level smaller than~$v$, or
\item the graph of levels~$\geq \ell(v)$ deprived of the vertex~$v$ has two
components, each of which has a marked pole in $H_p$ whose residue is not constrained to zero for all elements of~$\frakR$.
\end{enumerate}
\item For every horizontal edge~$e$ of~$\Gamma$ there is
\begin{enumerate}
\item a simple cycle based through~$e$ that does not pass through
any vertex of level smaller than~$\ell(e)$, or
\item the graph of levels~$\geq \ell(e)$ deprived of the edge~$e$ has two
components, each of which has a marked pole in $H_p$ whose residue is not constrained
to zero for all elements of~$\frakR$.
\end{enumerate}
\end{enumerate}
\end{prop}
\par
\begin{proof}
Consider first the case that $\widetilde{\Gamma}$ is connected.
We view $\widetilde{\Gamma}$ as an enhanced level graph of an auxiliary stratum $\widetilde{X}$ as follows: we add prongs to the new edges of $\widetilde{\Gamma}$ in accordance with the pole orders of the half-edges on $\Gamma$.
For each vertex $v$ at level $\infty$, we then extend the genus function by $g_v$ setting $2g_v-2$ as the sum of orders of the half-edges on $v$ induced by the newly added prongs, possibly adding an extra simple zero to fix parity issues. Note that all new vertices are of positive genus, so stability is not an issue.
Therefore, for the stratum $\widetilde{X}$, we are now reduced to the situation of \cite[Theorem~1]{MUW}.

In a product of strata, clearly a graph is admissible if and only if each component is admissible.
\end{proof}
See \autoref{subsec:R-GRC} for examples illustrating this criterion.

This criterion allows us to explicitly construct all graphs in a given stratum.
In fact, we can construct all graphs with no horizontal edges recursively from the two-level graphs.

\subsection{Constructing Level Graphs}
Recall the undegeneration maps $\delta_i$ \cite[\S 3.3]{strataEC}, contracting
all level crossing of an enhanced level graph $\Gamma$ without horizontal edges
except for the $i$-th level crossing, yielding a two-level graph.
The component of $\Gamma$ is contained in the product of the components of $\delta_i(\Gamma)$, which may be irreducible.

\begin{definition}\label{mathprofile}
Let $(\Gamma,\ell,\kappa)$ be an enhanced level graph with $L$ levels and without horizontal edges. We define the \emph{profile} of $\Gamma$ to be the tuple $(\delta_1(\Gamma),\dotsc,\delta_L(\Gamma))$.
An \emph{enhanced profile} is a profile together with a choice of irreducible component.
\end{definition}
Note that the association of a profile to a graph is not injective, see \autoref{reducibleprofile}.  
It is primarily useful to encode efficiently how graphs degenerate.
By definition, the information of an enhanced profile of $\Gamma$ is equivalent to the data of $\Gamma$.

Generating all non-horizontal graphs inside a stratum is thus equivalent to listing all non-empty enhanced profiles.
We do this recursively. For $X$ a generalised stratum, denote by $\BIC(X)$ the (non-horizontal) two-level graphs in the boundary of $X$.
\begin{definition}
Let $X$ be a generalised stratum. For each $\Gamma\in\BIC(X)$ we define:
\begin{enumerate}
\item the generalised strata $\Gamma^\top$ and $\Gamma^\bot$, the top and bottom levels of $\Gamma$;
\item a map $\beta_\Gamma^\top\colon\BIC(\Gamma^\top)\to\BIC(X)$ that associates to a graph $\Gamma'\in\BIC(\Gamma^\top)$ the graph $\delta_0(\Lambda)\in\BIC(X)$ where $\Lambda$ is the graph obtained by clutching $\Gamma'$ to $\Gamma^\bot$;
\item a map $\beta_\Gamma^\bot\colon\BIC(\Gamma^\bot)\to\BIC(X)$ that associates to a graph $\Gamma'\in\BIC(\Gamma^\bot)$ the graph $\delta_1(\Lambda)\in\BIC(X)$ where $\Lambda$ is the graph obtained by clutching $\Gamma^\top$ to $\Gamma'$.
\end{enumerate}
\end{definition}

The following proposition is an immediate consequence of \cite[Prop. 5.1]{strataEC}.

\begin{prop}\label{prop:beta}
Let $X$ be a generalised stratum and $\Gamma\in\BIC(X)$.
\begin{enumerate}
\item The images of $\beta^\top_\Gamma$ and $\beta^\bot_\Gamma$ are disjoint.
\item The profile $(\Gamma',\Gamma)$ is non-empty if and only if $\Gamma'$ is in the image of $\beta^\top_\Gamma$.
\item The profile $(\Gamma,\Gamma')$ is non-empty if and only if $\Gamma'$ is in the image of $\beta^\bot_\Gamma$.
\end{enumerate}
\end{prop}

As a consequence, we may define a partial order $\prec$ on $\BIC(X)$ by defining $\Gamma\prec\Gamma'$ if and only if $(\Gamma',\Gamma)$ is non-empty.

\begin{rem}
The maps $\beta$ are not necessarily injective. 
Indeed, whenever a profile $(\Gamma_1,\Gamma_2)$ is reducible, i.e. contains at least two distinct enhanced level graphs $\Lambda_1\neq\Lambda_2$, cutting the bottom level off these three-level graphs gives two distinct BICs $\tilde{\Lambda}_1, \tilde{\Lambda}_2\in\BIC(\Gamma_2^\top)$ (the ``cut'' edges correspond to edges of $\Gamma_2$ and become marked points in $\Gamma_2^\top$ and $\Gamma_2^\bot$) with $\beta_{\Gamma_2}^\top(\tilde{\Lambda}_1)=\beta_{\Gamma_2}^\top(\tilde{\Lambda}_2)=\Gamma_1$. See also \autoref{reducible_degen}.

However, non-injectivity does not imply reducibility. Degenerating a level in different ways may give \emph{isomorphic} graphs, see \autoref{reducible_aut} and \autoref{fig:vgraph}.
\end{rem}

This allows us to recursively compute all profiles.

\begin{prop}\label{prop:gengraphs}
Let $X$ be a generalised stratum.
\begin{enumerate}
\item The enhanced level graphs in $\BIC(X)$ can be listed explicitly.
\item All non-empty profiles in $X$ can be constructed recursively from $\BIC(X)$.
\item All graphs inside a profile can be constructed explicitly from the profile's components.
\end{enumerate}
\end{prop}
\begin{proof}
This is essentially the content of \autoref{sec:AllGraphs}: 
the effective construction of BICs is explained in detail in \autoref{subsec:bicgeneration}, in particular \autoref{BICgenalgo}.

The non-empty profiles are constructed recursively: a non-empty profile $(\Gamma_1,\dotsc,\Gamma_l)$ of length $l$ may be extended to a non-empty profile $(\Gamma_0,\Gamma_1,\dotsc,\Gamma_l)$ if and only if $\Gamma_1\prec\Gamma_0$.
Indeed, $\Gamma_1^\top$ is also the top level of any graph $\Lambda$ in $(\Gamma_1,\dotsc,\Gamma_l)$ and thus a preimage $(\beta_{\Gamma_1}^\top)^{-1}(\Gamma_0)$ exists in $\BIC(\Gamma_1^\top)$ and can be clutched to $\Lambda$ to yield a graph in $(\Gamma_0,\Gamma_1,\dotsc,\Gamma_l)$.
Similarly, the profile may be extended to $(\Gamma_1,\dotsc,\Gamma_l,\Gamma_{l+1})$ if and only if $\Gamma_{l+1}\prec\Gamma_l$.

The converse direction is simply squishing of a level.

To get all the graphs we follow the above procedure, noting that we might obtain several graphs in the same profile if $\beta$ is non-injective (we clutch each preimage with each other graph).
\end{proof}

The key observation is that BICs and three-level graphs are sufficient for constructing the entire stratum. In particular, all levels are seen by these.

\begin{rem}\label{leveltypes}
Let $X$ be a generalised stratum. Note that any level $L$ appearing in any graph in $X$ is one of the following three types:
\begin{enumerate}
\item a top level of a BIC,
\item a bottom level of a BIC, or
\item a middle level of a three-level graph.
\end{enumerate}
Indeed, given a graph $\Gamma$, level $l$ of $\Gamma$ remains unchanged by contracting any level crossing not adjacent to $l$. Contracting all non-adjacent levels results either in a BIC (if $l$ is top or bottom level) or in a three-level graph around level $l$.
\end{rem}

While the reducibility is recorded by the three-level graphs, determining the componentes of a profile is not straight-forward.
Indeed, given two non-empty profiles $(\Gamma_1,\Gamma_2)$ and $(\Gamma_2,\Gamma_3)$ the profile $(\Gamma_1,\Gamma_2,\Gamma_3)$ is non-empty, but in the reducible case it is not clear how the reducibility and the individual graphs are related.

\begin{rem}
Determining the reducibility of a profile is a delicate issue and is discussed in detail in \autoref{subsec:degenerations}.
In particular, each of the following may occur:
\begin{enumerate}
\item a degeneration of an irreducible profile can be reducible (\autoref{reducibleprofile});
\item a degeneration of a reducible profile can be irreducible (\autoref{reducible_degen});%
\item A reducible profile implies a non-injectivity of one of the maps $\beta$.
\end{enumerate}
The converse of the last statement is false in general, as the non-injectivity may stem from the presence of automorphisms.
\end{rem}

%% file: sec_levelgraphs.tex
\section{Level Graphs and Embeddings}
\label{sec:LG}

The \verb|diffstrata| package uses three basic classes to model the objects appearing in the boundary of $\bP\LMS$.

\begin{itemize}
\item \verb|GeneralisedStratum| models the strata $\bP\LMS$. All essential operations will be performed by an object of this type.
\item \verb|LevelGraph| is a low-level object that represents the actual underlying graph. While it is important for the underlying calculations, there should be little or no reason to construct an explicit \verb|LevelGraph| directly. It was originally modelled on the class \verb|stgraph| of \verb|admcycles|.
\item \verb|EmbeddedLevelGraph| is essentially a wrapper to encode how a \verb|LevelGraph| is embedded into a \verb|GeneralisedStratum|. A \verb|LevelGraph| will usually appear as an \verb|EmbeddedLevelGraph| and should be accessed through its \emph{profile}, see \autoref{sec:AllGraphs}.
\end{itemize}

Note that every level of an \verb|EmbeddedLevelGraph| is itself a \verb|GeneralisedStratum| of lower dimension, hence the recursive structure.

A \verb|GeneralisedStratum| is given by the following data:
\begin{itemize}
\item a list of \verb|Signature| objects, giving the signatures of the (possibly) meromorphic strata appearing as factors, and
\item optionally a list of residue conditions, signalling which residues add up to zero.
\end{itemize}

Consequently, the marked points of a stratum may be uniquely referred to by their coordinates inside the signature tuple. 

Note that there is also the \verb|Stratum| class, which is simply a frontend for the class \verb|GeneralisedStratum| with a simpler syntax: It can only be used to create connected strata with no residue conditions:
\begin{lstlisting}
sage: X=Stratum((2,))
sage: print(X)
Stratum: (2,)
with residue conditions: []
sage: isinstance(X, GeneralisedStratum)
True
\end{lstlisting}
Note that \verb|tuple|s with only one entry must be terminated by a \verb|,|!

\begin{rem}
A \verb|GeneralisedStratum| must be provided with a list of \verb|Signature| objects, a \verb|tuple| will raise an error! A \verb|Signature| object is initialised by a signature \verb|tuple| and makes all its intrinsic properties easily accessible. It is easiest understood via the package documentation:
\begin{lstlisting}
sage: Signature?
Init signature: Signature(sig)
Docstring:
   A signature of a stratum.

    Attributes:
        sig (tuple): signature tuple
        g (int): genus
        n (int): total number of points
        p (int): number of poles
        z (int): number of zeroes
        poles (tuple): tuple of pole orders
        zeroes (tuple): tuple of zero orders
        pole_ind (tuple): tuple of indices of poles
        zero_ind (tuple): tuple of indices of zeroes
        
   EXAMPLES

      sage: from admcycles.diffstrata.sig import Signature
      sage: sig=Signature((2,1,-1,0))
      sage: sig.g
      2
      sage: sig.n
      4
      sage: sig.poles
      (-1,)
      sage: sig.zeroes
      (2, 1)
      sage: sig.pole_ind
      (2,)
      sage: sig.zero_ind
      (0, 1)
      sage: sig.p
      1
      sage: sig.z
      2
Init docstring:
   Initialise signature

   Args:
      sig (tuple): signature tuple of integers adding up to 2g-2
\end{lstlisting}
\end{rem}

\subsection{Points on graphs and strata}
\label{pointtypes}
Analogous to \autoref{sec:GenStrata}, the \verb|diffstrata| package uses zeros and poles of a differential in two different contexts:
the \emph{marked points} of the stratum, as elements of $\bfmu$, have associated \emph{half-edges} of the graph.
To illustrate these, we need to first briefly explain how a \verb|LevelGraph| encodes the information of an enhanced level graph.

A point or \emph{leg} of a \verb|LevelGraph| is given by a positive integer. Each vertex of a \verb|LevelGraph| has a (possibly empty) list of legs associated to it and each edge is a \verb|tuple| of two legs: the command
\begin{lstlisting}
sage: L=LevelGraph([1, 0],[[1, 2], [3, 4, 5]],[(1, 4), (2, 5)],{1: 0, 2: 0, 3: 2, 4: -2, 5: -2},[0, -1])
\end{lstlisting}
encodes a graph $L$ with two vertices, one of genus $1$ (with internal name, the position in the list, \verb|0|) and one of genus $0$ (with internal name \verb|1|). The (nested) list of legs is in the same order as the list of vertices: the vertex \verb|0| has two legs, stored in the list \verb|[1, 2]| and the vertex \verb|1| has three legs, stored in \verb|[3, 4, 5]|. 

Furthermore, the graph has two edges, the first, \verb|(1, 4)|, connecting leg \verb|1| on vertex~\verb|0| with leg~\verb|4| on vertex \verb|1|, and the second, \verb|(2, 5)|, connecting leg \verb|2| on vertex~\verb|0| to leg~\verb|5| on vertex \verb|1|.

There are two more pieces of information needed to determine a \verb|LevelGraph|: a \verb|dict|, associating to every leg the order of the differential at this point (e.g. leg \verb|4| is a pole of order $-2$ and leg \verb|3| is a zero of order $2$; leg \verb|1| is simply a marked point (of order $0$)) and a \verb|list| of levels, \verb|[0, -1]|, indicating that vertex \verb|0| is on level $0$ and vertex \verb|1| is on level $-1$.

\begin{figure}
\[\begin{tikzpicture}[
		baseline={([yshift=-.5ex]current bounding box.center)},
		scale=2,very thick,
		bend angle=30,
		every loop/.style={very thick},
     		comp/.style={circle,fill,black,,inner sep=0pt,minimum size=5pt},
		order bottom left/.style={pos=.05,left,font=\tiny},
		order top left/.style={pos=.9,left,font=\tiny},
		order bottom right/.style={pos=.05,right,font=\tiny},
		order top right/.style={pos=.9,right,font=\tiny},
		circled number/.style={circle, draw, inner sep=0pt, minimum size=12pt},
		bottom right with distance/.style={below right,text height=10pt}]
\node[circled number] (T) [] {$1$}; 
\node[comp] (B) [below=of T] {}
	edge [bend left] 
		node [order bottom left] {$-2$} 
		node [order top left] {$0$} (T)
	edge [bend right] 
		node [order bottom right] {$-2$} 
		node [order top right] {$0$} (T);
\node [bottom right with distance] (B-2) at (B.south east) {$2$};
\path (B) edge [shorten >=4pt] (B-2.center);
\end{tikzpicture}\]
\caption{The ``banana'' graph in the boundary of $\Omega\cM_{2}(2)$.}
\label{fig:banana}
\end{figure}
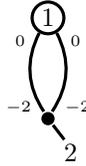

This uniquely determines an enhanced level graph in the sense of \autoref{sec:GenStrata}.
We see that $L$ describes the ``banana'' graph in the boundary of $\Omega\cM_{2}(2)$ depicted in \autoref{fig:banana}.

We can access all of this information about the \verb|LevelGraph| from within \verb|sage|:
\begin{lstlisting}
sage: L=LevelGraph([1, 0],[[1, 2], [3, 4, 5]],[(1, 4), (2, 5)],{1: 0, 2: 0, 3: 2, 4: -2, 5: -2},[0, -1])
sage: L.genus(0)
1
sage: L.legsatvertex(1)
[3, 4, 5]
sage: L.edges
[(1, 4), (2, 5)]
sage: L.vertex(1)  # vertex of leg
0
sage: L.levelofvertex(1)
-1
sage: L.orderatleg(3)
2
\end{lstlisting}

The orders are equivalent to the number of prong-matchings at the edges. The prongs are stored in a dictionary and their lcm can be calculated easily:
\begin{lstlisting}
sage: L=LevelGraph([0, 0, 0],[[1, 2, 3], [4, 5, 6], [7, 8, 9]],[(3, 6), (4, 8), (5, 9)],{1: 1, 2: -4, 3: 1, 4: 1, 5: 0, 6: -3, 7: 3, 8: -3, 9: -2},[0, -1, -2])
sage: L.prongs.items()
dict_items([((3, 6), 2), ((4, 8), 2), ((5, 9), 1)])
sage: L.prongs[(3,6)]
2
sage: lcm(L.prongs.values())
2
sage:
\end{lstlisting}

\begin{rem}
Note that, as mentioned above, \verb|LevelGraph|s should not be entered ``by hand'', but instead \verb|EmbeddedLevelGraph|s and profiles should be used.
\end{rem}

A marked point as an element of $\bfmu$ is considered by \verb|diffstrata| as a point of the \verb|GeneralisedStratum|.
An \verb|EmbeddedLevelGraph| is, essentially, a \verb|LevelGraph| together with the information, which of its legs correspond to marked points of the stratum.

This is encoded by a \verb|dict| usually denoted \verb|dmp| (dictionary of marked points) that identifies legs of the \verb|LevelGraph| (integers) with points of the stratum (\verb|tuple|s: index of component, index of point in signature, as in \autoref{sec:GenStrata}).

For technical reasons, the embedding also requires a dictionary of levels, \verb|dlevels|.

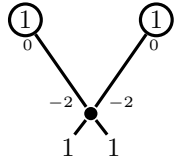
\begin{figure}
\[\begin{tikzpicture}[
		baseline={([yshift=-.5ex]current bounding box.center)},
		scale=2,very thick,
		bend angle=30,
		every loop/.style={very thick},
     		comp/.style={circle,fill,black,,inner sep=0pt,minimum size=5pt},
		order bottom left/.style={pos=.05,left,font=\tiny},
		order top left/.style={pos=.85,left,font=\tiny},
		order bottom right/.style={pos=.05,right,font=\tiny},
		order top right/.style={pos=.85,right,font=\tiny},
		circled number/.style={circle, draw, inner sep=0pt, minimum size=12pt},
		bottom right with distance/.style={below right,text height=10pt},
		bottom left with distance/.style={below left,text height=10pt}]
\begin{scope}[node distance=.5cm]
\node[circled number] (L) [] {$1$}; 
\node [] (T) [right=of L] {};
\node[circled number] (R) [right=of T] {$1$};
\end{scope}
\node[comp] (B) [below=of T] {}
	edge
		node [order bottom left] {$-2$} 
		node [order top left] {$0$} (L)
	edge
		node [order bottom right] {$-2$} 
		node [order top right] {$0$} (R);
\node [bottom right with distance] (B-11) at (B.south east) {$1$};
\path (B) edge [shorten >=4pt] (B-11.center);
\node [bottom left with distance] (B-12) at (B.south west) {$1$};
\path (B) edge [shorten >=4pt] (B-12.center);
\end{tikzpicture}\]
\caption{The ``V''-shaped graph in the boundary of the stratum $\Omega\cM_{2}(1,1)$.}
\label{fig:firstvgraph}
\end{figure}

In the above example, giving such an embedding is straight-forward:
\begin{lstlisting}
sage: L=LevelGraph([1, 0],[[1, 2], [3, 4, 5]],[(1, 4), (2, 5)],{1: 0, 2: 0, 3: 2, 4: -2, 5: -2},[0, -1])
sage: X=Stratum((2,))
sage: ELG=EmbeddedLevelGraph(X, L, dmp={3: (0,0)}, dlevels={0: 0, -1: -1})
\end{lstlisting}
Of course, our graph is isomorphic to one of the two BICs generated automatically by the stratum $X$:
\begin{lstlisting}
sage: any(ELG.is_isomorphic(B) for B in X.bics)
True
sage: len([B for B in X.bics if ELG.is_isomorphic(B)])
1
\end{lstlisting}

Our \verb|GeneralisedStratum| objects may have several connected components, as well as residue conditions, since we want to consider arbitrary levels of \verb|LevelGraph|s as strata,

Consider the ``V''-shaped graph in the boundary of the stratum $\Omega\cM_{2}(1,1)$ depicted in \autoref{fig:firstvgraph}.
In \verb|sage|, we may extract the two levels to see examples of more complicated strata:
\begin{lstlisting}
sage: X=Stratum((1,1))
sage: V=LevelGraph([1, 1, 0],[[1], [2], [3, 4, 5, 6]],[(1, 5), (2, 6)],{1: 0, 2: 0, 3: 1, 4: 1, 5: -2, 6: -2},[0, 0, -1])
sage: ELV=EmbeddedLevelGraph(X, V, dmp={3: (0, 0), 4: (0, 1)}, dlevels={0: 0, -1: -1})
sage: ELV.level(0)
LevelStratum(sig_list=[Signature((0,)), Signature((0,))],res_cond=[],leg_dict={1: (0, 0), 2: (1, 0)})
sage: ELV.level(1)
LevelStratum(sig_list=[Signature((1, 1, -2, -2))],res_cond=[[(0, 2)], [(0, 3)]],leg_dict={3: (0, 0), 4: (0, 1), 5: (0, 2), 6: (0, 3)})
\end{lstlisting}
Using \verb|print| gives more readable output:
\begin{lstlisting}
sage: print(ELV.level(0))
Product of Strata:
Signature((0,))
Signature((0,))
with residue conditions: 
dimension: 3
leg dictionary: {1: (0, 0), 2: (1, 0)}
leg orbits: [[(1, 0), (0, 0)]]

sage: print(ELV.level(1))
Stratum: Signature((1, 1, -2, -2))
with residue conditions: [(0, 2)] [(0, 3)]
dimension: 0
leg dictionary: {3: (0, 0), 4: (0, 1), 5: (0, 2), 6: (0, 3)}
leg orbits: [[(0, 0)], [(0, 1)], [(0, 3), (0, 2)]]
\end{lstlisting}
Observe that the extracted levels also remember the action of the automorphism group of the graph they were extracted from on their marked points.

The residue conditions are given by a nested list of points of the stratum. Marked poles whose residue adds up to $0$ are contained in the same list.

To illustrate this, let us compare the lower level of the Banana graph of above:
\begin{lstlisting}
sage: B=LevelGraph([1, 0],[[1, 2], [3, 4, 5]],[(1, 4), (2, 5)],{1: 0, 2: 0, 3: 2, 4: -2, 5: -2},[0, -1])
sage: Y=Stratum((2,))
sage: ELB=EmbeddedLevelGraph(Y, B, dmp={3: (0,0)}, dlevels={0: 0, -1: -1})
sage: ELB.level(1)
LevelStratum(sig_list=[Signature((2, -2, -2))],res_cond=[[(0, 1), (0, 2)]],leg_dict={3: (0, 0), 4: (0, 1), 5: (0, 2)})
sage: print(ELB.level(1))
Stratum: Signature((2, -2, -2))
with residue conditions: [(0, 1), (0, 2)]
dimension: 0
leg dictionary: {3: (0, 0), 4: (0, 1), 5: (0, 2)}
leg orbits: [[(0, 0)], [(0, 1), (0, 2)]]
\end{lstlisting}
In this case, the stratum records that the poles share a residue condition.

\begin{rem}\label{levelnumbers}
Note that while, mathematically, levels are usually indexed with negative numbers starting at $0$, in \verb|diffstrata| it is often much less confusing to work with positive level numbers. To ease this translation, \verb|LevelGraph|s come with a notion of ``internal'' versus ``relative'' level number:
\begin{lstlisting}
sage: B.internal_level_number(1)
-1
sage: B.level_number(-1)
1
\end{lstlisting}
In particular, the internal level numbers might not even be consecutive, while the \verb|level_number|s are guaranteed to run from $0$ to the number of levels.
\end{rem}

\verb|EmbeddedLevelGraph|s also come with an \verb|explain| method, intended to describe them in a more human-readable format. For the above examples:
\begin{lstlisting}
sage: ELV.explain()
LevelGraph embedded into stratum Stratum: (1, 1)
with residue conditions: []
 with:
On level 0:
* A vertex (number 0) of genus 1
* A vertex (number 1) of genus 1
On level 1:
* A vertex (number 2) of genus 0
The marked points are on level 1.
More precisely, we have:
* Marked point (0, 0) of order 1 on vertex 2 on level 1
* Marked point (0, 1) of order 1 on vertex 2 on level 1
Finally, we have 2 edges. More precisely:
* one edge between vertex 0 (on level 0) and vertex 2 (on level 1) with prong 1.
* one edge between vertex 1 (on level 0) and vertex 2 (on level 1) with prong 1.
sage: ELB.explain()
LevelGraph embedded into stratum Stratum: (2,)
with residue conditions: []
 with:
On level 0:
* A vertex (number 0) of genus 1
On level 1:
* A vertex (number 1) of genus 0
The marked points are on level 1.
More precisely, we have:
* Marked point (0, 0) of order 2 on vertex 1 on level 1
Finally, we have 2 edges. More precisely:
* 2 edges between vertex 0 (on level 0) and vertex 1 (on level 1) with prongs 1 and 1.
\end{lstlisting}

\subsection{Checking the $\mathfrak{R}$-GRC: points ``at level $\infty$''}
\label{subsec:R-GRC}

Checking the $\mathfrak{R}$-GRC is split into two steps. First, for \emph{any} \verb|LevelGraph|, we may check the classical GRC:

\begin{lstlisting}
sage: L=LevelGraph([1, 0],[[1, 2], [3, 4, 5]],[(1, 4), (2, 5)],{1: 0, 2: 0, 3: 2, 4: -2, 5: -2},[0, -1])
sage: L.is_legal()
True
\end{lstlisting}

\begin{example}\label{ex:GRC}
Consider the graphs in the stratum $(2,1,1)$ that are depicted in \autoref{fig:longzigzag}.
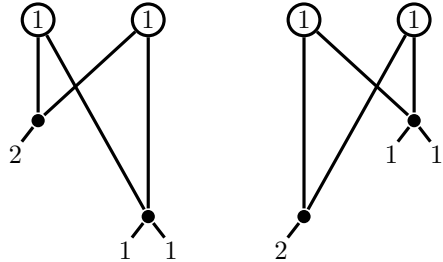
\begin{figure}
\[\begin{tikzpicture}[
		baseline={([yshift=-.5ex]current bounding box.center)},
		scale=2,very thick,
		bend angle=30,
		every loop/.style={very thick},
     		comp/.style={circle,fill,black,,inner sep=0pt,minimum size=5pt},
		prong left/.style={pos=.5,left,font=\small},
		prong right/.style={pos=.5,right,font=\small},
		circled number/.style={circle, draw, inner sep=0pt, minimum size=12pt},
		bottom right with distance/.style={below right,text height=10pt},
		bottom left with distance/.style={below left,text height=10pt}]
\node[circled number] (TL) {$1$};
\node[circled number] (TR) [right=of TL] {$1$};
\node[comp] (BL) [below=of TL] {}
	edge (TL)
	edge (TR);
\node (BR1) [below=of TR] {};
\node[comp] (BR2) [below=of BR1] {}
	edge (TL)
	edge (TR);
\node [bottom left with distance] (BL-2) at (BL.south west) {$2$};
\path (BL) edge [shorten >=4pt] (BL-2.center);
\node [bottom left with distance] (BR-11) at (BR2.south west) {$1$};
\path (BR2) edge [shorten >=4pt] (BR-11.center);
\node [bottom right with distance] (BR-12) at (BR2.south east) {$1$};
\path (BR2) edge [shorten >=4pt] (BR-12.center);
\end{tikzpicture}
\quad \quad \quad 
\begin{tikzpicture}[
		baseline={([yshift=-.5ex]current bounding box.center)},
		scale=2,very thick,
		bend angle=30,
		every loop/.style={very thick},
     		comp/.style={circle,fill,black,,inner sep=0pt,minimum size=5pt},
		prong left/.style={pos=.5,left,font=\small},
		prong right/.style={pos=.5,right,font=\small},
		circled number/.style={circle, draw, inner sep=0pt, minimum size=12pt},
		bottom right with distance/.style={below right,text height=10pt},
		bottom left with distance/.style={below left,text height=10pt}]
\node[circled number] (TL) {$1$};
\node[circled number] (TR) [right=of TL] {$1$};
\node (BL1) [below=of TL] {};
\node[comp] (BL2) [below=of BL1] {}
	edge (TL)
	edge (TR);
\node[comp] (BR) [below=of TR] {}
	edge (TL)
	edge (TR);
\node [bottom left with distance] (BL2-2) at (BL2.south west) {$2$};
\path (BL2) edge [shorten >=4pt] (BL2-2.center);
\node [bottom left with distance] (BR-11) at (BR.south west) {$1$};
\path (BR) edge [shorten >=4pt] (BR-11.center);
\node [bottom right with distance] (BR-12) at (BR.south east) {$1$};
\path (BR) edge [shorten >=4pt] (BR-12.center);
\end{tikzpicture}\]
\caption{An illegal (left) and legal (right) graph in the boundary of $\Omega\cM_3(2,1,1)$.}
\label{fig:longzigzag}
\end{figure}
Note that the left graph is illegal, as the bottom-left vertex is inconvenient as defined in \autoref{subsec:rescond}: It is a stratum with signature $(2,-2,-2)$ and both residues are forced zero, as there is no cycle or pole in the graph above to rectify this. This problem does not occur in the right graph.

In \verb|diffstrata|, we observe:
\begin{lstlisting}
sage: L=LevelGraph([1, 1, 0, 0],[[1, 2], [3, 4], [5, 6, 7], [8, 9, 10, 11]],[(1, 6), (3, 7), (4, 10), (2, 11)],{1: 0, 2: 0, 3: 0, 4: 0, 5: 2, 6: -2, 7: -2, 8: 1, 9: 1, 10: -2, 11: -2},[0, 0, -1, -2])
sage: L.is_legal()
Vertex 2 is illegal!
False
sage: R=LevelGraph([1, 1, 0, 0],[[1, 2], [3, 4], [5, 6, 7], [8, 9, 10, 11]],[(1, 6), (3, 7), (4, 10), (2, 11)],{1: 0, 2: 0, 3: 0, 4: 0, 5: 2, 6: -2, 7: -2, 8: 1, 9: 1, 10: -2, 11: -2},[0, 0, -2, -1])
sage: R.is_legal()
True
\end{lstlisting}
\end{example}

Similarly, if we start with an \verb|EmbeddedLevelGraph|, we may check the $\frakR$-GRC.

\begin{example}\label{zigzaggraph}
Consider the example from above as a degeneration of the bottom level of the ``zigzag-graph'', the BIC that is the common $\delta_1$ of the two graphs of \autoref{ex:GRC}, as depicted in \autoref{fig:zigzag}.
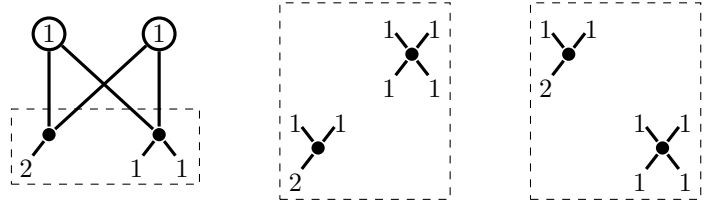
\begin{figure}
\[\begin{tikzpicture}[
		baseline={([yshift=-.5ex]current bounding box.center)},
		scale=2,very thick,
		bend angle=30,
		every loop/.style={very thick},
     		comp/.style={circle,fill,black,,inner sep=0pt,minimum size=5pt},
		prong left/.style={pos=.5,left,font=\small},
		prong right/.style={pos=.5,right,font=\small},
		circled number/.style={circle, draw, inner sep=0pt, minimum size=12pt},
		bottom right with distance/.style={below right,text height=10pt},
		bottom left with distance/.style={below left,text height=10pt}]
\node[circled number] (TL) {$1$};
\node[circled number] (TR) [right=of TL] {$1$};
\node[comp] (BL) [below=of TL] {}
	edge (TL)
	edge (TR);
\node[comp] (BR) [below=of TR] {}
	edge (TL)
	edge (TR);
\node [bottom left with distance] (BL-2) at (BL.south west) {$2$};
\path (BL) edge [shorten >=4pt] (BL-2.center);
\node [bottom left with distance] (BR-11) at (BR.south west) {$1$};
\path (BR) edge [shorten >=4pt] (BR-11.center);
\node [bottom right with distance] (BR-12) at (BR.south east) {$1$};
\path (BR) edge [shorten >=4pt] (BR-12.center);
\draw[thin, dashed] (-.25,-.5) rectangle (1,-1);
\end{tikzpicture}
\quad \quad \quad 
\begin{tikzpicture}[
		baseline={([yshift=-.5ex]current bounding box.center)},
		scale=2,very thick,
		bend angle=30,
		every loop/.style={very thick},
     		comp/.style={circle,fill,black,,inner sep=0pt,minimum size=5pt},
		prong left/.style={pos=.5,left,font=\small},
		prong right/.style={pos=.5,right,font=\small},
		circled number/.style={circle, draw, inner sep=0pt, minimum size=12pt},
		bottom right with distance/.style={below right,text height=10pt},
		bottom left with distance/.style={below left,text height=10pt},
		above right with distance/.style={above right,text height=10pt},
		above left with distance/.style={above left,text height=10pt}]
\node (BL1) {};
\node[comp] (BL2) [below=of BL1] {};
\node[comp] (BR) [right=of BL1] {};
\node [bottom left with distance] (BL2-2) at (BL2.south west) {$2$};
\path (BL2) edge [shorten >=4pt] (BL2-2.center);
\node [above left with distance] (BL2-11) at (BL2.north west) {$1$};
\path (BL2) edge [shorten >=4pt] (BL2-11.center);
\node [above right with distance] (BL2-12) at (BL2.north east) {$1$};
\path (BL2) edge [shorten >=4pt] (BL2-12.center);
\node [bottom left with distance] (BR-11) at (BR.south west) {$1$};
\path (BR) edge [shorten >=4pt] (BR-11.center);
\node [bottom right with distance] (BR-12) at (BR.south east) {$1$};
\path (BR) edge [shorten >=4pt] (BR-12.center);
\node [above left with distance] (BR-13) at (BR.north west) {$1$};
\path (BR) edge [shorten >=4pt] (BR-13.center);
\node [above right with distance] (BR-14) at (BR.north east) {$1$};
\path (BR) edge [shorten >=4pt] (BR-14.center);
\draw[thin,dashed] (current bounding box.north west) rectangle (current bounding box.south east);
\end{tikzpicture}
\quad \quad \quad 
\begin{tikzpicture}[
		baseline={([yshift=-.5ex]current bounding box.center)},
		scale=2,very thick,
		bend angle=30,
		every loop/.style={very thick},
     		comp/.style={circle,fill,black,,inner sep=0pt,minimum size=5pt},
		prong left/.style={pos=.5,left,font=\small},
		prong right/.style={pos=.5,right,font=\small},
		circled number/.style={circle, draw, inner sep=0pt, minimum size=12pt},
		bottom right with distance/.style={below right,text height=10pt},
		bottom left with distance/.style={below left,text height=10pt},
		above right with distance/.style={above right,text height=10pt},
		above left with distance/.style={above left,text height=10pt}]
\node[comp] (BL2) {};
\node (BT) [right=of BL2] {};
\node[comp] (BR) [below=of BT] {};
\node [bottom left with distance] (BL2-2) at (BL2.south west) {$2$};
\path (BL2) edge [shorten >=4pt] (BL2-2.center);
\node [above left with distance] (BL2-11) at (BL2.north west) {$1$};
\path (BL2) edge [shorten >=4pt] (BL2-11.center);
\node [above right with distance] (BL2-12) at (BL2.north east) {$1$};
\path (BL2) edge [shorten >=4pt] (BL2-12.center);
\node [bottom left with distance] (BR-11) at (BR.south west) {$1$};
\path (BR) edge [shorten >=4pt] (BR-11.center);
\node [bottom right with distance] (BR-12) at (BR.south east) {$1$};
\path (BR) edge [shorten >=4pt] (BR-12.center);
\node [above left with distance] (BR-13) at (BR.north west) {$1$};
\path (BR) edge [shorten >=4pt] (BR-13.center);
\node [above right with distance] (BR-14) at (BR.north east) {$1$};
\path (BR) edge [shorten >=4pt] (BR-14.center);
\draw[thin,dashed] (current bounding box.north west) rectangle (current bounding box.south east);
\end{tikzpicture}\]
\caption{The ``zigzag-graph'' (left) in the boundary of the stratum $\Omega\cM_3(2,1,1)$ and an legal (center) and illegal (right) degeneration of its bottom level.}
\label{fig:zigzag}
\end{figure}
Note that there are residue conditions intertwining the simple zeros on the two components. These imply that the middle graph is legal, while the right one is not.

To check this in \verb|diffstrata|, we begin by generating the stratum and entering the graph:
\begin{lstlisting}
sage: X=Stratum((2,1,1))
sage: LG=LevelGraph([1, 1, 0, 0],[[1, 2], [3, 4], [5, 6, 7], [8, 9, 10, 11]],[(1, 6), (3, 7), (4, 10), (2, 11)],{1: 0, 2: 0, 3: 0, 4: 0, 5: 2, 6: -2, 7: -2, 8: 1, 9: 1, 10: -2, 11: -2},[0, 0, -1, -1])
sage: ELG=EmbeddedLevelGraph(X,LG,dmp={5: (0, 0), 8: (0, 1), 9: (0, 2)},dlevels={0: 0, -1: -1})
\end{lstlisting}
We may convince ourselves that the graph is as depicted in \autoref{fig:zigzag}:
\begin{lstlisting}
sage: ELG.explain()
LevelGraph embedded into stratum Stratum: (2, 1, 1)
with residue conditions: []
 with:
On level 0:
* A vertex (number 0) of genus 1
* A vertex (number 1) of genus 1
On level 1:
* A vertex (number 2) of genus 0
* A vertex (number 3) of genus 0
The marked points are on level 1.
More precisely, we have:
* Marked point (0, 0) of order 2 on vertex 2 on level 1
* Marked point (0, 1) of order 1 on vertex 3 on level 1
* Marked point (0, 2) of order 1 on vertex 3 on level 1
Finally, we have 4 edges. More precisely:
* one edge between vertex 0 (on level 0) and vertex 2 (on level 1) with prong 1.
* one edge between vertex 1 (on level 0) and vertex 2 (on level 1) with prong 1.
* one edge between vertex 1 (on level 0) and vertex 3 (on level 1) with prong 1.
* one edge between vertex 0 (on level 0) and vertex 3 (on level 1) with prong 1.
\end{lstlisting}
We may extract the bottom level of this \verb|EmbeddedLevelGraph|:
\begin{lstlisting}
sage: L=ELG.level(1); L  # extract bottom level
LevelStratum(sig_list=[Signature((2, -2, -2)), Signature((1, 1, -2, -2))],res_cond=[[(0, 1), (1, 3)], [(0, 2), (1, 2)]],leg_dict={5: (0, 0), 6: (0, 1), 7: (0, 2), 8: (1, 0), 9: (1, 1), 10: (1, 2), 11: (1, 3)})
\end{lstlisting}
Embedding (into \verb|L|!) the two degenerations depicted in \autoref{fig:zigzag}, we confirm that one is legal, while the other is not:
\begin{lstlisting}
sage: M=EmbeddedLevelGraph(L, LevelGraph([0, 0],[[1, 2, 3], [4, 5, 6, 7]],[],{1: 2, 2: -2, 3: -2, 4: 1, 5: 1, 6: -2, 7: -2},[-1, 0]),dmp={1: (0, 0), 2: (0, 1), 3: (0, 2), 4: (1, 0), 5: (1, 1), 6: (1, 2), 7: (1, 3)},dlevels={-1: -1, 0: 0})
sage: M.is_legal()
True
sage: R=EmbeddedLevelGraph(L, LevelGraph([0, 0],[[1, 2, 3], [4, 5, 6, 7]],[],{1: 2, 2: -2, 3: -2, 4: 1, 5: 1, 6: -2, 7: -2},[0, -1]),dmp={1: (0, 0), 2: (0, 1), 3: (0, 2), 4: (1, 0), 5: (1, 1), 6: (1, 2), 7: (1, 3)},dlevels={-1: -1, 0: 0})
sage: R.is_legal()
False
\end{lstlisting}
\end{example}

\subsection{The Underlying Graph}
\label{subsec:underlying_graph}

Each \verb|EmbeddedLevelGraph| has an associated ``underlying graph'', which is a Sage \verb|Graph|. It's main use is checking the graph-theoretic conditions for the $\mathfrak{R}$-GRC corresponding, essentially, to the auxiliary graph $\widetilde{\Gamma}$ constructed in the proof of \autoref{prop:realizewithR}.
To encode all necessary data, it has a slightly more involved format, which we briefly explain.

The underlying graph is, formally, associated to the underlying \verb|LevelGraph|. The vertices are formally \verb|tuple|s encoding the vertex number and its genus:
\begin{lstlisting}
sage: LG=LevelGraph([1, 1, 0, 0],[[1, 2], [3, 4], [5, 6, 7], [8, 9, 10, 11]],[(1, 6), (3, 7), (4, 10), (2, 11)],{1: 0, 2: 0, 3: 0, 4: 0, 5: 2, 6: -2, 7: -2, 8: 1, 9: 1, 10: -2, 11: -2},[0, 0, -1, -1])
sage: LG.genus(0)
1
sage: LG.UG_vertex(0)
(0, 1, 'LG')
sage: LG.underlying_graph
Looped multi-graph on 4 vertices
sage: LG.underlying_graph.vertices()
[(0, 1, 'LG'), (1, 1, 'LG'), (2, 0, 'LG'), (3, 0, 'LG')]
sage: LG.underlying_graph.edges()
[((0, 1, 'LG'), (2, 0, 'LG'), (1, 6)), ((0, 1, 'LG'), (3, 0, 'LG'), (2, 11)), ((1, 1, 'LG'), (2, 0, 'LG'), (3, 7)), ((1, 1, 'LG'), (3, 0, 'LG'), (4, 10))]
\end{lstlisting}
The entry \verb|'LG'| indicates that this vertex comes from a vertex of the \verb|LevelGraph|.

In the case of an \verb|EmbeddedLevelGraph|, the underlying \verb|LevelGraph| may be endowed with extra ``points at $\infty$'' to encode residue conditions. To illustrate this, we consider again the situation of \autoref{zigzaggraph}.
\begin{lstlisting}
sage: ELG=EmbeddedLevelGraph(X, LG, dmp={5: (0, 0), 8: (0, 1), 9: (0, 2)}, dlevels={0: 0, -1: -1})
sage: X=Stratum((2,1,1))
sage: LG=LevelGraph([1, 1, 0, 0],[[1, 2], [3, 4], [5, 6, 7], [8, 9, 10, 11]],[(1, 6), (3, 7), (4, 10), (2, 11)],{1: 0, 2: 0, 3: 0, 4: 0, 5: 2, 6: -2, 7: -2, 8: 1, 9: 1, 10: -2, 11: -2},[0, 0, -1, -1])
sage: ELG=EmbeddedLevelGraph(X, LG, dmp={5: (0, 0), 8: (0, 1), 9: (0, 2)}, dlevels={0: 0, -1: -1})
sage: L=ELG.level(1); print(L)
Product of Strata:
Signature((2, -2, -2))
Signature((1, 1, -2, -2))
with residue conditions: [(0, 1), (1, 3)] [(0, 2), (1, 2)]
dimension: 1
leg dictionary: {5: (0, 0), 6: (0, 1), 7: (0, 2), 8: (1, 0), 9: (1, 1), 10: (1, 2), 11: (1, 3)}
leg orbits: [[(0, 0)], [(0, 1), (0, 2)], [(1, 0)], [(1, 1)], [(1, 2), (1, 3)]]

sage: L.bics[0]  # numbering of list items might differ
EmbeddedLevelGraph(LG=LevelGraph([0, 0, 0],[[1, 2, 3], [4, 5, 6], [7, 8, 9]],[(6, 9)],{1: 2, 2: -2, 3: -2, 4: -2, 5: -2, 6: 2, 7: 1, 8: 1, 9: -4},[0, 0, -1],True),dmp={1: (0, 0), 2: (0, 1), 3: (0, 2), 4: (1, 2), 5: (1, 3), 7: (1, 0), 8: (1, 1)},dlevels={0: 0, -1: -1})
sage: L.bics[0].LG.underlying_graph.vertices()
[(0, 0, 'LG'), (0, 0, 'res'), (1, 0, 'LG'), (1, 0, 'res'), (2, 0, 'LG')]
sage: L.bics[0].LG.underlying_graph.edges()
[((0, 0, 'LG'), (0, 0, 'res'), 'res0edge2'), ((0, 0, 'LG'), (1, 0, 'res'), 'res1edge3'), ((0, 0, 'res'), (1, 0, 'LG'), 'res0edge5'), ((1, 0, 'LG'), (1, 0, 'res'), 'res1edge4'), ((1, 0, 'LG'), (2, 0, 'LG'), (6, 9))]
\end{lstlisting}
Each residue condition corresponds to one vertex of the underlying graph, labelled with \verb|'res'|. This vertex is connected to all poles sharing this residue condition.

Note that the residue conditions of a stratum are stored most transparently in the ``smooth'' ($0$-level) graph inside the stratum. In the above situation:
\begin{lstlisting}
sage: L.smooth_LG
EmbeddedLevelGraph(LG=LevelGraph([0, 0],[[1, 2, 3], [4, 5, 6, 7]],[],{1: 2, 2: -2, 3: -2, 4: 1, 5: 1, 6: -2, 7: -2},[0, 0],True),dmp={1: (0, 0), 2: (0, 1), 3: (0, 2), 4: (1, 0), 5: (1, 1), 6: (1, 2), 7: (1, 3)},dlevels={0: 0})
sage: L.smooth_LG.LG.underlying_graph.vertices()
[(0, 0, 'LG'), (0, 0, 'res'), (1, 0, 'LG'), (1, 0, 'res')]
sage: L.smooth_LG.LG.underlying_graph.edges()
[((0, 0, 'LG'), (0, 0, 'res'), 'res0edge2'), ((0, 0, 'LG'), (1, 0, 'res'), 'res1edge3'), ((0, 0, 'res'), (1, 0, 'LG'), 'res0edge7'), ((1, 0, 'LG'), (1, 0, 'res'), 'res1edge6')]
\end{lstlisting}

\subsection{Checking the $\mathfrak{R}$-GRC}

To check the Global Residue Condition, we use \autoref{prop:realizewithR}.
This is checked on the level of vertices (and horizontal edges, although we do not care for them here).

More precisely, \verb|LevelGraph| uses the following check for legality:
\begin{algorithm}[{name=Vertex Legality (GRC)}]
The method \verb|is_illegal_vertex| checks the legality of a vertex $v$ in the following steps:
\begin{description}
\item[Step 1] We first check if $v$ is \emph{inconvenient} in the sense of \autoref{sec:GenStrata}, that is $g_v = 0$, there are no simple poles on $v$ and there exists a pole order that is greater than the difference of the sum of the pole orders and the number of poles, i.e., denoting by $m_i$ the pole orders and by $p$ the number of poles, there exists an $m_i$ such that $m_i > \sum m_j -p-1$.

If $v$ is not inconvenient, $v$ is legal.
\item[Step 2] We next check for loops above $v$: If there are less than two edges going up from $v$, there can be no loop and $v$ is illegal. Otherwise we build the subgraph of the underlying graph that consists of vertices above the level of $v$ and generate a list of its connected components (this is easy, as this are sage \verb|Graph|s).
\item[Step 3] For each of these connected components, we check if there are two edges connecting $v$ to this connected component. (there is a technical subtlety here: the subgraph described above does not contain $v$ and hence non of the edges we are interested in. We also can't restrict to the edges of the \verb|LevelGraph| as will become clear when describing the $\frakR$-GRC below. Therefore, we actually work with two subgraphs, the \verb|abovegraph| consisting of all edges with level $\geq$ the level of $v$ and the subgraph yielding the connected components, which consists of the vertices above the level of $v$ in the connected component of \verb|abovegraph| containing $v$.) 
If this is the case, $v$ is legal.
\item[Step 4] In the case of a meromorphic stratum, $v$ is legal if there are at least two ``free'' poles on the connected component of $v$. We thus count all free poles at $v$ and in the connected components above.
\item[Step 5] If all of this fails, $v$ is illegal.
\end{description}
\end{algorithm}
In the case that there are residue conditions, we have to slightly modify Step~3 above. More precisely, also \verb|EmbeddedLevelGraph| includes an \verb|is_legal| method that checks the $\frakR$-GRC in a stratum with residue conditions.

Moreover, \verb|is_legal| also checks if any of the levels of the graph are empty (i.e. calls the appropriate \verb|is_empty| methods, described below).

Then a list of \emph{all} poles (as points on the graph) contained in \emph{any} residue condition passed to the enveloping stratum is compiled. In particular, \emph{poles are not included in this list, if they are only constrained by the Residue Theorem on some component}. This list is used in Step 4 of the above algorithm: a pole is free if it is not contained in this list, i.e. not involved in any residue condition passed to the stratum.

Note that, as we are working with the underlying graph, the vertices ``at $\infty$'' are considered for all connectivity issues (once these have been set up correctly by the \verb|EmbeddedLevelGraph|).

\begin{rem}
While most of \verb|diffstrata| works only for graphs with no horizontal edges (in particular the notion of profile does not immediately extend) \verb|LevelGraph|s may have horizontal edges and \verb|is_legal| can be used to check legality of edges, cf. \autoref{prop:realizewithR}.
\end{rem}

\subsection{Residue Matrices and Dimension}

To see how the residue conditions interact, it is often helpful to consider the matrix given by their linear equations. Moreover, for correctly degenerating residue conditions, it is important to distinguish residue conditions imposed ``from the stratum'' from those imposed on a particular graph by the Residue Theorem on each vertex.

Given a list of residue conditions and a \verb|GeneralisedStratum|, one can transform these into a matrix using the \verb|matrix_from_res_conditions| method:
\begin{lstlisting}
sage: X=GeneralisedStratum([Signature((2,-2,-2)),Signature((2,-2,-2))])
sage: X.matrix_from_res_conditions([[(0,1),(0,2),(1,2)],[(0,1),(1,1)],[(1,1),(1,2)]])
[1 1 0 1]
[1 0 1 0]
[0 0 1 1]
\end{lstlisting}
The \verb|residue_matrix| of a stratum is this matrix for the residue conditions of the stratum.

For calculating the dimension of a stratum, we need to consider \emph{all} imposed residue conditions, i.e. also the conditions imposed by the Residue Theorem on each component. This can be done for \emph{any} \verb|EmbeddedLevelGraph| by using the method \verb|residue_matrix_from_RT|, the combined matrix (with all residue conditions) is available via \verb|full_residue_matrix|.

\begin{example}
We illustrate the difference by considering the bottom level of the V-graph, the Banana graph and the zigzag graph introduced above.

Again, we input the graphs a \verb|LevelGraph|s and embed them into their corresponding strata:
\begin{lstlisting}
sage: V=LevelGraph([1, 1, 0],[[1], [2], [3, 4, 5, 6]],[(1, 5), (2, 6)],{1: 0, 2: 0, 3: 1, 4: 1, 5: -2, 6: -2},[0, 0, -1])
sage: X=Stratum((1,1))
sage: ELV=EmbeddedLevelGraph(X, V, dmp={3: (0, 0), 4: (0, 1)}, dlevels={0: 0, -1: -1})
sage: B=LevelGraph([1, 0],[[1, 2], [3, 4, 5]],[(1, 4), (2, 5)],{1: 0, 2: 0, 3: 2, 4: -2, 5: -2},[0, -1])
sage: Y=Stratum((2,))
sage: ELB=EmbeddedLevelGraph(Y, B, dmp={3: (0,0)}, dlevels={0: 0, -1: -1})
sage: Z=Stratum((2,1,1))
sage: LG=LevelGraph([1, 1, 0, 0],[[1, 2], [3, 4], [5, 6, 7], [8, 9, 10, 11]],[(1, 6), (3, 7), (4, 10), (2, 11)],{1: 0, 2: 0, 3: 0, 4: 0, 5: 2, 6: -2, 7: -2, 8: 1, 9: 1, 10: -2, 11: -2},[0, 0, -1, -1])
sage: ELZ=EmbeddedLevelGraph(Z, LG, dmp={5: (0, 0), 8: (0, 1), 9: (0, 2)}, dlevels={0: 0, -1: -1})
\end{lstlisting}
We extract the bottom levels of each graph:
\begin{lstlisting}
sage: V_bottom=ELV.level(1)
sage: B_bottom=ELB.level(1)
sage: Z_bottom=ELZ.level(1)
\end{lstlisting}
Inspecting these, we see the different residue conditions imposed on the poles:
\begin{lstlisting}
sage: print(V_bottom)
Stratum: Signature((1, 1, -2, -2))
with residue conditions: [(0, 2)] [(0, 3)]
dimension: 0
leg dictionary: {3: (0, 0), 4: (0, 1), 5: (0, 2), 6: (0, 3)}
leg orbits: [[(0, 0)], [(0, 1)], [(0, 3), (0, 2)]]

sage: print(B_bottom)
Stratum: Signature((2, -2, -2))
with residue conditions: [(0, 1), (0, 2)]
dimension: 0
leg dictionary: {3: (0, 0), 4: (0, 1), 5: (0, 2)}
leg orbits: [[(0, 0)], [(0, 1), (0, 2)]]

sage: print(Z_bottom)
Product of Strata:
Signature((2, -2, -2))
Signature((1, 1, -2, -2))
with residue conditions: [(0, 1), (1, 3)] [(0, 2), (1, 2)]
dimension: 1
leg dictionary: {5: (0, 0), 6: (0, 1), 7: (0, 2), 8: (1, 0), 9: (1, 1), 10: (1, 2), 11: (1, 3)}
leg orbits: [[(0, 0)], [(0, 1), (0, 2)], [(1, 0)], [(1, 1)], [(1, 2), (1, 3)]]
\end{lstlisting}
We can now calculate the residue matrices of the strata:
\begin{lstlisting}
sage: V_bottom.residue_matrix()
[1 0]
[0 1]
sage: B_bottom.residue_matrix()
[1 1]
sage: Z_bottom.residue_matrix()
[1 0 0 1]
[0 1 1 0]
\end{lstlisting}
and compare these to the \verb|full_residue_matrix| of the \verb|smooth_LG| inside each stratum.
\begin{lstlisting}
sage: B_bottom.smooth_LG.residue_matrix_from_RT
[1 1]
sage: V_bottom.smooth_LG.residue_matrix_from_RT
[1 1]
sage: Z_bottom.smooth_LG.residue_matrix_from_RT
[1 1 0 0]
[0 0 1 1]
sage: B_bottom.smooth_LG.full_residue_matrix
[1 1]
[1 1]
sage: V_bottom.smooth_LG.full_residue_matrix
[1 0]
[0 1]
[1 1]
sage: Z_bottom.smooth_LG.full_residue_matrix
[1 0 0 1]
[0 1 1 0]
[1 1 0 0]
[0 0 1 1]
\end{lstlisting}
\end{example}

The rank of the \verb|full_residue_matrix| is the rank of the residue space of the stratum. It is vital for calculating the dimension of a stratum. Recall \cite[Remark 4.1]{strataEC} that
\[\dim X = \sum_{i=1}^k (2g_i + n_i - 1) - \rk(M) - 1,\]
where $k$ is the number of connected components and $M$ is the full residue matrix. We use this to calculate the correct dimensions:
\begin{lstlisting}
sage: B_bottom.dim()
0
sage: V_bottom.dim()
0
sage: Z_bottom.dim()
1
\end{lstlisting}

\begin{example}
We can check that the top and bottom level of every BIC in a stratum do indeed add up to one less than the dimension:
\begin{lstlisting}
sage: X4=GeneralisedStratum([Signature((4,))])
sage: all(B.level(0).dim() + B.level(1).dim() == X4.dim() - 1 for B in X4.bics)
True
\end{lstlisting}
\end{example}

\begin{rem}
We can also check if the residue at a pole is forced $0$ by the residue conditions (i.e. if the rank of the \verb|full_residue_matrix| changes by adding this condition). In the case that a stratum has simple poles, we use this to determine if the stratum is empty (if a simple pole is forced $0$):
\begin{lstlisting}
sage: I=Stratum((1,-1))
sage: I.smooth_LG.residue_zero((0,1))
True
sage: I.is_empty()
True
\end{lstlisting}
\end{rem}

\subsection{Level Extraction}

We have already seen level extraction in action and describe the process now in some more detail.

Note that, strictly speaking, an extracted level is a \verb|LevelStratum|, which inherits from \verb|GeneralisedStratum|.

Revisiting the V-graph from above:
\begin{lstlisting}
sage: ELV.level(1)
LevelStratum(sig_list=[Signature((1, 1, -2, -2))],res_cond=[[(0, 2)], [(0, 3)]],leg_dict={3: (0, 0), 4: (0, 1), 5: (0, 2), 6: (0, 3)})
sage: type(ELV.level(1))
<class 'admcycles.diffstrata.levelstratum.LevelStratum'>
sage: isinstance(ELV.level(1), GeneralisedStratum)
True
\end{lstlisting}

The main extra piece of information it holds is \verb|leg_dict|, a dictionary mapping the legs of the \verb|LevelGraph| we extracted the level of (in this case \verb|V|) to the marked points of the \verb|LevelStratum|. Also, the orbits of the automorphism group are recorded:
\begin{lstlisting}
sage: ELV
EmbeddedLevelGraph(LG=LevelGraph([1, 1, 0],[[1], [2], [3, 4, 5, 6]],[(1, 5), (2, 6)],{1: 0, 2: 0, 3: 1, 4: 1, 5: -2, 6: -2},[0, 0, -1],True),dmp={3: (0, 0), 4: (0, 1)},dlevels={0: 0, -1: -1})
sage: print(ELV.level(1))
Stratum: Signature((1, 1, -2, -2))
with residue conditions: [(0, 2)] [(0, 3)]
dimension: 0
leg dictionary: {3: (0, 0), 4: (0, 1), 5: (0, 2), 6: (0, 3)}
leg orbits: [[(0, 0)], [(0, 1)], [(0, 3), (0, 2)]]
\end{lstlisting}

The main part of the extraction happens again on the level of the \verb|LevelGraph|, more precisely in the method \verb|stratum_from_level|.
\begin{algorithm}[Level Extraction] 
The method \verb|stratum_from_level| is implemented as follows:
\begin{description}
\item[Step 1] Extraction of the relevant data from the graph: vertices at the current level as well as the legs and their orders. Here, we must be careful with the level numbering, cf. \autoref{levelnumbers}. As we sort the marked points into their new signatures, we create \verb|leg_dict|.
\item[Step 2] Adding the residue conditions is implemented analogous to the $\frakR$-GRC check: the poles on the level are sorted by the connected component of the underlying graph they are attached to. Recall that the underlying graph includes the ``level at $\infty$''. Again, components containing free poles must be ignored. Moreover, the poles need to translated into their respective strata points for storing the residue conditions.
\end{description}
\end{algorithm}
When a level is extracted from an \verb|EmbeddedLevelGraph| (as is usually the case), first the free poles determined by the residue conditions of the ``big'' graph, as described above and then \verb|stratum_from_level| is called on the underlying \verb|LevelGraph| with this information.

The resulting \verb|LevelStratum| is then equipped with the automorphism orbits of the legs and cached.

\begin{rem}\label{topbottom}
In the case of a BIC, i.e. an \verb|EmbeddedLevelGraph| with exactly two levels, the levels are also stored in the \verb|top| and \verb|bot| attributes:
\begin{lstlisting}
sage: ELV.is_bic()
True
sage: ELV.top == ELV.level(0)
True
sage: ELV.bot == ELV.level(1)
True
\end{lstlisting}
\end{rem}

%% file: sec_generating.tex
\section{Generating all Level Graphs}
\label{sec:AllGraphs}

In this section, we give a detailed description of the explicit construction of \autoref{prop:gengraphs}.

Generating all level graphs (as \verb|EmbeddedLevelGraph|s) inside a generalised stratum is done recursively: first, we generate all two-level graphs (BICs) and then we clutch them together in all possible ways. This gives the degeneration graph as a by-product.

We can inspect this data for any stratum:
\begin{lstlisting}
sage: X=Stratum((4,))
sage: X.info()
Stratum: (4,)
with residue conditions: []

Genus: [3]
Dimension: 5
Boundary Graphs (without horizontal edges):
Codimension 0: 1 graph
Codimension 1: 8 graphs
Codimension 2: 19 graphs
Codimension 3: 16 graphs
Codimension 4: 4 graphs
Total graphs: 48
\end{lstlisting}
Note that the number of graphs increases quickly:
\begin{lstlisting}
sage: X=Stratum((1,1,1,1))
sage: X.info()
Stratum: (1, 1, 1, 1)
with residue conditions: []

Genus: [3]
Dimension: 8
Boundary Graphs (without horizontal edges):
Codimension 0: 1 graph
Codimension 1: 102 graphs
Codimension 2: 1100 graphs
Codimension 3: 4222 graphs
Codimension 4: 7531 graphs
Codimension 5: 6708 graphs
Codimension 6: 2856 graphs
Codimension 7: 456 graphs
Total graphs: 22976
\end{lstlisting}
\begin{rem}
For genus $4$, strata having several million boundary components appear. Already the generation of the strata requires significant computational and memory power.
\end{rem}

The key observation is that, by \autoref{prop:beta}, each three-level graph arises by either clutching a top level of a BIC to a BIC of its bottom level or a BIC of the top level to the bottom level.
On the other hand, each three-level graph is contained in the intersection of two (different) BICs of the Stratum.

More explicitly, any three-level graph is contained in a profile $(\beta^\top_\Gamma(\Gamma'),\Gamma)$ for some $\Gamma\in\BIC(X)$ and $\Gamma'\in\BIC(\Gamma^\top)$ (or equivalently $(\Lambda,\beta^\bot_{\Lambda}(\Lambda'))$ for some $\Lambda\in\BIC(X)$ and $\Lambda'\in\BIC(\Lambda^\bot)$.

We can thereby express the clutching of a product of BICs in the top
and bottom components of a BIC in our stratum as a product of BICs of
our stratum. Hence the procedure is recursive.

Therefore, to generate all graphs, we need to generate only the BICs
together with, for each BIC, top and bottom components and the two maps.

More precisely, the degeneration graph of a \verb|GeneralisedStratum| $X$
is determined by the following information:
\begin{itemize}
\item The set $\BIC(X)$ corresponding to the \verb|list| \verb|X.bics|.
\item For each $\Gamma\in\BIC(X)$, its top and bottom component $\Gamma^\top$ and $\Gamma^\bot$ (each a \verb|GeneralisedStratum| together with a dictionary mapping Stratum points to \verb|LevelGraph| points).
If $\Gamma$ corresponds to an \verb|EmbeddedLevelGraph| \verb|B|, $\Gamma^\top$ is given by \verb|B.top| and $\Gamma^\bot$ is given by \verb|B.bot|.
\item For each $\Gamma\in\BIC(X)$, the sets $\BIC(\Gamma^\top)$ and $\BIC(\Gamma^\bot)$ together with the maps $\beta_\Gamma^\top$ and $\beta_\Gamma^\bot$.
The maps $\beta$ are given as dictionaries of indices of \verb|list|s: if $\Gamma$ corresponds to \verb|X.bics[i]|, then
\begin{itemize}
\item \verb|top_to_bic(i)| is a \verb|dict| mapping indices of \verb|X.bics[i].top.bics| to indices of \verb|X.bics| such that \verb|X.DG.top_to_bic(i)[j]| corresponds to $\beta^\top_\Gamma(\Gamma')$ if $\Gamma'$ corresponds to \verb|X.bics[i].top.bics[j]|
and 
\item \verb|bot_to_bic(i)| is a \verb|dict| mapping indices of \verb|X.bics[i].bot.bics| to indices of \verb|X.bics| such that \verb|X.DG.bot_to_bic(i)[j]| corresponds to $\beta^\bot_\Gamma(\Gamma')$ if $\Gamma'$ corresponds to \verb|X.bics[i].bot.bics[j]|.
\end{itemize}
\end{itemize}
See \autoref{ex:toptobic} for details.

\begin{rem}\label{numbering}
For caching purposes, it is essential for \verb|diffstrata| to refer to a BIC by its \emph{index} in the \verb|list| \verb|X.bics|.

Correspondingly, a profile is encoded as a \verb|tuple| of \verb|int|s, an enhanced profile is encoded as a nested \verb|tuple|, \verb|(p, i)|, where \verb|p| is a profile and \verb|i| is an \verb|int| determining the component.
As described above, the maps $\beta$ are also given by \verb|dict|s mapping indices to indices.

Note, however, that starting with Sage 9, the numbering of profiles can (and will!) change between every \verb|sage| session!
\end{rem}

We can now calculate the degeneration graph.
\begin{algorithm}[Degeneration Graph]\mbox{}\label{alg:tbtobic}
\begin{description}
\item[Step 1] Calculate all BICs in a \verb|GeneralisedStratum|.
\item[Step 2] Separate these into top and bottom component.
\item[Step 3] Calculate all BICs in every top and bottom component.
\item[Step 4] Calculate \verb|top_to_bic| and \verb|bot_to_bic| for each BIC in the
Stratum (as dictionaries mapping an index of \verb|bics| of \verb|top|/\verb|bot| to
an index of \verb|bics| of the stratum).
\end{description}
\end{algorithm}

In particular, this yields a  recursive algorithm for
the \verb|EmbeddedLevelGraph| of an arbitrary product of BICs in
the stratum as follows:
\begin{algorithm}[Building a Graph from BICs]\mbox{}\label{alg:profilesfrombics}
\begin{description}
\item[INPUT] Product of BICs.
\item[OUTPUT] \verb|EmbeddedLevelGraph|.
\item[Step 1] Choose a BIC $B$ from the product (e.g. the first).
\item[Step 2] Find the preimages of the other BICs in the product under
\verb|top_to_bic| and \verb|bot_to_bic| of $B$.

This gives (possibly multiple) products of BICs in the top and bottom
stratum of B.
\item[Step 3] Apply to product in top an bottom to get two \verb|EmbeddedLevelGraphs|.
\item[Step 4] Return the clutching of the top and bottom graph.
\end{description}
\end{algorithm}

Moreover, we can generate the ``lookup list'', consisting of the non-empty
profiles in each stratum.

Before we discuss the details of the implementation, we give an overview of the relevant interface and some examples.

\subsection{Interface and Examples}
\label{deg:interface}

Let \verb|X| be a \verb|GeneralisedStratum|. To access an \verb|EmbeddedLevelGraph| inside \verb|X|, we need to provide an {enhanced profile}, that is
\begin{itemize}
\item a \verb|tuple| of indices of \verb|X.bics|, the {profile}, and
\item an \verb|int| as a choice of component, in case this profile is reducible.
\end{itemize}
The relevant methods are
\begin{itemize}
\item \verb|X.lookup| with argument a profile \verb|tuple| and returns a \verb|list| of all the \verb|EmbeddedLevelGraph|s matching this profile (even if the profile is irreducible, this returns a \verb|list| of length 1!).
\item \verb|X.lookup_graph| with argument a profile \verb|tuple| and optional argument a choice of connected component to be provided as an \verb|int| index of \verb|X.lookup|, which defaults to \verb|0|.
\item \verb|X.lookup_list| is a nested \verb|list|, where \verb|X.lookup_list[codim]| is a \verb|list| of profiles of length \verb|codim| (i.e. the corresponding level graph is of codimension \verb|codim|).
\item \verb|X.enhanced_profiles_of_length| with argument an \verb|int| returns a list of all enhanced profiles of this length.
\end{itemize}

\begin{example}\label{H2}
Let us consider the stratum $\Omega\cM_{2}(2)$. The degeneration graph consists of the three graphs depicted in \autoref{degenerationgraphH2}
\begin{figure}
\begin{tikzpicture}[
		baseline={([yshift=-.5ex]current bounding box.center)},
		scale=2,very thick,
		bend angle=30,
		every loop/.style={very thick},
     		comp/.style={circle,fill,black,,inner sep=0pt,minimum size=5pt},
		order bottom left/.style={pos=.05,left,font=\tiny},
		order top left/.style={pos=.9,left,font=\tiny},
		order bottom right/.style={pos=.05,right,font=\tiny},
		order top right/.style={pos=.9,right,font=\tiny},
		circled number/.style={circle, draw, inner sep=0pt, minimum size=12pt},
		bottom right with distance/.style={below right,text height=10pt}]

\begin{scope}[local bounding box = l]
\node[circled number] (T) [] {$1$}; 
\node[comp] (B) [below=of T] {}
	edge [bend left] 
		node [order bottom left] {$-2$} 
		node [order top left] {$0$} (T)
	edge [bend right] 
		node [order bottom right] {$-2$} 
		node [order top right] {$0$} (T);
\node [bottom right with distance] (B-2) at (B.south east) {$2$};
\path (B) edge [shorten >=4pt] (B-2.center);
\end{scope}

\begin{scope}[shift={($(l.base) + (1.5,0)$)}, local bounding box = ct]
\node[circled number] (T) [] {$1$}; 
\node[circled number] (B) [below=of T] {$1$}
	edge 
		node [order bottom left] {$-2$} 
		node [order top left] {$0$} (T);
\node [bottom right with distance] (B-2) at (B.south east) {$2$};
\path (B) edge [shorten >=4pt] (B-2.center);
\end{scope}

\begin{scope}[shift={($(l.base) + (.75,-1.2)$)}, local bounding box = ctl]
\node[circled number] (T) [] {$1$}; 
\node[comp] (B1) [below=of T] {}
	edge 
		node [order bottom left] {$-2$} 
		node [order top left] {$0$} (T);
\node[comp] (B2) [below=of B1] {}
	edge [bend left] 
		node [order bottom left] {$-2$} 
		node [order top left] {$0$} (B1)
	edge [bend right] 
		node [order bottom right] {$-2$} 
		node [order top right] {$0$} (B1);
\node [bottom right with distance] (B2-2) at (B2.south east) {$2$};
\path (B2) edge [shorten >=4pt] (B2-2.center);
\end{scope}
\draw[thin] (l) edge (ctl);
\draw[thin] (ct) edge (ctl);
\end{tikzpicture}
\caption{The degeneration graph of $\Omega\cM_2(2)$.}
\label{degenerationgraphH2}
\end{figure}
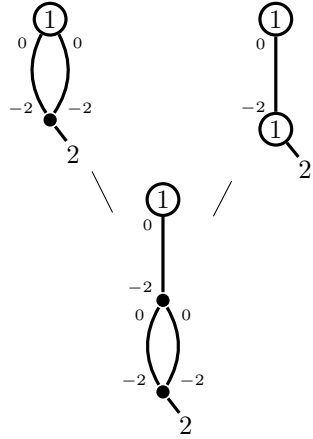

In \verb|diffstrata|, the situation is as follows. We first generate the stratum and may then inspect the BICs:
\begin{lstlisting}
sage: X=Stratum((2,))
sage: X.bics  # order might differ
[EmbeddedLevelGraph(LG=LevelGraph([1, 0],[[1, 2], [3, 4, 5]],[(1, 4), (2, 5)],{1: 0, 2: 0, 3: 2, 4: -2, 5: -2},[0, -1],True),dmp={3: (0, 0)},dlevels={0: 0, -1: -1}),
 EmbeddedLevelGraph(LG=LevelGraph([1, 1],[[1], [2, 3]],[(1, 3)],{1: 0, 2: 2, 3: -2},[0, -1],True),dmp={2: (0, 0)},dlevels={0: 0, -1: -1})]
\end{lstlisting}
As these are \verb|EmbeddedLevelGraph|s, we may use the \verb|explain| method to relate them to \autoref{degenerationgraphH2}:
\begin{lstlisting}
sage: X.bics[0].explain()
LevelGraph embedded into stratum Stratum: (2,)
with residue conditions: []
 with:
On level 0:
* A vertex (number 0) of genus 1
On level 1:
* A vertex (number 1) of genus 0
The marked points are on level 1.
More precisely, we have:
* Marked point (0, 0) of order 2 on vertex 1 on level 1
Finally, we have 2 edges. More precisely:
* 2 edges between vertex 0 (on level 0) and vertex 1 (on level 1) with prongs 1 and 1.
sage: X.bics[1].explain()
LevelGraph embedded into stratum Stratum: (2,)
with residue conditions: []
 with:
On level 0:
* A vertex (number 0) of genus 1
On level 1:
* A vertex (number 1) of genus 1
The marked points are on level 1.
More precisely, we have:
* Marked point (0, 0) of order 2 on vertex 1 on level 1
Finally, we have one edge. More precisely:
* one edge between vertex 0 (on level 0) and vertex 1 (on level 1) with prong 1.
\end{lstlisting}
The profiles are stored in the \verb|lookup_list| of \verb|X|:
\begin{lstlisting}
sage: X.lookup_list
[[()], [(0,), (1,)], [(1, 0)]]
\end{lstlisting}
We can inspect the graph in a profile:
\begin{lstlisting}
sage: X.lookup_graph((1,0))   # this might be (0,1) !!
EmbeddedLevelGraph(LG=LevelGraph([1, 0, 0],[[1], [2, 3, 4], [5, 6, 7]],[(1, 4), (2, 6), (3, 7)],{1: 0, 2: 0, 3: 0, 4: -2, 5: 2, 6: -2, 7: -2},[0, -1, -2],True),dmp={5: (0, 0)},dlevels={0: 0, -1: -1, -2: -2})
sage: X.lookup((1,0))  # list
[EmbeddedLevelGraph(LG=LevelGraph([1, 0, 0],[[1], [2, 3, 4], [5, 6, 7]],[(1, 4), (2, 6), (3, 7)],{1: 0, 2: 0, 3: 0, 4: -2, 5: 2, 6: -2, 7: -2},[0, -1, -2],True),dmp={5: (0, 0)},dlevels={0: 0, -1: -1, -2: -2})]
\end{lstlisting}
Again, the \verb|explain| method relates this graph to \autoref{degenerationgraphH2}:
\begin{lstlisting}
sage: X.lookup_graph((1,0)).explain()
LevelGraph embedded into stratum Stratum: (2,)
with residue conditions: []
 with:
On level 0:
* A vertex (number 0) of genus 1
On level 1:
* A vertex (number 1) of genus 0
On level 2:
* A vertex (number 2) of genus 0
The marked points are on level 2.
More precisely, we have:
* Marked point (0, 0) of order 2 on vertex 2 on level 2
Finally, we have 3 edges. More precisely:
* one edge between vertex 0 (on level 0) and vertex 1 (on level 1) with prong 1.
* 2 edges between vertex 1 (on level 1) and vertex 2 (on level 2) with prongs 1 and 1.
\end{lstlisting}
Note that we should always keep \autoref{numbering} in mind when working with (enhanced) profiles!
We convince ourselves, that the profile works as described, i.e. are compatible with the maps $\delta_i$:
\begin{lstlisting}
sage: G=X.lookup_graph((1,0))
sage: G.delta(1)
EmbeddedLevelGraph(LG=LevelGraph([1, 1],[[1], [4, 5]],[(1, 4)],{1: 0, 4: -2, 5: 2},[0, -1],True),dmp={5: (0, 0)},dlevels={0: 0, -1: -1})
sage: G.delta(1).is_isomorphic(X.bics[1])
True
sage: G.delta(2)
EmbeddedLevelGraph(LG=LevelGraph([1, 0],[[2, 3], [5, 6, 7]],[(2, 6), (3, 7)],{2: 0, 3: 0, 5: 2, 6: -2, 7: -2},[0, -2],True),dmp={5: (0, 0)},dlevels={0: 0, -2: -1})
sage: G.delta(2).is_isomorphic(X.bics[0])
True\end{lstlisting}

We can also list the enhanced profiles and convince ourselves that all profiles are irreducible in this case.
\begin{lstlisting}
sage: X.enhanced_profiles_of_length(0)
(((), 0),)
sage: X.enhanced_profiles_of_length(1)
(((0,), 0), ((1,), 0))
sage: X.enhanced_profiles_of_length(2)
(((1, 0), 0),)
\end{lstlisting}
Note that the empty profile always corresponds to the $0$-level graph, \verb|smooth_LG|:
\begin{lstlisting}
sage: X.lookup_graph(()) == X.smooth_LG
True
\end{lstlisting}
\end{example}

While the profile gives a good first approximation of an \EmbeddedLevelGraph{}, there can be distinct level graphs sharing the same profile.

\begin{example}\label{reducibleprofile}
Consider the following level graphs in the boundary of $\Omega\cM_3(4)$ depicted in \autoref{fig:reducible}.

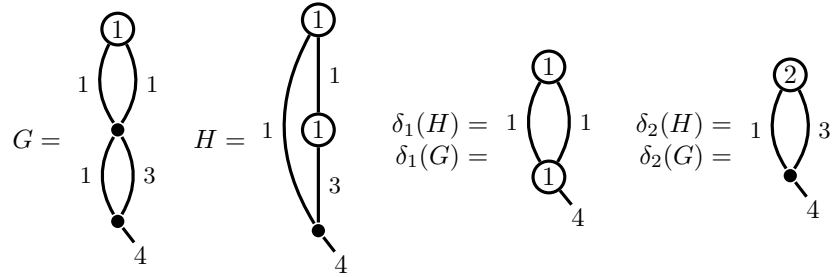
\begin{figure}
\[
G = \begin{tikzpicture}[
		baseline={([yshift=-.5ex]current bounding box.center)},
		scale=2,very thick,
		bend angle=30,
		every loop/.style={very thick},
     		comp/.style={circle,fill,black,,inner sep=0pt,minimum size=5pt},
		prong left/.style={pos=.5,left,font=\small},
		prong right/.style={pos=.5,right,font=\small},
		circled number/.style={circle, draw, inner sep=0pt, minimum size=12pt},
		bottom right with distance/.style={below right,text height=10pt}]
\node[circled number] (T) {$1$};
\node[comp] (B1) [below=of T] {}
	edge [bend left] node [prong left] {$1$} (T)
	edge [bend right] node [prong right] {$1$} (T);
\node[comp] (B2) [below=of B1] {}
	edge [bend left] node [prong left] {$1$} (B1)
	edge [bend right] node [prong right] {$3$} (B1);
\node [bottom right with distance] (B2-4) at (B2.south east) {$4$};
\path (B2) edge [shorten >=4pt] (B2-4.center);
\end{tikzpicture}
\quad
H = \begin{tikzpicture}[
		baseline={([yshift=-.5ex]current bounding box.center)},
		scale=2,very thick,
		bend angle=30,
		every loop/.style={very thick},
     		comp/.style={circle,fill,black,,inner sep=0pt,minimum size=5pt},
		prong left/.style={pos=.5,left,font=\small},
		prong right/.style={pos=.5,right,font=\small},
		circled number/.style={circle, draw, inner sep=0pt, minimum size=12pt},
		bottom right with distance/.style={below right,text height=10pt}]
\node[circled number] (T) {$1$};
\node[circled number] (B1) [below=of T] {$1$}
	edge node [prong right] {$1$} (T);
\node[comp] (B2) [below=of B1] {}
	edge [bend left] node [prong left] {$1$} (T)
	edge node [prong right] {$3$} (B1);
\node [bottom right with distance] (B2-4) at (B2.south east) {$4$};
\path (B2) edge [shorten >=4pt] (B2-4.center);
\end{tikzpicture}
\quad
\newlength\deltawidth
\settowidth\deltawidth{$\delta_1(H)=$ }
\parbox{\deltawidth}{\raggedleft$\delta_1(H)=$\\$\delta_1(G)=$}\ 
\begin{tikzpicture}[
		baseline={([yshift=-.5ex]current bounding box.center)},
		scale=2,very thick,
		bend angle=30,
		every loop/.style={very thick},
     		comp/.style={circle,fill,black,,inner sep=0pt,minimum size=5pt},
		prong left/.style={pos=.5,left,font=\small},
		prong right/.style={pos=.5,right,font=\small},
		circled number/.style={circle, draw, inner sep=0pt, minimum size=12pt},
		bottom right with distance/.style={below right,text height=10pt}]
\node[circled number] (T) {$1$};
\node[circled number] (B) [below=of T] {$1$}
	edge [bend left] node [prong left] {$1$} (T)
	edge [bend right] node [prong right] {$1$} (T);
\node [bottom right with distance] (B-4) at (B.south east) {$4$};
\path (B) edge [shorten >=4pt] (B-4.center);
\end{tikzpicture}
\quad
\settowidth\deltawidth{$\delta_2(H)=$ }
\parbox{\deltawidth}{\raggedleft$\delta_2(H)=$\\$\delta_2(G)=$}\ 
\begin{tikzpicture}[
		baseline={([yshift=-.5ex]current bounding box.center)},
		scale=2,very thick,
		bend angle=30,
		every loop/.style={very thick},
     		comp/.style={circle,fill,black,,inner sep=0pt,minimum size=5pt},
		prong left/.style={pos=.5,left,font=\small},
		prong right/.style={pos=.5,right,font=\small},
		circled number/.style={circle, draw, inner sep=0pt, minimum size=12pt},
		bottom right with distance/.style={below right,text height=10pt}]
\node[circled number] (T) {$2$};
\node[comp] (B) [below=of T] {}
	edge [bend left] node [prong left] {$1$} (T)
	edge [bend right] node [prong right] {$3$} (T);
\node [bottom right with distance] (B-4) at (B.south east) {$4$};
\path (B) edge [shorten >=4pt] (B-4.center);
\end{tikzpicture}
\]
\caption{An example of a reducible profile: the two non-isomorphic three-level graphs $G$ and $H$ have the same two-level graphs as undegenerations.}
\label{fig:reducible}
\end{figure}

In light of \autoref{numbering}, it is always a good idea to use intrinsic properties of the graphs to safely retrieve their profile:
\begin{lstlisting}
sage: X=GeneralisedStratum([Signature((4,))])
sage: left_graph = [ep for ep in X.enhanced_profiles_of_length(2)
....: if set(X.lookup_graph(*ep).LG.prongs.values()) == set([1,3]) \
....:    and len(X.lookup_graph(*ep).LG.edges) == 4 \
....:    and not X.lookup_graph(*ep).LG.has_long_edge()]
sage: assert len(left_graph) == 1
sage: left_graph = left_graph[0]
sage: right_graph = [ep for ep in X.enhanced_profiles_of_length(2)
....: if set(X.lookup_graph(*ep).LG.prongs.values()) == set([1,3]) \
....:    and len(X.lookup_graph(*ep).LG.edges) == 3 \
....:    and X.lookup_graph(*ep).LG.has_long_edge()]
sage: assert len(right_graph) == 1
sage: right_graph = right_graph[0]
\end{lstlisting}
We check that the undegenerations agree as claimed:
\begin{lstlisting}
sage: assert right_graph[0] == left_graph[0]  # compare profiles
sage: assert X.lookup_graph(*right_graph).delta(1).is_isomorphic(X.bics[right_graph[0][0]])   # check profiles are compatible with delta as claimed
sage: assert X.lookup_graph(*right_graph).delta(2).is_isomorphic(X.bics[right_graph[0][1]])
sage: assert X.lookup_graph(*left_graph).delta(1).is_isomorphic(X.bics[left_graph[0][0]])
sage: assert X.lookup_graph(*left_graph).delta(2).is_isomorphic(X.bics[left_graph[0][1]])
\end{lstlisting}
We confirm that the found graphs correspond to \autoref{fig:reducible}:
\begin{lstlisting}
sage: X.lookup_graph(*left_graph).explain()
LevelGraph embedded into stratum Stratum: (4,)
with residue conditions: []
 with:
On level 0:
* A vertex (number 0) of genus 1
On level 1:
* A vertex (number 1) of genus 0
On level 2:
* A vertex (number 2) of genus 0
The marked points are on level 2.
More precisely, we have:
* Marked point (0, 0) of order 4 on vertex 2 on level 2
Finally, we have 4 edges. More precisely:
* 2 edges between vertex 0 (on level 0) and vertex 1 (on level 1) with prongs 1 and 1.
* 2 edges between vertex 1 (on level 1) and vertex 2 (on level 2) with prongs 3 and 1.
sage: X.lookup_graph(*right_graph).explain()
LevelGraph embedded into stratum Stratum: (4,)
with residue conditions: []
 with:
On level 0:
* A vertex (number 0) of genus 1
On level 1:
* A vertex (number 1) of genus 1
On level 2:
* A vertex (number 2) of genus 0
The marked points are on level 2.
More precisely, we have:
* Marked point (0, 0) of order 4 on vertex 2 on level 2
Finally, we have 3 edges. More precisely:
* one edge between vertex 0 (on level 0) and vertex 1 (on level 1) with prong 1.
* one edge between vertex 1 (on level 1) and vertex 2 (on level 2) with prong 3.
* one edge between vertex 0 (on level 0) and vertex 2 (on level 2) with prong 1.
\end{lstlisting}

One can (and should) inspect the enhanced profile stored in \verb|right_graph| and \verb|left_graph| but, as mentioned above, the actual numbers will differ between \verb|sage| sessions, while the above code should always retrieve the same graphs.

For the record, in $\Omega\cM_3(4)$, among the $19$ three-level graphs, $15$ have an irreducible profile, while there are two pairs sharing a profile:
\begin{lstlisting}
sage: len(X.enhanced_profiles_of_length(2))
19
sage: len([p for p in X.lookup_list[2] if len(X.lookup(p)) == 1])
15
sage: len([p for p in X.lookup_list[2] if len(X.lookup(p)) == 2])
2
\end{lstlisting}
\end{example}

We end this section with a brief summary of the various methods around undegeneration. They come, broadly speaking, in two groups: the first are methods of \verb|EmbeddedLevelGraph| and thus use the explicit graphs; the second are methods of \verb|GeneralisedStratum| and deal exclusively with (enhanced) profiles.

In \verb|EmbeddedLevelGraph|:
\begin{itemize}
\item \verb|squish_vertical| returns the \verb|EmbeddedLevelGraph| where the supplied level-crossing has been squished.
\item \verb|delta| returns the BIC (as an \verb|EmbeddedLevelGraph|) where all but the supplied level-crossings have been squished.
\end{itemize}

In \verb|GeneralisedStratum|:
\begin{itemize}
\item \verb|squish| may be applied to an \emph{enhanced profile} and performs the vertical squish in the language of enhanced profiles.
\item \verb|is_degeneration| may be applied to two \emph{enhanced profiles} and answers if the first is a degeneration of the second.
\item \verb|lies_over| may be applied to two BICs \verb|i| and \verb|j| (as indices of \verb|bics|, i.e. \verb|int|s) and answers if \verb|(i, j)| is a non-empty profile.
\item \verb|merge_profiles| takes two \emph{profiles} (i.e. \verb|tuple|s of indices of \verb|bics|) and merges them with respect to the ordering of \verb|lies_over|.
\item \verb|common_degenerations| may be applied to a pair of \emph{enhanced profiles} and returns a \verb|list| of \emph{enhanced profiles} that occur as common degenerations.
\item \verb|codim_one_degenerations| may be applied to an \emph{enhanced profile} and returns a \verb|list| of \emph{enhanced profiles} corresponding to all degenerations with exactly one more level.
\item \verb|codim_one_common_undegenerations| may be applied to three \emph{enhanced profiles} and gives a list of codimension one degenerations of the third, that are also degenerations of the other two \emph{enhanced profiles} provided.
\item \verb|minimal_common_undegeneration| given two \emph{enhanced profiles} return the \emph{enhanced profile} of the graph of minimal dimension that is a common degeneration of both.
\end{itemize}
More details can be found in the \verb|docstring|s of the methods.

We now explain the relevant steps of the implementation in more detail.

\subsection{Generating all BICs}
\label{subsec:bicgeneration}

The generation of all BICs is again accomplished by \verb|diffstrata| in several steps: first, a list of all combinatorially possible $2$-level graphs is generated for each component of the stratum, then the GRC is checked for each graph and these are then combined to give all possible BICs for a \verb|GeneralisedStratum|. Finally, they are sorted by isomorphism class and checked against the $\frakR$-GRC.

For a \verb|GeneralisedStratum| $X$, the BICs are automatically generated (as a list of \verb|EmbeddedLevelGraph|s) and accessed through \verb|X.bics|. Note that in Sage 9 the numbering of this list changes with every \verb|sage| session. The BICs are generated by \verb|GeneralisedStratum.gen_bics()|.

\begin{lstlisting}
sage: X=GeneralisedStratum([Signature((0,0))])
sage: X.bics
[EmbeddedLevelGraph(LG=LevelGraph([1, 0],[[1], [2, 3, 4]],[(1, 4)],{1: 0, 2: 0, 3: 0, 4: -2},[0, -1],True),dmp={2: (0, 0), 3: (0, 1)},dlevels={0: 0, -1: -1})]
sage: X=GeneralisedStratum([Signature((2,))])
sage: X.bics # order will change!
[EmbeddedLevelGraph(LG=LevelGraph([1, 0],[[1, 2], [3, 4, 5]],[(1, 4), (2, 5)],{1: 0, 2: 0, 3: 2, 4: -2, 5: -2},[0, -1],True),dmp={3: (0, 0)},dlevels={0: 0, -1: -1}),
 EmbeddedLevelGraph(LG=LevelGraph([1, 1],[[1], [2, 3]],[(1, 3)],{1: 0, 2: 2, 3: -2},[0, -1],True),dmp={2: (0, 0)},dlevels={0: 0, -1: -1})]
sage: assert X.bics == X.gen_bic()
\end{lstlisting}

\subsubsection*{Combinatorial enumeration}

The combinatorial enumeration occurs inside the \verb|bic| module of \verb|diffstrata|. More precisely, the raw list of \verb|LevelGraph|s is generated by \verb|bic.bic_alt_noiso| (which needs to be provided a \verb|Signature| or signature tuple as an argument).

\begin{lstlisting}
sage: from admcycles.diffstrata.bic import bic_alt_noiso
sage: bic_alt_noiso((1,1))
[LevelGraph([1, 1],[[3], [1, 2, 4]],[(3, 4)],{1: 1, 2: 1, 4: -2, 3: 0},[0, -1],True),
 LevelGraph([2, 0],[[3], [1, 2, 4]],[(3, 4)],{1: 1, 2: 1, 4: -4, 3: 2},[0, -1],True),
 LevelGraph([1, 1, 0],[[3], [4], [1, 2, 5, 6]],[(3, 5), (4, 6)],{1: 1, 2: 1, 5: -2, 6: -2, 3: 0, 4: 0},[0, 0, -1],True),
 LevelGraph([1, 1, 0],[[4], [3], [1, 2, 5, 6]],[(3, 5), (4, 6)],{1: 1, 2: 1, 5: -2, 6: -2, 3: 0, 4: 0},[0, 0, -1],True),
 LevelGraph([1, 0],[[3, 4], [1, 2, 5, 6]],[(3, 5), (4, 6)],{1: 1, 2: 1, 5: -2, 6: -2, 3: 0, 4: 0},[0, -1],True)]
\end{lstlisting}

Note that in Sage 9 the ordering of these is no longer deterministic, i.e. the BICs will be produced in a different ordering in every \verb|sage| session.

To describe the algorithm, we denote by $z$ the number of zeros and by $\bz$ the corresponding vector of these zeros. By analogy, we denote by $p$ the number of poles, by $\bp$ the vector of poles and by $n$ the number of marked points (the marked points are all of order zero and thus do not need to be distinguished here). Furthermore, we denote the associated genus by $g$.

We recall some elementary bounds: We denote the maximal sum of genera of vertices on bottom level by \verb|g_bot_max| and the minimal sum of genera of vertices on top level by \verb|g_top_min|.
As every top-level vertex $v$ with $g_v=0$ requires at least one pole, $\verb|g_bot_max| = g-1$ and $\verb|g_top_min|=1$ for holomorphic strata and correspondingly $\verb|g_bot_max| = g$ and $\verb|g_top_min|=0$ for meromorphic strata.

We now summarise the actual algorithm. The steps correspond mostly to nested loops.
\begin{algorithm}[BIC Generation]\mbox{}\label{BICgenalgo}
\begin{description}
\item[Step 1] We begin by iterating over \verb|bot_comp_len|, the possible number of vertices on bottom level. As every bottom-level vertex needs at minimum either a zero or (if it's genus $0$ and has only one edge going up) two marked points, we note that $z+n$ is an upper bound for \verb|bot_comp_len|.
\item[Step 2] Next, we distribute the zeros between upper and lower level by iterating over all 2-length partitions of $\bz$ (and also include the case where all zeros are on bottom level), ensuring at each step that we have enough zeros to satisfy \verb|bot_comp_len|.
\item[Step 3] The zeros are distributed onto the bottom components. If there are no marked points, every component needs at least one zero, otherwise we are more flexible. We use the appropriate helper function, \verb|_distribute_fully| (partitions composed with permutations) or \verb|_distribute_points| (essentially a powerset, implemented via computing $b$-ary representations of numbers up to $b^l$ and using this to place the points on the components).
\item[Step 4] The genus is partitioned into the contribution from the top vertices, the bottom vertices and the graph. For this, we iterate over \verb|total_g_bot| and \verb|total_g_bot| (using the bounds described above).
\item[Step 5] \verb|total_g_bot| is distributed onto the bottom components. For this, we use the helper function \verb|_distribute_part_ordered|, which is essentially a wrapper for partitions of a fixed length filled with zeros. At this point, we have added all zeros and will have to add (at least) double poles for the edges. Thus, if the orders on any vertex $v$ sum up to less than $2g_v$, we can move on to the next iteration.
\item[Step 6] We now distribute the poles $\bp$. This is again achieved via $2$-length partitions. Note that every $g=0$ vertex on top needs at least one pole (to compensate the edge(s) going down), so this gives an immediate check for the partitions.
\item[Step 7] Next, we distribute the poles among the bottom components. As poles are optional, we use \verb|_distribute_points|. At this point, we also save the difference of $2g_v-2$ and the orders distributed to $v$, i.e. the ``space'' left for half-edges. We test that this is at least $-2$ on each component, so that there is potential for at least one edge going up for every vertex.
\item[Step 8] \emph{Now we consider the top level for the first time.} We iterate through the number of top-level vertices, \verb|top_comp_len|, which is bounded by the sum of \verb|total_g_bot| and the number of poles on top-level. As we know the number of vertices and the distribution of genus, the Euler formula determines the number of edges, \verb|num_of_edges|. This gives a ``global'' check for the ``spaces'' left on the bottom components: they must sum up to $-2\cdot\verb|num_of_edges|$.
\item[Step 9] Similar to above, we now distribute genus, poles and zeros on top level: we use again \verb|_distribute_part_ordered| and \verb|_distribute_points|, as any zeros and poles are optional on top-level. We also run the obvious tests on the orders and record the spaces on top (there are much fewer constraints here, as the spaces may well be zero).
\item[Step 10] We now place the half-edges. We again start on bottom-level, because the poles give stronger constraints. Moreover, the half-edge orders on bottom determine those on top. The distribution of the half-edges into the spaces is accomplished by \verb|_place_legs_on_bot|, which uses partitions to recursively place the poles, and \verb|_place_legs_on_top|, which works similar but uses the orders distributed on bottom-level.
\item[Step 11] To assert that the graph we have created is connected, we create a Sage \verb|Graph| consisting of the edges we have placed. Note that every vertex is attached to an edge and we don't care for multi-edges when checking connectedness, so a very basic Sage \verb|Graph| suffices here.
\item[Step 12] The only thing missing are the marked points. These are distributed via \verb|_distribute_points| among all the vertices. We now check for stability.
\item[Step 13] In the final step, the legs are renumbered for consistency and the data is used to create a \verb|LevelGraph|. At this point, we check the GRC (via \verb|LevelGraph.is_legal()|). We also remove duplicates from the list of generated graphs.
\end{description}
\end{algorithm}

While we do not claim that algorithm is optimal, it certainly runs in a
reasonable time even for large strata:

\begin{lstlisting}
sage: from admcycles.diffstrata.bic import bic_alt_noiso
sage: %
CPU times: user 678 ms, sys: 9.38 ms, total: 688 ms
Wall time: 691 ms
384
\end{lstlisting}

\subsubsection*{Uniting {\upshape\texttt{LevelGraph}}s}

The BICs of a \verb|GeneralisedStratum| are the products of the BICs inside each component, subject to the residue conditions. Moreover, if there are several connected components, we also need to include the smooth stratum on each level.

Therefore, \verb|gen_bic| starts by running \verb|bic_alt_noiso| on each connected component and generating the \verb|dmp| for the embedding into the enveloping stratum. This is possible, because the BIC algorithm described above places the points on specific legs\footnote{This is only true for the \texttt{LevelGraph}s generated on a single component. In the final \texttt{EmbeddedLevelGraph}, the legs will be renumbered and this assumption is no longer valid!}: on each BIC, the $i$-th point of the signature is the point $i+1$.

Then, potentially the \verb|smooth_LG| of each component is added and we iterate over the product of these \verb|EmbeddedLevelGraph|s. 
Building the ``product graph'' is now simply a renumbering issue and accomplished by \verb|unite_embedded_graphs|.

The result is a list of valid \verb|EmbeddedLevelGraph|s (\verb|dmp| is now surjective), on which we can check the $\frakR$-GRC via \verb|is_legal|. Finally, we sort this list into isomorphism classes.

All BICs of a stratum are stored in the list \verb|bics|, which is generated automatically on first call:
\begin{lstlisting}
sage: P=GeneralisedStratum([Signature((0,0)),Signature((0,))])
sage: P.bics  # order might change
[EmbeddedLevelGraph(LG=LevelGraph([1, 0, 1],[[1], [2, 3, 4], [5]],[(1, 4)],{1: 0, 2: 0, 3: 0, 4: -2, 5: 0},[0, -1, 0],True),dmp={2: (0, 0), 3: (0, 1), 5: (1, 0)},dlevels={0: 0, -1: -1}),
 EmbeddedLevelGraph(LG=LevelGraph([1, 0, 1],[[1], [2, 3, 4], [5]],[(1, 4)],{1: 0, 2: 0, 3: 0, 4: -2, 5: 0},[0, -1, -1],True),dmp={2: (0, 0), 3: (0, 1), 5: (1, 0)},dlevels={0: 0, -1: -1}),
 EmbeddedLevelGraph(LG=LevelGraph([1, 1],[[1, 2], [3]],[],{1: 0, 2: 0, 3: 0},[0, -1],True),dmp={1: (0, 0), 2: (0, 1), 3: (1, 0)},dlevels={0: 0, -1: -1}),
 EmbeddedLevelGraph(LG=LevelGraph([1, 1],[[1, 2], [3]],[],{1: 0, 2: 0, 3: 0},[-1, 0],True),dmp={1: (0, 0), 2: (0, 1), 3: (1, 0)},dlevels={-1: -1, 0: 0})]
\end{lstlisting}

\subsection{Generating Profiles and Graphs}
\label{subsec:graphsandprofiles}
We now make \autoref{alg:profilesfrombics} more explicit.

\subsubsection*{Listing all Profiles}

We begin by constructing the \verb|lookup_list| of all profiles in a \verb|GeneralisedStratum|. The codimension $0$ and $1$ lists have already been treated: the first consists of the empty profile, the second of (the indices of) all BICs. 

The key tool for working with profiles are the maps \verb|top_to_bic| and \verb|bot_to_bic| described above. These are implemented as dictionary objects. 

More precisely, let \verb|X| be a\verb|GeneralisedStratum|. Then, for each index \verb|i| of \verb|X.bics|, \verb|X.DG.top_to_bic(i)| is a \verb|dict| mapping indices of \verb|X.bics[i].top.bics| to indices of \verb|X.bics|. Similarly, \verb|X.DG.bot_to_bic(i)| is a \verb|dict| mapping indices of \verb|X.bics[i].bot.bics| to indices of \verb|X.bics|.

\begin{example}\label{ex:toptobic}
We illustrate this in $\Omega\cM_2(2)$, cf. \autoref{H2} and note \autoref{numbering}:
\begin{lstlisting}
sage: X=Stratum((2,))
sage: X.bics  # order might differ
[EmbeddedLevelGraph(LG=LevelGraph([1, 0],[[1, 2], [3, 4, 5]],[(1, 4), (2, 5)],{1: 0, 2: 0, 3: 2, 4: -2, 5: -2},[0, -1],True),dmp={3: (0, 0)},dlevels={0: 0, -1: -1}),
 EmbeddedLevelGraph(LG=LevelGraph([1, 1],[[1], [2, 3]],[(1, 3)],{1: 0, 2: 2, 3: -2},[0, -1],True),dmp={2: (0, 0)},dlevels={0: 0, -1: -1})]
sage: X.DG.top_to_bic(0)
{0: 1}
sage: X.DG.bot_to_bic(0)
{}
sage: X.DG.top_to_bic(1)
{}
sage: X.DG.bot_to_bic(1)
{0: 0}
sage: X.lookup_list[0]
[()]
sage: X.lookup_list[1]
[(0,), (1,)]
sage: X.lookup_list[2]
[(1, 0)]
\end{lstlisting}
For completeness, we list the divisors of the top and bottom strata:
\begin{lstlisting}
sage: X.bics[0].top.bics
[EmbeddedLevelGraph(LG=LevelGraph([1, 0],[[1], [2, 3, 4]],[(1, 4)],{1: 0, 2: 0, 3: 0, 4: -2},[0, -1],True),dmp={2: (0, 0), 3: (0, 1)},dlevels={0: 0, -1: -1})]
sage: X.bics[0].bot.bics
[]
sage: X.bics[1].top.bics
[]
sage: X.bics[1].bot.bics
[EmbeddedLevelGraph(LG=LevelGraph([0, 0],[[1, 2, 3], [4, 5, 6]],[(2, 5), (3, 6)],{1: -2, 2: 0, 3: 0, 4: 2, 5: -2, 6: -2},[0, -1],True),dmp={1: (0, 1), 4: (0, 0)},dlevels={0: 0, -1: -1})]
\end{lstlisting}
Note that these are of course not graphs in \verb|X| (there are more marked points in \verb|dmp|!) but instead the \verb|top| and \verb|bot| strata of the BICs:
\begin{lstlisting}
sage: X.bics[0].top
LevelStratum(sig_list=[Signature((0, 0))],res_cond=[],leg_dict={1: (0, 0), 2: (0, 1)})
sage: X.bics[0].bot
LevelStratum(sig_list=[Signature((2, -2, -2))],res_cond=[[(0, 1), (0, 2)]],leg_dict={3: (0, 0), 4: (0, 1), 5: (0, 2)})
sage: X.bics[1].top
LevelStratum(sig_list=[Signature((0,))],res_cond=[],leg_dict={1: (0, 0)})
sage: X.bics[1].bot
LevelStratum(sig_list=[Signature((2, -2))],res_cond=[[(0, 1)]],leg_dict={2: (0, 0), 3: (0, 1)})
\end{lstlisting}
\end{example}

For examples of \verb|top_to_bic| and \verb|bot_to_bic| failing to be injective, see \autoref{reducible_degen} and \autoref{reducible_aut} below.

\begin{rem}
  Similar to \verb|top_to_bic| and \verb|bot_to_bic| there is also
  the function \verb|X.DG.middle_to_bic| describing degenerations of the middle levels of $3$-level graphs. In light of \autoref{leveltypes}, these give all ways of recursively extending a profile.

While \verb|middle_to_bic| is not needed for generating all \verb|EmbeddedLevelGraphs| it is essential for pulling back classes from a level, cf. \autoref{clutching}.
\end{rem}

With the help of \verb|top_to_bic| and \verb|bot_to_bic|, it is not difficult to construct the \verb|lookup_list| recursively:
\begin{algorithm}[Profile Generation]\mbox{}
\begin{description}
\item[Step 1] Construct \verb|lookup_list[0]| and \verb|lookup_list[1]|.
\item[Step 2] Constructing profiles of length $l+1\geq 2$: for every profile $p=(p_1,\dotsc,p_{l})$ of length $l$, add any profile of the form $(p_0,\dotsc,p_{l})$, where $p_0$ lies in the \verb|top_to_bic|-image of $p_1$.
If $l>1$, i.e. $p_1\neq p_{l}$, we also add the profiles $(p_1,\dotsc,p_{l+1})$, where $p_{l+1}$ lies in the \verb|bot_to_bic|-image of $p_{l}$.
\item[Step 3] Duplicates are removed and the \verb|list| constructed in the previous step is appended to \verb|lookup_list|.
\item[Step 4] Steps 2 and 3 are repeated until the list produced in Step 3 is empty.
\end{description}
\end{algorithm}

\begin{rem}
A few notes on the validity of the above algorithm:
\begin{enumerate}
\item The algorithm terminates, as there are no graphs with more levels than the dimension of the stratum. Indeed, each level increases the codimension of the corresponding boundary contribution by $1$ and thus also the dimension of the levels must eventually decrease and they will stop having BICs.
\item The entries $p_i$ of a profile are distinct, cf. \cite[\S 5]{strataEC}. Moreover, a BIC can appear only as a top \emph{or} a bottom degeneration of another BIC. Hence $p_0$ and $p_l$ cannot have been contained in the previous profile.
\item In the case $l=1$, the bottom-level degenerations $(p_1,p_2)$ of $p_1$ are also obtained as top-level degenerations of $p_2$.
\end{enumerate}
\end{rem}

All that is needed in the above construction are the maps \verb|top_to_bic| and \verb|bot_to_bic|, which are constructed as follows (for \verb|top_to_bic(i)|, \verb|bot_to_bic(i)| uses the same algorithm, with the obvious changes).
\begin{algorithm}[Construction of \texttt{top\char`_to\char`_bic} and \texttt{bot\char`_to\char`_bic}]\mbox{}
\begin{description}
\item[Step 1] Loop through all BICs \verb|X.bics[i].top|. 

We denote \verb|X.bics[i].top.bics[j]| by \verb|Bt|.
\item[Step 2] Clutch \verb|Bt| to \verb|X.bics[i].bot| to obtain a ($3$-level) graph \verb|G|.
\item[Step 3] Squish the bottom level of \verb|G| to obtain the BIC \verb|G.delta(1)|.
\item[Step 4] Go through the list \verb|X.bics| and check which of these is isomorphic to \verb|G.delta(1)|. The index of this graph (in \verb|X.bics|) is stored as the image of \verb|j| in \verb|top_to_bic(i)|.
\end{description}
\end{algorithm}

\begin{rem}
The implementation includes a small optimisation: Using the fact that a BIC can only appear as a top- \emph{or} a bottom-level degeneration of another BIC, we remove the BICs used for \verb|top_to_bic(i)| from the candidates for \verb|bot_to_bic(i)| and vice versa.
\end{rem}

Two operations are used in the above algorithm: Clutching and \verb|delta|. Clutching and splitting graphs is a delicate issue and is described in more detail in \autoref{clutching}.

We briefly describe the process of squishing graphs and the implementation of \verb|delta|.

While we may apply \verb|delta| to an \verb|EmbeddedLevelGraph|, the actual work happens on the underlying \verb|LevelGraph| (the \verb|dmp| is not changed by applying \verb|delta|).

On the \verb|LevelGraph|, \verb|delta| is implemented by consecutively contracting all levels, except for one. Because the graph has fewer levels on each iteration, it is easier to keep track of the numbering when iterating through the levels starting at bottom level.

\begin{rem}
Every \verb|LevelGraph| object should be considered \emph{immutable}. In particular, all these operations \emph{do not} change the \verb|LevelGraph|, but instead always return a new object.
\end{rem}

Squishing a single level is achieved by \verb|squish_vertical(i)|. 

While level-crossings can be very non-trivial, squishing horizontal edges is straight-forward. There are only two possibilities:
\begin{enumerate}
\item The edge is a loop. In this case, remove the edge and increase the genus of the vertex by one.
\item The edge connects two distinct vertices $v_1$ and $v_2$. In this case, remove the edge, add the genus of $v_2$ to the genus of $v_1$ and all the legs of $v_2$ to the legs of $v_1$ and remove $v_2$.
\end{enumerate}

We can fall back to this for vertical squishing:
\begin{algorithm}[Vertical Squishing]\mbox{}
\begin{description}
\item[Step 1] Determine the set of vertices \verb|vv| on the next lower level, as well as the edges, \verb|ee|, passing from level \verb|i| to the next lower level.
\item[Step 2] Make a copy of the information needed to generate the \verb|LevelGraph| (i.e. genera, legs, edges, pole orders and levels).
\item[Step 3] Raise the level of every vertex in \verb|vv| to \verb|i| in the copied data.
\item[Step 4] Create a new \verb|LevelGraph| from this data.
\item[Step 5] Squish each of the (now horizontal!) edges in \verb|ee| on this graph.
\end{description}
\end{algorithm}

\begin{rem}
This is not very efficient, because a new graph has to be created for every edge crossing the level. However, this implementation can still be found in the method \verb|LevelGraph.squish_vertical_slow|. The actual \verb|squish_vertical| avoids this, by reimplementing the book-keeping that is avoided by the recursive graph creation, creates the ``correct'' data in one go and thus creates only one new \verb|LevelGraph|.
\end{rem}

\begin{example}
We briefly illustrate squishing a $4$-level graph and how this relates to the numbering of the profile. Note that we use \verb|any| so that we make no assumptions about the irreducibility of the profiles.
\begin{lstlisting}
sage: X=GeneralisedStratum([Signature((4,))])
sage: p = X.enhanced_profiles_of_length(4)[0][0]
sage: g = X.lookup_graph(p)
sage: assert any(g.squish_vertical(0).is_isomorphic(G) 
									for G in X.lookup(p[1:]))
sage: assert any(g.squish_vertical(1).is_isomorphic(G) 
									for G in X.lookup(p[:1]+p[2:]))
sage: assert any(g.squish_vertical(2).is_isomorphic(G) 
									for G in X.lookup(p[:2]+p[3:]))
sage: assert any(g.squish_vertical(3).is_isomorphic(G) 
									for G in X.lookup(p[:3]))
\end{lstlisting}
\end{example}

\begin{rem}
There are some subtleties regarding the level numbering.
\begin{enumerate}
\item The arguments of \verb|delta| and \verb|squish_vertical| are shifted by one, that is \verb|delta(i)| squishes all level passages except the one from $i-1$ to $i$, while \verb|squish_vertical(i)| squishes the level passage from $i$ to $i+1$.
\item Note that the \verb|squish_vertical| method exists both for \verb|LevelGraph|s and \verb|EmbeddedLevelGraph|s. Of course, the \verb|EmbeddedLevelGraph| version should be used (the level dictionary needs to be adapted) and the argument is the \emph{relative} level number (while the \verb|LevelGraph| version requires the internal level number).
\end{enumerate}
\end{rem}

\subsubsection*{Building Graphs from Profiles} 

Now that we have generated all possible profiles, the question remains how to build the graph(s) associated to a profile. This is achieved by recursive clutching. The idea is to pick a BIC $B$ from the profile and sort the remainder into profiles in $B^\top$ an $B^\bot$, build the graphs there and then clutch these together.

More precisely, the \verb|lookup| method of a \verb|GeneralisedStratum| \verb|X| is implemented, for a profile $(p_0,\dotsc,p_l)$, as follows:
\begin{algorithm}[Graph lookup]\mbox{}
\begin{description}
\item[Step 1] Denote by \verb|B| the BIC with index $p_0$ in \verb|X.bics|.
\item[Step 2] Create a new nested \verb|list| \verb|bot_lists|.
\item[Step 3] For every \verb|j| in the (remaining) profile $(p_1,\dotsc,p_l)$, find all occurrences of \verb|j| in the values of \verb|bot_to_bic|.
\item[Step 4] For each occurrence of \verb|j| in the values of \verb|bot_to_bic|, append the corresponding key to (a copy of) each \verb|list| in \verb|bot_lists|.
\item[Step 5] Replace \verb|bot_lists| by this nested list and continue the loop at Step 3.
\item[Step 6] Apply \verb|B.bot.lookup| to each profile in \verb|bot_lists| to obtain a \verb|list| of \verb|EmbeddedLevelGraph|s in \verb|B.bot|.
\item[Step 7] Clutch each of these graphs to the graph \verb|B.top| to obtain a \verb|list| of \verb|EmbeddedLevelGraph|s in \verb|X|.
\item[Step 8] Return this \verb|list|, sorted by isomorphism classes.
\end{description}
\end{algorithm}
As mentioned above, clutching and splitting are subtle issues and will be addressed in detail in \autoref{clutching}.

\begin{rem}
The ``branching out'' in Step 4 can cause reducible profiles (if \verb|bot_to_bic| is not injective and the clutching also results in non-isomorphic graphs in \verb|X|).

The ``inverse'' dictionaries of \verb|top_to_bic| and \verb|bot_to_bic| needed in Step 3 are also readily available: Revisiting \autoref{ex:toptobic}, we see that
\begin{lstlisting}
sage: X.DG.top_to_bic(0)
{0: 1}
sage: X.DG.top_to_bic_inv(0)
{1: [0]}
sage: X.DG.bot_to_bic(0)
{}
sage: X.DG.bot_to_bic_inv(0)
{}
\end{lstlisting}
Obviously, the values of \verb|top_to_bic_inv| and \verb|bot_to_bic_inv| have to be \verb|list|s, as the dictionaries are not necessarily injective.
\end{rem}

Now, \verb|lookup_graph| simply returns elements of the \verb|list| generated by \verb|lookup| (by default the first).

\begin{rem}
Note that the \verb|lookup| algorithm actually checks both \verb|top_to_bic| and \verb|bot_to_bic| and thus works for an arbitrary permutation of the profile (cf. \cite[Prop. 5.1]{strataEC}).

As the ordering given by \verb|delta| is important at many other places, for useful comparison and caching, it makes sense to restrict to the ``ordered'' profile.

For example (note that the profile values are ``arbitrary'', cf. \autoref{numbering}):
\begin{lstlisting}
sage: X=Stratum((4,))
sage: X.lookup_list[2][0]
(2, 7)
sage: X.lookup((2,7))
[EmbeddedLevelGraph(LG=LevelGraph([1, 0, 0],[[1], [2, 3, 4, 5], [6, 7, 8, 9]],[(1, 5), (2, 7), (3, 8), (4, 9)],{1: 0, 2: 0, 3: 0, 4: 0, 5: -2, 6: 4, 7: -2, 8: -2, 9: -2},[0, -1, -2],True),dmp={6: (0, 0)},dlevels={0: 0, -1: -1, -2: -2})]
sage: X.lookup((7,2))
[EmbeddedLevelGraph(LG=LevelGraph([1, 0, 0],[[1], [2, 3, 4, 5], [6, 7, 8, 9]],[(3, 7), (4, 8), (5, 9), (1, 2)],{1: 0, 2: -2, 3: 0, 4: 0, 5: 0, 6: 4, 7: -2, 8: -2, 9: -2},[0, -1, -2],True),dmp={6: (0, 0)},dlevels={0: 0, -1: -1, -2: -2})]
sage: X.lookup((7,2))[0].is_isomorphic(X.lookup((2,7))[0])
True
\end{lstlisting}
However, \verb|lookup_graph| does not permit this and requires the profile to be ordered:
\begin{lstlisting}
sage: X.lookup_graph((2,7))
EmbeddedLevelGraph(LG=LevelGraph([1, 0, 0],[[1], [2, 3, 4, 5], [6, 7, 8, 9]],[(1, 5), (2, 7), (3, 8), (4, 9)],{1: 0, 2: 0, 3: 0, 4: 0, 5: -2, 6: 4, 7: -2, 8: -2, 9: -2},[0, -1, -2],True),dmp={6: (0, 0)},dlevels={0: 0, -1: -1, -2: -2})
sage: X.lookup_graph((7,2))
\end{lstlisting}
\end{rem}

\subsection{Checking Degenerations}
\label{subsec:degenerations}

Clearly, for a graph $\Gamma'$ to be a degeneration of $\Gamma$, it is necessary for the profile of $\Gamma$ to be an (ordered) subset of the profile of $\Gamma'$. However, if $\Gamma$ is reducible, this criterion is no longer sufficient.

To check whether a graph associated to an enhanced profile is the degeneration of the graph associated to another enhanced profile, \verb|diffstrata| provides the method \verb|is_degeneration|:
\begin{lstlisting}
sage: X=Stratum((2,))
sage: X.is_degeneration(((1,), 0), ((), 0))
True
sage: X.is_degeneration(((1,), 0), ((0,), 0))
False
sage: X.is_degeneration(((), 0), ((0,), 0))
False
\end{lstlisting}
Note that the arguments are enhanced profiles, i.e. \verb|tuple|s of \verb|tuple|s.

Before describing the implementation, we continue to examine some of the phenomena in $\Omega\cM_3(4)$, the minimal stratum in genus $3$, 
continuing \autoref{reducibleprofile}.

Keeping \autoref{numbering} in mind, we nonetheless work with concrete enhanced profiles to ease readability.

\begin{example}\label{reducible_degen}
We revisit the situation of \autoref{reducibleprofile}. It is not difficult to find the BICs $\delta_1(G)$ and $\delta_2(G)$ in the \verb|list| of (8) BICs in $\Omega\cM_3(4)$. In our case, we see that $\delta_1(G)$ corresponds to BIC \verb|1| and $\delta_2(G)$ corresponds to BIC \verb|6| (note the values of the prongs!).
\begin{lstlisting}
sage: X=Stratum((4,))
sage: len(X.bics)
8
sage: X.bics[1].explain()
LevelGraph embedded into stratum Stratum: (4,)
with residue conditions: []
 with:
On level 0:
* A vertex (number 0) of genus 1
On level 1:
* A vertex (number 1) of genus 1
The marked points are on level 1.
More precisely, we have:
* Marked point (0, 0) of order 4 on vertex 1 on level 1
Finally, we have 2 edges. More precisely:
* 2 edges between vertex 0 (on level 0) and vertex 1 (on level 1) with prongs 1 and 1.
sage: X.bics[6].explain()
LevelGraph embedded into stratum Stratum: (4,)
with residue conditions: []
 with:
On level 0:
* A vertex (number 0) of genus 2
On level 1:
* A vertex (number 1) of genus 0
The marked points are on level 1.
More precisely, we have:
* Marked point (0, 0) of order 4 on vertex 1 on level 1
Finally, we have 2 edges. More precisely:
* 2 edges between vertex 0 (on level 0) and vertex 1 (on level 1) with prongs 3 and 1.
\end{lstlisting}
Consequently, the profile \verb|(1, 6)| is reducible and contains the graphs $G$ and $H$. Here, $H$ has enhanced profile \verb|((1, 6), 0)| and $G$ is \verb|((1, 6), 1)| (note the genera!):
\begin{lstlisting}
sage: X.lookup((1,6))
[EmbeddedLevelGraph(LG=LevelGraph([1, 1, 0],[[1, 2], [3, 4], [5, 6, 7]],[(2, 4), (3, 6), (1, 7)],{1: 0, 2: 0, 3: 2, 4: -2, 5: 4, 6: -4, 7: -2},[0, -1, -2],True),dmp={5: (0, 0)},dlevels={0: 0, -1: -1, -2: -2}),
 EmbeddedLevelGraph(LG=LevelGraph([1, 0, 0],[[1, 2], [3, 4, 5, 6], [7, 8, 9]],[(1, 5), (2, 6), (3, 8), (4, 9)],{1: 0, 2: 0, 3: 2, 4: 0, 5: -2, 6: -2, 7: 4, 8: -4, 9: -2},[0, -1, -2],True),dmp={7: (0, 0)},dlevels={0: 0, -1: -1, -2: -2})]
\end{lstlisting}
We see the reducibility reflected in \verb|top_to_bic| and \verb|bot_to_bic| (illustrated best by their inverse dictionaries):
\begin{lstlisting}
sage: X.DG.bot_to_bic_inv(1)
{3: [0, 2], 6: [1, 4, 7], 0: [3], 4: [5], 2: [6]}
sage: X.DG.top_to_bic_inv(6)
{1: [0, 4], 4: [1], 5: [2], 3: [3]}
\end{lstlisting}
The four-level graphs degenerating from this profile can be seen in \autoref{fourlevelred}.
\begin{figure}
\[\begin{tikzpicture}[
		baseline={([yshift=-.5ex]current bounding box.center)},
		scale=2,very thick,
		bend angle=30,
		every loop/.style={very thick},
     		comp/.style={circle,fill,black,,inner sep=0pt,minimum size=5pt},
		prong left/.style={pos=.5,left,font=\small},
		prong right/.style={pos=.5,right,font=\small},
		circled number/.style={circle, draw, inner sep=0pt, minimum size=12pt},
		bottom right with distance/.style={below right,text height=10pt}]
\node[circled number] (VT) {$1$};
\node[comp] (T) [below=of VT] {}
	edge node [prong right] {$1$} (VT);
\node[comp] (B1) [below=of T] {}
	edge [bend left] node [prong left] {$1$} (T)
	edge [bend right] node [prong right] {$1$} (T);
\node[comp] (B2) [below=of B1] {}
	edge [bend left] node [prong left] {$1$} (B1)
	edge [bend right] node [prong right] {$3$} (B1);
\node [bottom right with distance] (B2-4) at (B2.south east) {$4$};
\path (B2) edge [shorten >=4pt] (B2-4.center);
\end{tikzpicture}
\quad
\begin{tikzpicture}[
		baseline={([yshift=-.5ex]current bounding box.center)},
		scale=2,very thick,
		bend angle=30,
		every loop/.style={very thick},
     		comp/.style={circle,fill,black,,inner sep=0pt,minimum size=5pt},
		prong left/.style={pos=.5,left,font=\small},
		prong right/.style={pos=.5,right,font=\small},
		circled number/.style={circle, draw, inner sep=0pt, minimum size=12pt},
		bottom right with distance/.style={below right,text height=10pt}]
\node[circled number] (VT) {$1$};
\node[comp] (T) [below=of VT] {}
	edge node [prong right] {$1$} (VT);
\node[circled number] (B1) [below=of T] {$1$}
	edge node [prong right] {$1$} (T);
\node[comp] (B2) [below=of B1] {}
	edge [bend left] node [prong left] {$1$} (T)
	edge node [prong right] {$3$} (B1);
\node [bottom right with distance] (B2-4) at (B2.south east) {$4$};
\path (B2) edge [shorten >=4pt] (B2-4.center);
\end{tikzpicture}
\quad
\begin{tikzpicture}[
		baseline={([yshift=-.5ex]current bounding box.center)},
		scale=2,very thick,
		bend angle=30,
		every loop/.style={very thick},
     		comp/.style={circle,fill,black,,inner sep=0pt,minimum size=5pt},
		prong left/.style={pos=.5,left,font=\small},
		prong right/.style={pos=.5,right,font=\small},
		circled number/.style={circle, draw, inner sep=0pt, minimum size=12pt},
		bottom right with distance/.style={below right,text height=10pt}]
\node[circled number] (T) {$1$};
\node[comp] (B1) [below=of T] {}
	edge [bend left] node [prong left] {$1$} (T)
	edge [bend right] node [prong right] {$1$} (T);
\node[comp] (M) [below=of B1] {}
	edge node [prong right] {$1$} (B1);
\node[comp] (B2) [below=of M] {}
	edge [bend left] node [prong left] {$1$} (M)
	edge [bend right] node [prong right] {$3$} (M);
\node [bottom right with distance] (B2-4) at (B2.south east) {$4$};
\path (B2) edge [shorten >=4pt] (B2-4.center);
\end{tikzpicture}
\quad
\begin{tikzpicture}[
		baseline={([yshift=-.5ex]current bounding box.center)},
		scale=2,very thick,
		bend angle=30,
		every loop/.style={very thick},
     		comp/.style={circle,fill,black,,inner sep=0pt,minimum size=5pt},
		prong left/.style={pos=.5,left,font=\small},
		prong right/.style={pos=.5,right,font=\small},
		circled number/.style={circle, draw, inner sep=0pt, minimum size=12pt},
		bottom right with distance/.style={below right,text height=10pt}]
\node[circled number] (T) {$1$};
\node[comp] (B1) [below=of T] {}
	edge node [prong right] {$1$} (T);
\node[comp] (M) [below=of B1] {}
	edge [bend left] node [prong left] {$1$} (B1)
	edge [bend right] node [prong right] {$1$} (B1);
\node[comp] (B2) [below=of M] {}
	edge [bend left] node [prong left] {$1$} (T)
	edge node [prong right] {$3$} (M);
\node [bottom right with distance] (B2-4) at (B2.south east) {$4$};
\path (B2) edge [shorten >=4pt] (B2-4.center);
\end{tikzpicture}
\quad
\begin{tikzpicture}[
		baseline={([yshift=-.5ex]current bounding box.center)},
		scale=2,very thick,
		bend angle=30,
		every loop/.style={very thick},
     		comp/.style={circle,fill,black,,inner sep=0pt,minimum size=5pt},
		prong left/.style={pos=.5,left,font=\small},
		prong right/.style={pos=.5,right,font=\small},
		circled number/.style={circle, draw, inner sep=0pt, minimum size=12pt},
		bottom right with distance/.style={below right,text height=10pt}]
\node[circled number] (T) {$1$};
\node[comp] (B1) [below=of T] {}
	edge node [prong right] {$1$} (T);
\node[comp] (M) [below=of B1] {}
	edge [bend left] node [prong left] {$1$} (T)
	edge node [prong left, left=-3pt] {$1$} (B1);
\node[comp] (B2) [below=of M] {}
	edge [bend right] node [prong right] {$1$} (B1)
	edge node [prong left] {$3$} (M);
\node [bottom right with distance] (B2-4) at (B2.south east) {$4$};
\path (B2) edge [shorten >=4pt] (B2-4.center);
\end{tikzpicture}
\]
\caption{The four-level graphs degenerating from the \emph{reducible} profile \texttt{(1, 6)}. From left to right: \texttt{((5,1, 6), 0)}, \texttt{((5, 1, 6), 1)}, \texttt{((1, 4, 6), 0)}, \texttt{((1, 3, 6), 0)} and \texttt{((1, 3, 6), 1)}. Note that \texttt{(1, 4, 6)} is \emph{irreducible}.}
\label{fourlevelred}
\end{figure}
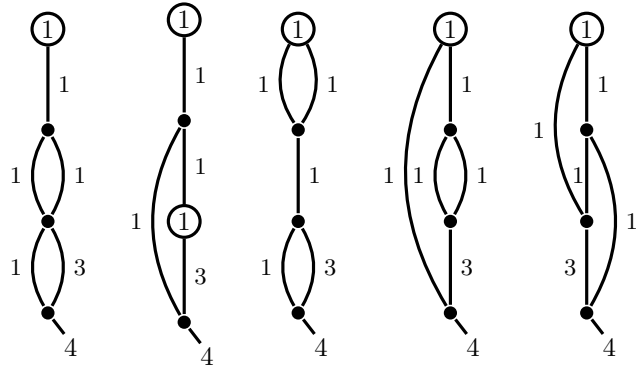
Applying a compact-type degeneration on top level, the profile remains reducible:
\begin{lstlisting}
sage: X.lookup((5,1,6))
[EmbeddedLevelGraph(LG=LevelGraph([1, 0, 1, 0],[[1], [2, 3, 4], [5, 6], [7, 8, 9]],[(1, 4), (3, 6), (5, 8), (2, 9)],{1: 0, 2: 0, 3: 0, 4: -2, 5: 2, 6: -2, 7: 4, 8: -4, 9: -2},[0, -1, -2, -3],True),dmp={7: (0, 0)},dlevels={0: 0, -1: -1, -2: -2, -3: -3}),
 EmbeddedLevelGraph(LG=LevelGraph([1, 0, 0, 0],[[1], [2, 3, 4], [5, 6, 7, 8], [9, 10, 11]],[(1, 4), (2, 7), (3, 8), (5, 10), (6, 11)],{1: 0, 2: 0, 3: 0, 4: -2, 5: 2, 6: 0, 7: -2, 8: -2, 9: 4, 10: -4, 11: -2},[0, -1, -2, -3],True),dmp={9: (0, 0)},dlevels={0: 0, -1: -1, -2: -2, -3: -3})]
\end{lstlisting}
We thus see the necessity of considering the \emph{enhanced} profile when considering degenerations:
\begin{lstlisting}
sage: X.is_degeneration(((5,1,6),0), ((1,6),0))
True
sage: X.is_degeneration(((5,1,6),0), ((1,6),1))
False
sage: X.is_degeneration(((5,1,6),1), ((1,6),0))
False
sage: X.is_degeneration(((5,1,6),1), ((1,6),1))
True
\end{lstlisting}
Also, $G$ admits a compact-type degeneration on middle level that is not possible in $H$ for stability reasons. The profile \verb|(1, 4, 6)| is irreducible:
\begin{lstlisting}
sage: X.lookup((1,4,6))
[EmbeddedLevelGraph(LG=LevelGraph([1, 0, 0, 0],[[1, 2], [3, 4, 5], [6, 7, 8], [9, 10, 11]],[(1, 4), (2, 5), (3, 8), (6, 10), (7, 11)],{1: 0, 2: 0, 3: 2, 4: -2, 5: -2, 6: 2, 7: 0, 8: -4, 9: 4, 10: -4, 11: -2},[0, -1, -2, -3],True),dmp={9: (0, 0)},dlevels={0: 0, -1: -1, -2: -2, -3: -3})]
sage: X.is_degeneration(((1,4,6),0), ((1,6),0))
False
sage: X.is_degeneration(((1,4,6),0), ((1,6),1))
True
\end{lstlisting}
The final two degenerations come from intersecting with BIC \verb|3|, the ``triple banana'' (cf. \autoref{triplebanana} and \autoref{fig:triplebanana}). Again, the profile remains reducible but note that both graphs have long edges now.
\begin{lstlisting}
sage: X.lookup((1,3,6))
[EmbeddedLevelGraph(LG=LevelGraph([1, 0, 0, 0],[[1, 2], [3, 4, 5], [6, 7, 8], [9, 10, 11]],[(2, 5), (3, 7), (4, 8), (6, 10), (1, 11)],{1: 0, 2: 0, 3: 0, 4: 0, 5: -2, 6: 2, 7: -2, 8: -2, 9: 4, 10: -4, 11: -2},[0, -1, -2, -3],True),dmp={9: (0, 0)},dlevels={0: 0, -1: -1, -2: -2, -3: -3}),
 EmbeddedLevelGraph(LG=LevelGraph([1, 0, 0, 0],[[1, 2], [3, 4, 5], [6, 7, 8], [9, 10, 11]],[(2, 5), (1, 7), (4, 8), (6, 10), (3, 11)],{1: 0, 2: 0, 3: 0, 4: 0, 5: -2, 6: 2, 7: -2, 8: -2, 9: 4, 10: -4, 11: -2},[0, -1, -2, -3],True),dmp={9: (0, 0)},dlevels={0: 0, -1: -1, -2: -2, -3: -3})]
sage: X.is_degeneration(((1,3,6),0), ((1,6),0))
True
sage: X.is_degeneration(((1,3,6),1), ((1,6),0))
False
sage: X.is_degeneration(((1,3,6),1), ((1,6),1))
True
sage: X.is_degeneration(((1,3,6),0), ((1,6),1))
False
\end{lstlisting}
Note that the two graphs in the profile \verb|(1, 3, 6)| cannot be distinguished by their levels.
\end{example}

\begin{rem}
We summarise the following observations from \autoref{reducible_degen}.
\begin{enumerate}
\item Extending an irreducible profile can make it reducible.
\item Extending a reducible profile can make it irreducible.
\item A reducible profile implies a non-injectivity of \verb|top_to_bic| and \verb|bot_to_bic|.
\end{enumerate}
\end{rem}

The converse of the last statement is not true in general.

\begin{example}\label{triplebanana}
We continue in the notation of \autoref{reducible_degen}. Consider the ``triple banana'', the BIC with a non-trivial $S_3$ action (cf. \autoref{fig:triplebanana}).
It is again not difficult to find in the \verb|list| of BICs, here it is number \verb|3|:
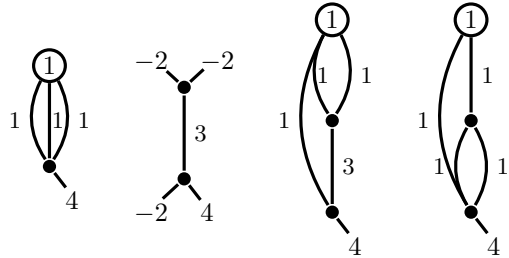
\begin{figure}
\[\begin{tikzpicture}[
		baseline={([yshift=-.5ex]current bounding box.center)},
		scale=2,very thick,
		bend angle=30,
		every loop/.style={very thick},
     		comp/.style={circle,fill,black,,inner sep=0pt,minimum size=5pt},
		prong left/.style={pos=.5,left,font=\small},
		prong right/.style={pos=.5,right,font=\small},
		circled number/.style={circle, draw, inner sep=0pt, minimum size=12pt},
		bottom right with distance/.style={below right,text height=10pt}]
\node[circled number] (T) {$1$};
\node[comp] (B) [below=of T] {}
	edge [bend left] node [prong left] {$1$} (T)
	edge node [prong right, right=-3pt] {$1$} (T)
	edge [bend right] node [prong right] {$1$} (T);
\node [bottom right with distance] (B-4) at (B.south east) {$4$};
\path (B) edge [shorten >=4pt] (B-4.center);
\end{tikzpicture}
\quad
\begin{tikzpicture}[
		baseline={([yshift=-.5ex]current bounding box.center)},
		scale=2,very thick,
		bend angle=30,
		every loop/.style={very thick},
     		comp/.style={circle,fill,black,,inner sep=0pt,minimum size=5pt},
		prong left/.style={pos=.5,left,font=\small},
		prong right/.style={pos=.5,right,font=\small},
		circled number/.style={circle, draw, inner sep=0pt, minimum size=12pt},
		bottom right with distance/.style={below right,text height=10pt},
		bottom left with distance/.style={below left,text height=10pt},
		above right with distance/.style={above right,text height=10pt},
		above left with distance/.style={above left,text height=10pt}]
\node[comp] (T) {};
\node[comp] (B) [below=of T] {}
	edge node [prong right] {$3$} (T);
\node [bottom right with distance] (B-4) at (B.south east) {$4$};
\path (B) edge [shorten >=4pt] (B-4.center);
\node [bottom left with distance] (B-2) at (B.south west) {$-2$};
\path (B) edge [shorten >=6pt] (B-2.center);
\node [above right with distance] (T-21) at (T.north east) {$-2$};
\path (T) edge [shorten >=7pt] (T-21.center);
\node [above left with distance] (T-22) at (T.north west) {$-2$};
\path (T) edge [shorten >=8pt] (T-22.center);
\end{tikzpicture}
\quad
\begin{tikzpicture}[
		baseline={([yshift=-.5ex]current bounding box.center)},
		scale=2,very thick,
		bend angle=30,
		every loop/.style={very thick},
     		comp/.style={circle,fill,black,,inner sep=0pt,minimum size=5pt},
		prong left/.style={pos=.5,left,font=\small},
		prong right/.style={pos=.5,right,font=\small},
		circled number/.style={circle, draw, inner sep=0pt, minimum size=12pt},
		bottom right with distance/.style={below right,text height=10pt}]
\node[circled number] (T) {$1$};
\node[comp] (B1) [below=of T] {}
	edge [bend left] node [prong right, right=-3pt] {$1$} (T)
	edge [bend right] node [prong right] {$1$} (T);
\node[comp] (B2) [below=of B1] {}
	edge [bend left] node [prong left] {$1$} (T)
	edge node [prong right] {$3$} (B1);
\node [bottom right with distance] (B2-4) at (B2.south east) {$4$};
\path (B2) edge [shorten >=4pt] (B2-4.center);
\end{tikzpicture}
\quad
\begin{tikzpicture}[
		baseline={([yshift=-.5ex]current bounding box.center)},
		scale=2,very thick,
		bend angle=30,
		every loop/.style={very thick},
     		comp/.style={circle,fill,black,,inner sep=0pt,minimum size=5pt},
		prong left/.style={pos=.5,left,font=\small},
		prong right/.style={pos=.5,right,font=\small},
		circled number/.style={circle, draw, inner sep=0pt, minimum size=12pt},
		bottom right with distance/.style={below right,text height=10pt}]
\node[circled number] (T) {$1$};
\node[comp] (B1) [below=of T] {}
	edge node [prong right] {$1$} (T);
\node[comp] (B2) [below=of B1] {}
	edge [bend left] node [prong left] {$1$} (T)
	edge [bend left] node [prong left] {$1$} (B1)
	edge [bend right] node [prong right] {$1$} (B1);
\node [bottom right with distance] (B2-4) at (B2.south east) {$4$};
\path (B2) edge [shorten >=4pt] (B2-4.center);
\end{tikzpicture}
\]
\caption{From left to right: the ``triple banana'', BIC \texttt{3}; the compact type BIC inside $\Omega\cM_0(4, -2, -2, -2)$, the bottom level of the ``triple banana''; the unique graph in the profile \texttt{(3, 6)}; the unique graph in the profile \texttt{(1, 3)}.}
\label{fig:triplebanana}
\end{figure}
\begin{lstlisting}
sage: X.bics[3].explain()
LevelGraph embedded into stratum Stratum: (4,)
with residue conditions: []
 with:
On level 0:
* A vertex (number 0) of genus 1
On level 1:
* A vertex (number 1) of genus 0
The marked points are on level 1.
More precisely, we have:
* Marked point (0, 0) of order 4 on vertex 1 on level 1
Finally, we have 3 edges. More precisely:
* 3 edges between vertex 0 (on level 0) and vertex 1 (on level 1) with prongs 1, 1 and 1.
sage: len(X.bics[3].automorphisms)
6
\end{lstlisting}
The bottom level is the stratum $\Omega\cM_0(4, -2, -2, -2)$, which is one-dimensional and contains three BICs of compact type, distinguished only by the numbering of their marked points (\verb|dmp|!):
\begin{lstlisting}
sage: B=X.bics[3].bot; print(B)
Stratum: Signature((4, -2, -2, -2))
with residue conditions: [(0, 1), (0, 2), (0, 3)]
dimension: 1
leg dictionary: {4: (0, 0), 5: (0, 1), 6: (0, 2), 7: (0, 3)}
leg orbits: [[(0, 0)], [(0, 1), (0, 3), (0, 2)]]

sage: B.bics
[EmbeddedLevelGraph(LG=LevelGraph([0, 0],[[1, 2, 3], [4, 5, 6]],[(3, 6)],{1: -2, 2: -2, 3: 2, 4: 4, 5: -2, 6: -4},[0, -1],True),dmp={1: (0, 1), 2: (0, 3), 4: (0, 0), 5: (0, 2)},dlevels={0: 0, -1: -1}),
 EmbeddedLevelGraph(LG=LevelGraph([0, 0],[[1, 2, 3], [4, 5, 6]],[(3, 6)],{1: -2, 2: -2, 3: 2, 4: 4, 5: -2, 6: -4},[0, -1],True),dmp={1: (0, 2), 2: (0, 3), 4: (0, 0), 5: (0, 1)},dlevels={0: 0, -1: -1}),
 EmbeddedLevelGraph(LG=LevelGraph([0, 0],[[1, 2, 3], [4, 5, 6]],[(3, 6)],{1: -2, 2: -2, 3: 2, 4: 4, 5: -2, 6: -4},[0, -1],True),dmp={1: (0, 1), 2: (0, 2), 4: (0, 0), 5: (0, 3)},dlevels={0: 0, -1: -1})]
\end{lstlisting}
Checking \verb|bot_to_bic|, one might expect the profile \verb|(3, 6)| to be reducible, but it turns out that, no matter how we glue the compact type graph into the bottom level, the resulting graphs are always isomorphic:
\begin{lstlisting}
sage: X.DG.bot_to_bic_inv(3)
{6: [0, 1, 2]}
sage: X.lookup((3,6))
[EmbeddedLevelGraph(LG=LevelGraph([1, 0, 0],[[1, 2, 3], [4, 5, 6], [7, 8, 9]],[(2, 5), (3, 6), (4, 8), (1, 9)],{1: 0, 2: 0, 3: 0, 4: 2, 5: -2, 6: -2, 7: 4, 8: -4, 9: -2},[0, -1, -2],True),dmp={7: (0, 0)},dlevels={0: 0, -1: -1, -2: -2})]
\end{lstlisting}
Recall, however, that after intersecting with BIC \verb|1|, the profile \verb|(1, 3, 6)| \emph{is} reducible! In other words, gluing the three BICs of \verb|B| into the bottom level of \verb|(1, 3)| yields two non-isomorphic graphs, cf. \autoref{reducible_degen} and \autoref{fourlevelred}.
\end{example}

\begin{example}
The last point of \autoref{triplebanana} is worth emphasizing. Using the notation of above, we observe the following situation: the profile \verb|(1, 3)| is irreducible, the profile \verb|(3, 6)| is irreducible, but the profile \verb|(1, 3, 6)| is reducible!
\begin{lstlisting}
sage: X.lookup((1,3))
[EmbeddedLevelGraph(LG=LevelGraph([1, 0, 0],[[1, 2], [3, 4, 5], [6, 7, 8, 9]],[(2, 5), (1, 7), (3, 8), (4, 9)],{1: 0, 2: 0, 3: 0, 4: 0, 5: -2, 6: 4, 7: -2, 8: -2, 9: -2},[0, -1, -2],True),dmp={6: (0, 0)},dlevels={0: 0, -1: -1, -2: -2})]
sage: X.lookup((3,6))
[EmbeddedLevelGraph(LG=LevelGraph([1, 0, 0],[[1, 2, 3], [4, 5, 6], [7, 8, 9]],[(2, 5), (3, 6), (4, 8), (1, 9)],{1: 0, 2: 0, 3: 0, 4: 2, 5: -2, 6: -2, 7: 4, 8: -4, 9: -2},[0, -1, -2],True),dmp={7: (0, 0)},dlevels={0: 0, -1: -1, -2: -2})]
sage: X.lookup((1,3,6))
[EmbeddedLevelGraph(LG=LevelGraph([1, 0, 0, 0],[[1, 2], [3, 4, 5], [6, 7, 8], [9, 10, 11]],[(2, 5), (3, 7), (4, 8), (6, 10), (1, 11)],{1: 0, 2: 0, 3: 0, 4: 0, 5: -2, 6: 2, 7: -2, 8: -2, 9: 4, 10: -4, 11: -2},[0, -1, -2, -3],True),dmp={9: (0, 0)},dlevels={0: 0, -1: -1, -2: -2, -3: -3}),
 EmbeddedLevelGraph(LG=LevelGraph([1, 0, 0, 0],[[1, 2], [3, 4, 5], [6, 7, 8], [9, 10, 11]],[(2, 5), (1, 7), (4, 8), (6, 10), (3, 11)],{1: 0, 2: 0, 3: 0, 4: 0, 5: -2, 6: 2, 7: -2, 8: -2, 9: 4, 10: -4, 11: -2},[0, -1, -2, -3],True),dmp={9: (0, 0)},dlevels={0: 0, -1: -1, -2: -2, -3: -3})]
\end{lstlisting}
Recall also \autoref{fourlevelred} and \autoref{fig:triplebanana}.

Determining how and when a profile becomes reducible is still very mysterious. It is also not immediately obvious from \verb|bot_to_bic| and \verb|top_to_bic|:
\begin{lstlisting}
sage: X.DG.bot_to_bic_inv(1)
{3: [0, 2], 6: [1, 4, 7], 0: [3], 4: [5], 2: [6]}
sage: X.DG.bot_to_bic_inv(3)
{6: [0, 1, 2]}
sage: X.DG.top_to_bic_inv(6)
{1: [0, 4], 4: [1], 5: [2], 3: [3]}
sage: X.DG.top_to_bic_inv(3)
{1: [0, 1, 2], 5: [3]}
\end{lstlisting}
Moreover, one should not expect the reducibility of a profile to be determined by the reducibility of the consecutive length-2 profiles appearing in it. In particular, while the product of the number of components of all length-2-profiles appearing as undegenerations gives an upper bound, this is very coarse and it is not clear how it can be improved (e.g. how to see the irreducibility of \verb|(1, 4, 6)| in this example).
\end{example}

As these examples illustrate, working with enhanced profiles can be quite subtle and is not yet completely understood. Therefore, \verb|diffstrata| has to occasionally work with the underlying graph. Most methods in \verb|GeneralisedStratum| that give relationships between the underlying \verb|LevelGraph|s of profiles start with \verb|explicit|. The backbone of \verb|is_degeneration| is \verb|explicit_leg_maps|, which raises a \verb|UserWarning| if no map is found, i.e. the enhanced profiles are not degenerations of each other. Using the notation from the above examples, we see e.g.
\begin{lstlisting}
sage: X.lookup_graph((5,))
EmbeddedLevelGraph(LG=LevelGraph([1, 2],[[1], [2, 3]],[(1, 3)],{1: 0, 2: 4, 3: -2},[0, -1],True),dmp={2: (0, 0)},dlevels={0: 0, -1: -1})
sage: X.lookup_graph((5,6))
EmbeddedLevelGraph(LG=LevelGraph([1, 1, 0],[[1], [2, 3, 4], [5, 6, 7]],[(1, 4), (2, 6), (3, 7)],{1: 0, 2: 2, 3: 0, 4: -2, 5: 4, 6: -4, 7: -2},[0, -1, -2],True),dmp={5: (0, 0)},dlevels={0: 0, -1: -1, -2: -2})
sage: X.explicit_leg_maps(((5,),0), ((5,6),0))
[{1: 1, 2: 5, 3: 4}]
\end{lstlisting}
Of course, if there are automorphisms involved, there will be many leg maps:
\begin{lstlisting}
sage: X.lookup_graph((1,))
EmbeddedLevelGraph(LG=LevelGraph([1, 1],[[1, 2], [3, 4, 5]],[(1, 4), (2, 5)],{1: 0, 2: 0, 3: 4, 4: -2, 5: -2},[0, -1],True),dmp={3: (0, 0)},dlevels={0: 0, -1: -1})
sage: X.lookup_graph((1,3))
EmbeddedLevelGraph(LG=LevelGraph([1, 0, 0],[[1, 2], [3, 4, 5], [6, 7, 8, 9]],[(2, 5), (1, 7), (3, 8), (4, 9)],{1: 0, 2: 0, 3: 0, 4: 0, 5: -2, 6: 4, 7: -2, 8: -2, 9: -2},[0, -1, -2],True),dmp={6: (0, 0)},dlevels={0: 0, -1: -1, -2: -2})
sage: X.explicit_leg_maps(((1,),0), ((1,3),0))
[{2: 1, 1: 2, 3: 6, 5: 7, 4: 5}, {2: 2, 1: 1, 3: 6, 5: 5, 4: 7}]
sage: len(X.lookup_graph((1,)).automorphisms)
2
sage: len(X.lookup_graph((1,3)).automorphisms)
2
\end{lstlisting}
In fact, \verb|explicit_leg_maps| first squishes the larger graph at the appropriate places to land in the profile of the smaller graph.
Recall from \autoref{subsec:graphsandprofiles} that this removes legs and vertices but leaves the numbering of the remaining graph untouched!
Then it goes through this profile and checks if any component is isomorphic to the squished graph. When one is found, all isomorphisms are returned.
Now, \verb|is_degeneration| simply returns a boolean if there exists a leg map.

Now, implementing the methods described in \autoref{deg:interface} from these building blocks is fairly straight-forward.

\subsection{Isomorphisms}
\label{subsec:isomorphisms}

The final step in the discussion of graphs and their degenerations is determining, when two graphs are isomorphic. Obviously, an isomorphism of \verb|EmbeddedLevelGraph|s must respect the embedding: recall, e.g., the three compact type BICs in the stratum $\Omega\cM_0(4,-2,-2,-2)$ from \autoref{triplebanana}. In that case, the underlying \verb|LevelGraph|s were identical and the isomorphism classes depended only on the embedding.

On the other hand, an isomorphism of \verb|EmbeddedLevelGraph|s can permute the marked points of a \verb|LevelGraph| if this is ``corrected'' by the embedding \verb|dmp|.

\begin{definition}
An \emph{isomorphism} of \verb|EmbeddedLevelGraph|s is a \verb|tuple| of \verb|dict|s, \verb|(isom_vertices, isom_legs)|, satisfying the following compatibility conditions:
\begin{enumerate}
\item \verb|isom_vertices| maps vertices to vertices (as indices of the corresponding \verb|list|s), respecting the genera and levels.
\item \verb|isom_legs| maps half-edges to half-edges, respecting the edges, the pole-orders and the marked points of the stratum (via \verb|dmp|).
\end{enumerate}
\end{definition}

\begin{rem}
The leg maps of \autoref{subsec:degenerations} consisted only of the \verb|isom_legs| component.
\end{rem}

The methods regarding isomorphisms are quite straight-forward to use. However, as the isomorphism classes are already encoded in the enhanced profiles, this ``higher-level'' access should be preferred. The number of automorphisms is, obviously, still important.

Automorphisms are recorded in a \verb|list|. We briefly illustrate that they respect the prongs:
\begin{lstlisting}
sage: X=Stratum((4,))
sage: symmetric_banana=EmbeddedLevelGraph(X,LG=LevelGraph([2, 0],[[1, 2], [3, 4, 5]],[(1, 4), (2, 5)],{1: 1, 2: 1, 3: 4, 4: -3, 5: -3},[0, -1],True),dmp={3: (0, 0)},dlevels={0: 0, -1: -1})
sage: symmetric_banana.automorphisms
[({0: 0, 1: 1}, {2: 1, 1: 2, 3: 3, 5: 4, 4: 5}),
 ({0: 0, 1: 1}, {2: 2, 1: 1, 3: 3, 5: 5, 4: 4})]
sage: asymmetric_banana=EmbeddedLevelGraph(X,LG=LevelGraph([2, 0],[[1, 2], [3, 4, 5]],[(1, 4), (2, 5)],{1: 2, 2: 0, 3: 4, 4: -4, 5: -2},[0, -1],True),dmp={3: (0, 0)},dlevels={0: 0, -1: -1})
sage: asymmetric_banana.automorphisms
[({0: 0, 1: 1}, {2: 2, 1: 1, 3: 3, 5: 5, 4: 4})]
\end{lstlisting}
The cardinality of the automorphism group is seen easily using \verb|len|:
\begin{lstlisting}
sage: len(symmetric_banana.automorphisms)
2
sage: len(asymmetric_banana.automorphisms)
1
sage: set([len(B.automorphisms) for B in X.bics])
{1, 2, 6}
\end{lstlisting}
We can also simply check if two graphs are isomorphic:
\begin{lstlisting}
sage: symmetric_banana.is_isomorphic(asymmetric_banana)
False
sage: symmetric_banana.is_isomorphic(symmetric_banana)
True
\end{lstlisting}
Note that the isomorphisms are generators (for checking, it is enough to construct the first). If we want to list all, we should convert them to a \verb|list|:
\begin{lstlisting}
sage: symmetric_banana.isomorphisms(asymmetric_banana)
<generator object EmbeddedLevelGraph.isomorphisms at 0x1cf5112a0>
sage: list(symmetric_banana.isomorphisms(asymmetric_banana))
[]
sage: list(symmetric_banana.isomorphisms(symmetric_banana))
[({0: 0, 1: 1}, {2: 1, 1: 2, 3: 3, 5: 4, 4: 5}),
 ({0: 0, 1: 1}, {2: 2, 1: 1, 3: 3, 5: 5, 4: 4})]
\end{lstlisting}

For the implementation, we decided to construct isomorphisms level by level.
An isomorphism of \verb|EmbeddedLevelGraph|s is then a set of compatible level isomorphisms.
We iterate through the isomorphisms on each level and yield whenever we find
compatible (i.e. respecting the edges) level isomorphisms for all levels.
Note that we use \verb|dlevels| for this, as these should be compatible.

An \emph{isomorphism of levels} now consists of a map (\verb|dict|) vertices to vertices and a map (\verb|dict|) legs to legs, respecting the genus, the number of legs on every vertex, the order at every leg and the marked points of the stratum (via \verb|dmp|).

These are constructed by the following algorithm:
\begin{algorithm}[Level Isomorphisms]\mbox{}
\begin{description}
\item[Step 1] Extract the information about the current level from the \verb|LevelGraph|. Note that we do not use \verb|level| to avoid all the overhead (residue conditions, etc.).
\item[Step 2] If the number of vertices, legs, legs per vertex, or the genera do not match, we are done: there can be no isomorphism.
\item[Step 3] The same marked points have to be on this level. Also, this gives the first part of the maps, as their legs (and thus also their vertices) must be mapped to each other. Note that we must also ensure that the vertices are compatible (same genera, numbers of legs, orders) and the marked points are split among vertices the same way.
\item[Step 4] For each genus $g$ appearing, we map the vertices of genus $g$ to the vertices of genus $g$ on the target level. For this, we use a simple recursive algorithm, enumerating all legal maps.
\item[Step 5] For each of these vertex maps, we construct all legal leg maps in a similar fashion (vertex by vertex).
\item[Step 6] Finally, we take the product of all constructed level isomorphisms and check which of these are compatible with the edges.
\end{description}
\end{algorithm}

This allows us to construct the degeneration graph of the isomorphism classes non-horizontal level graphs of any generalised stratum and reference the objects in a fairly efficient way.

%% file: sec_AGs.tex
\section{Additive Generators, Tautological Classes and Evaluation}
\label{sec:AGs}

The purpose of the \verb|diffstrata| package is to facilitate calculations in the tautological ring of strata. In light of \cite[Thm. 1.5]{strataEC}, any tautological class may be expressed as a formal sum of $\psi$-classes on graphs.

The \verb|diffstrata| package uses two classes to model the situation.
Denote by \verb|X| a \verb|GeneralisedStratum|.
\begin{itemize}
\item An \verb|AdditiveGenerator| encodes a product of $\psi$-classes on
  (various levels of) an \verb|EmbeddedLevelGraph| inside \verb|X|.
\item An \verb|ELGTautClass| is a formal sum of \verb|AdditiveGenerator|s, all on the same stratum \verb|X|.
\end{itemize}
Any \verb|ELGTautClass|es on \verb|X| can be added, multiplied and evaluated (i.e. integrated against the class of \verb|X|), although this will give $0$ if the class is not of top degree. The evaluation works by breaking down the expression into $\psi$-products on meromorphic strata and using the \verb|admcycles| package \cite{admcycles} to evaluate these, see \autoref{subsec:eval} for details.

Moreover, the class $\xi=c_1(\cO(-1))$ of \verb|X| is encoded using Sauvaget's relation (cf. \cite[Thm. 6(1)]{sauvaget}, \cite[Prop. 8.2]{strataEC} for the adaption to this setting) by \verb|X.xi|:
\begin{lstlisting}
sage: X=Stratum((2,))
sage: print(X.xi)
Tautological class on Stratum: (2,)
with residue conditions: []

3 * Psi class 1 with exponent 1 on level 0 * Graph ((), 0) +
-1 * Graph ((0,), 0) +
-1 * Graph ((1,), 0) +

sage: print(X.xi^X.dim())
Tautological class on Stratum: (2,)
with residue conditions: []

27 * Psi class 1 with exponent 3 on level 0 * Graph ((), 0) +
-9 * Psi class 1 with exponent 1 on level 0 * Psi class 2 with exponent 1 on level 1 * Graph ((0,), 0) +
-2 * Psi class 1 with exponent 1 on level 0 * Psi class 3 with exponent 1 on level 1 * Graph ((0,), 0) +
3 * Psi class 1 with exponent 1 on level 0 * Graph ((0, 1), 0) +
-1/2 * Psi class 1 with exponent 1 on level 0 * Psi class 2 with exponent 1 on level 0 * Graph ((1,), 0) +
-1/4 * Psi class 1 with exponent 2 on level 0 * Graph ((1,), 0) +
-1/4 * Psi class 2 with exponent 2 on level 0 * Graph ((1,), 0) +

sage: (X.xi^X.dim()).evaluate()
-1/640
\end{lstlisting}

We now describe these classes in some more detail and explain the evaluation process.

\subsection{AdditiveGenerators}

An \verb|AdditiveGenerator| requires three pieces of information:
\begin{itemize}
\item a \verb|GeneralisedStratum|, \verb|X|;
\item an \emph{enhanced profile} describing a level graph in \verb|X| (cf. \autoref{deg:interface}) and
\item a \verb|dict| encoding the powers of psi classes on this graph.
\end{itemize}
We explain the last item: the enhanced profile has an \verb|EmbeddedLevelGraph| associated to it and each leg of the underlying \verb|LevelGraph| is determined by a number (\verb|int|). The \verb|leg_dict| of an \verb|AdditiveGenerator| is a \verb|dict| with entries of the from \verb|{l : e}|, where \verb|l| is the number of a leg and \verb|e| is the exponent of the $\psi$-class at that leg (more precisely: the pullback of the $\psi$-class on the level the leg is on).

\begin{rem}
For caching purposes, \verb|AdditiveGenerator|s should never be created directly. Instead, the \verb|additive_generator| method of \verb|GeneralisedStratum| should be used:
\begin{lstlisting}
sage: A = X.additive_generator(((0,), 0), {2 : 1})
sage: print(A)
Psi class 2 with exponent 1 on level 1 * Graph ((0,), 0)
\end{lstlisting}
\end{rem}

\begin{rem}
\verb|AdditiveGenerator|s should be considered immutable, i.e. once generated they should not be changed. Indeed, each \verb|AdditiveGenerator| has a unique hash, consisting of the enhanced profile and the $\psi$-dictionary that is used for fast lookup and comparison.
\end{rem}

\subsection{Tautological Classes}
\label{subsec:tautclasses}

An \verb|ELGTautClass| consists of
\begin{itemize}
\item a \verb|GeneralisedStratum|, \verb|X| and
\item a \verb|list|, \verb|psi_list|, of \verb|tuple|s, each consisting of a coefficient (usually a rational number) and an \verb|AdditiveGenerator| on \verb|X|.
\end{itemize}
\verb|ELGTautClass| has a \verb|reduce| method that is called upon initialisation and after every operation. This method combines elements of \verb|psi_list| having the same \verb|AdditiveGenerator| and removes elements with coefficient zero or that vanish for dimension reasons.

\verb|ELGTautClass|es associated to the same stratum can be added (the \verb|psi_list|s are united and the \verb|reduce| method is applied) and multiplied (see \autoref{sec:mult}). Each of these operations yields again an \verb|ELGTautClass|.

\begin{rem}
For summing \verb|ELGTautClass|es, one may use the \verb|sum| function. However, this calls \verb|reduce| at every step of the summation. It is therefore often preferable to use \verb|GeneralisedStratum|'s \verb|ELGsum| method, which first concatenates all \verb|psi_list|s and only calls \verb|reduce| on the total list.
\end{rem}

Any \verb|AdditiveGenerator| may be converted to its associated \verb|ELGTautClass| using its \verb|as_taut| method. Any graph in \verb|X| can be converted to its associated \verb|ELGTautClass| using \verb|X.taut_from_graph|.
Moreover, every \verb|GeneralisedStratum|, \verb|X|, comes with a set of ``builtin'' tautological classes:
\begin{itemize}
\item \verb|X.ZERO| is the class with empty \verb|psi_list|.
\item \verb|X.ONE| is the class that consists only of the zero-level
  graph \verb|X.smooth_LG| representing smooth curves.
\item \verb|X.psi(l)| is the $\psi$-class associated to the marked
  point~\verb|l| of \verb|X|. Note that \verb|l| is a leg number of \verb|X.smooth_LG|. This should only be used for connected strata.
\item \verb|X.xi| is the class of $\xi$ in \verb|X| using Sauvaget's relation (cf. \cite[Thm. 6(1)]{sauvaget}, \cite[Prop. 8.2]{strataEC}). Note that this involves a choice of leg. By default, we choose the one resulting in the expression with the least number of summands. To choose a leg manually, use \verb|X.xi_with_leg| and provide a stratum point (cf. \autoref{pointtypes}).
\item \verb|X.exp_xi| is the class $\exp(\xi)$ on \verb|X|.
\item \verb|X.c1_E| is the first Chern class of $\Omega^1(\log)$ on \verb|X| (cf. \cite[Thm. 1.1]{strataEC}).
\item \verb|X.chern_class(n)| gives the \verb|n|-th Chern class of $\Omega^1(\log)$ on \verb|X| (cf. \cite[Thm. 9.10]{strataEC}).
\item \verb|X.calL()| (or more generally, \verb|X.calL(ep, i)| for an enhanced profile \verb|ep| associated to a graph $\Gamma$ in \verb|X| and a level $i$ of $\Gamma$) gives the class $\cL$ or more generally $\cL_{\Gamma}^{[i]}$, see \autoref{subsec:NB}.
\end{itemize}
Note that we can multiply an \verb|ELGTautClass| not only with any rational number, but also with any \verb|sage| object that can be multiplied with a rational number. In particular, this allows us to use \verb|sage| symbolic expressions:
\begin{lstlisting}
sage: X=Stratum((2,))
sage: var('a', 'b')
(a, b)
sage: T=a*X.ONE + b*X.psi(1)
sage: print(T^2)
Tautological class on Stratum: (2,)
with residue conditions: []

a^2 * Graph ((), 0) +
2*a*b * Psi class 1 with exponent 1 on level 0 * Graph ((), 0) +
b^2 * Psi class 1 with exponent 2 on level 0 * Graph ((), 0) +

sage: (T^3).evaluate()
1/1920*b^3
\end{lstlisting}

\subsection{Evaluation}
\label{subsec:eval}

Any \verb|ELGTautClass| can be evaluated, \verb|AdditiveGenerator|s that are not top degree will evaluate to $0$ without complaining:
\begin{lstlisting}
sage: X=Stratum((2,))
sage: (X.psi(1)).evaluate()
0
sage: (X.psi(1)^3).evaluate()
1/1920
\end{lstlisting}
Note that \verb|ELGTautClass| offers the check \verb|is_equidimensional| and the methods \verb|degree| and \verb|list_by_degree|:
\begin{lstlisting}
sage: X.xi.is_equidimensional()
True
sage: T=X.ONE + X.psi(1)
sage: T.is_equidimensional()
False
sage: print(T.degree(1))
Tautological class on Stratum: (2,)
with residue conditions: []

1 * Psi class 1 with exponent 1 on level 0 * Graph ((), 0) +

sage: print(T.degree(0))
Tautological class on Stratum: (2,)
with residue conditions: []

1 * Graph ((), 0) +
\end{lstlisting}
Observe that \verb|list_by_degree| is a \verb|list| of length \verb|X.dim()| filled up with \verb|X.ZERO|:
\begin{lstlisting}
sage: T.list_by_degree()
[ELGTautClass(X=GeneralisedStratum(sig_list=[Signature((2,))],res_cond=[]),psi_list=[(1, AdditiveGenerator(X=GeneralisedStratum(sig_list=[Signature((2,))],res_cond=[]),enh_profile=((), 0),leg_dict={}))]),
 ELGTautClass(X=GeneralisedStratum(sig_list=[Signature((2,))],res_cond=[]),psi_list=[(1, AdditiveGenerator(X=GeneralisedStratum(sig_list=[Signature((2,))],res_cond=[]),enh_profile=((), 0),leg_dict={1: 1}))]),
 ELGTautClass(X=GeneralisedStratum(sig_list=[Signature((2,))],res_cond=[]),psi_list=[]),
 ELGTautClass(X=GeneralisedStratum(sig_list=[Signature((2,))],res_cond=[]),psi_list=[])]
\end{lstlisting}

The key to evaluating an \verb|ELGTautClass| is to split it into pieces that can be evaluated by \verb|admcycles| using the expression of the class of a stratum in the tautological ring of $\overline{\cM}_{g,n}$ by Sauvaget \cite{sauvaget}. For this purpose, the \verb|stack_factor| of an \verb|AdditiveGenerator| is defined as follows: let \verb|G| be the associated graph, then the \verb|stack_factor| is the quotient of the product of the prongs of \verb|G| and the product of \verb|B.ell| for every \verb|B| appearing in the profile of \verb|G| and the number of automorphisms of \verb|G|. For a BIC \verb|B|,
the number \verb|B.ell| is the $\lcm$ of the prongs of \verb|B|. By
\cite[Lemma 9.12]{strataEC} this factor is necessary to convert an integral
over a boundary stratum $D_\Gamma$ into the product of level-wise evaluations.
\par
\begin{algorithm}[Evaluation]\mbox{}
\begin{description}
\item[Step 1] Go through \verb|psi_list| and take the sum of the evaluation of the \verb|AdditiveGenerators| multiplied with their respective coefficients.
\item[Step 2] Each \verb|AdditiveGenerator| is evaluated by sorting the $\psi$-classes by level and taking the product of the evaluations of these $\psi$-polynomials on each level. This product is then multiplied with the \verb|stack_factor|.
\item[Step 3] If the residue space~$\frakR$ is empty, we evaluate as follows. If the level is disconnected, it evaluates to $0$; if it is $0$-dimensional, it evaluates to $1$; otherwise we use the \verb|Strataclass| function of \verb|admcycles| to generate the tautological class of the stratum on $\overline{\cM}_{g,n}$ and integrate this against the $\psi$-classes, using \verb|admcycle|'s \verb|evaluate|.
\item[Step 4] If the level splits as a product (with residue conditions, i.e. on the level of the underlying graph, cf. \autoref{subsec:underlying_graph}), it evaluates to $0$, since the fiber dimension to the product of moduli spaces is
positive and the $\psi$-classes are pullbacks from there.
\par
Otherwise, we create a new \verb|GeneralisedStratum| with one residue condition removed. We repeat until this condition is non-trivial (or $\frakR = \emptyset$) and evaluate the product of the (original) $\psi$-expression and the class cut out by this residue condition (using \cite[Prop. 8.3]{strataEC}). The class cut out by a residue condition inside a \verb|GeneralisedStratum| can be obtained by \verb|res_stratum_class|.
\end{description}
\end{algorithm}
\par
\begin{rem}
We provide some more details on the above algorithm.
\begin{enumerate}
\item If a level of an \verb|AdditiveGenerator| evaluates to $0$, the evaluation returns immediately and lower levels are not evaluated.
\item The evaluations in Step 3 are based on the algorithms in \cite{admcycles}.
In particular the \verb|Strataclass| function is based on the description
of fundamental classes of (non-generalised) strata
conjectured in \cite{FP} and \cite{schmittDimTh} and proven recently
in \cite{BHPSS} based on results from \cite{HolmesSchmitt}.
\end{enumerate}
\end{rem}
\par
\begin{example}
We illustrate the calculation of a class cut out by a residue condition:
\begin{lstlisting}
sage: X=GeneralisedStratum([Signature((23,5,-13,-17))])
sage: assert X.res_stratum_class([(0,2)]).evaluate() == 5
\end{lstlisting}
In fact, this stratum has three boundary points corresponding
to graphs $\Gamma_1$, $\Gamma_2$ and~$\Gamma_3$ that have
respectively the marked points of order $(23,5)$, $(23,-13)$ and $(23,-17)$ on
lower level. By \cite[Prop.~8.2]{strataEC} we can express~$\xi$
using the first $\psi$-class as $\int_X \xi = 24-29-11-7 = -23$.
Now in \cite[Prop. 8.3]{strataEC} the contributing boundary
graphs are~$\Gamma_2$ (because the point with order $-13$ is on lower
level) and~$\Gamma_3$ (because the zero residue at the point of
order $-13$ on upper level is automatic), but not~$\Gamma_1$. We
find that the evaluation of the boundary stratum is
$-(-23) - 11 - 7 = 5$, as claimed.
\end{example}

%% file: sec_mult.tex
\section{Normal Bundles, Pullback and Multiplication}
\label{sec:mult}

While adding two \verb|AdditiveGenerator|s is straight-forward, expressing the intersection of two \verb|AdditiveGenerator|s again as a sum of \verb|AdditiveGenerator|s is subtle, in particular if the intersection is not transversal.

A first approximation is finding common degenerations of two graphs, but if these graphs have a common undegeneration, there is a normal bundle contribution.
The precise answer is the excess intersection formula \cite[Prop. 8.1]{strataEC}.

\begin{example}\label{multiplication}
Consider the minimal stratum in genus two, $\Omega\cM_2(2)$. The class of $\xi$ can be expressed, by Sauvaget's relation, as a sum of \verb|AdditiveGenerator|s (cf. \autoref{sec:AGs}).
\begin{lstlisting}
sage: X=Stratum((2,))
sage: print(X.xi)
Tautological class on Stratum: (2,)
with residue conditions: []

3 * Psi class 1 with exponent 1 on level 0 * Graph ((), 0) +
-1 * Graph ((0,), 0) +
-1 * Graph ((1,), 0) +
\end{lstlisting}
As $\dim\Omega\cM_2(2)=3$, we can evaluate $\xi^3$. However, calculating this class requires several applications of the excess intersection formula described above:
\begin{lstlisting}
sage: X.dim()
3
sage: print(X.xi^3)
Tautological class on Stratum: (2,)
with residue conditions: []

27 * Psi class 1 with exponent 3 on level 0 * Graph ((), 0) +
-9 * Psi class 1 with exponent 1 on level 0 * Psi class 2 with exponent 1 on level 1 * Graph ((0,), 0) +
-2 * Psi class 1 with exponent 1 on level 0 * Psi class 3 with exponent 1 on level 1 * Graph ((0,), 0) +
3 * Psi class 1 with exponent 1 on level 0 * Graph ((0, 1), 0) +
-1/2 * Psi class 1 with exponent 1 on level 0 * Psi class 2 with exponent 1 on level 0 * Graph ((1,), 0) +
-1/4 * Psi class 1 with exponent 2 on level 0 * Graph ((1,), 0) +
-1/4 * Psi class 2 with exponent 2 on level 0 * Graph ((1,), 0) +
\end{lstlisting}
The resulting top degree class can now be evaluated (cf. \autoref{subsec:eval}) to give the expected result:
\begin{lstlisting}
sage: (X.xi^3).evaluate()
-1/640
\end{lstlisting}
\end{example}

Before going into the details of the implementation, we illustrate the use and necessity of multiplication inside an ambient stratum. The standard operations \verb|*| and \verb|^| are performed in the Chow ring of the stratum as illustrated above. We can instead work in the Chow ring of any substratum $D_\Gamma$ by specifying as \verb|ambient| the enhanced profile of $\Gamma$.

\begin{example}\label{topxieval}
Consider the stratum $\Omega\cM_2(1,1)$ and the boundary divisor represented by a compact-type graph with a genus two component on top and a genus zero component on bottom level. By inspecting the \verb|list| of BICs we find that it is BIC \verb|3| (cf. \autoref{numbering}):
\begin{lstlisting}
sage: Y=Stratum((1,1))
sage: Y.bics[3]
EmbeddedLevelGraph(LG=LevelGraph([2, 0],[[1], [2, 3, 4]],[(1, 4)],{1: 2, 2: 1, 3: 1, 4: -4},[0, -1],True),dmp={2: (0, 0), 3: (0, 1)},dlevels={0: 0, -1: -1})
\end{lstlisting}
As the stratum \verb|Y| is four-dimensional, we can multiply this graph with $\xi^3$ to obtain a top-degree class that we may evaluate to find a number:
\begin{lstlisting}
sage: Y.dim()
4
sage: (Y.xi^3*Y.additive_generator(((3,),0))).evaluate()
-1/640
\end{lstlisting}
This matches the observation that the top level is the minimal stratum in genus two where $\xi^3$ evaluates to $-\frac{1}{640}$, the bottom level is a point and there are no prongs or automorphisms involved, compare %
\cite[\S 4.3 and Lemma 9.12]{strataEC}:
\begin{lstlisting}
sage: Y.bics[3].top
LevelStratum(sig_list=[Signature((2,))],res_cond=[],leg_dict={1: (0, 0)})
sage: (Y.bics[3].top.xi^3).evaluate()   # calculating on Y.bics[3].top
-1/640
\end{lstlisting}
We can perform the same calculation in \verb|Y| using \verb|xi_at_level|, see \autoref{xiatlevel}. However, as the class $\xi^{[i]}_{\Gamma}$ correspoding to \verb|xi_at_level| lives not on \verb|Y| but on the BIC \verb|3|,
performing the normal multiplication with \verb|*| or \verb|^| will not yield the correct result: we must use our BIC $\Gamma$ as the ambient stratum, as we want to 
perform the multiplication in $\CH(D_\Gamma)$, i.e. \emph{before} pushing forward to \verb|Y|.

Indeed, while the evaluation of the cube of \verb|xi_at_level| performed in \verb|Y| is
\begin{lstlisting}
sage: (Y.xi_at_level(0, ((3,),0))^3).evaluate()
0
\end{lstlisting}
we obtain the expected result when using the correct ambient stratum, the BIC \verb|3|:
\begin{lstlisting}
sage: CT_xi_top = Y.xi_at_level(0, ((3,),0))
sage: (Y.intersection(Y.intersection(CT_xi_top, CT_xi_top, ((3,),0)), CT_xi_top, ((3,),0))).evaluate()
-1/640
\end{lstlisting}
In the special case of taking exponents of \verb|xi_at_level|, we can use the method \verb|xi_at_level_pow|, which computes the exponent with the correct ambient graph:
\begin{lstlisting}
sage: (Y.xi_at_level_pow(0, ((3,),0), 3)).evaluate()
-1/640
\end{lstlisting}
Indeed, in this case we may even use the method \verb|top_xi_at_level| which computes and evaluates the top-power of $\xi$:
\begin{lstlisting}
sage: Y.top_xi_at_level(((3,),0), 0)
-1/640
\end{lstlisting}
In fact, whenever possible this method should be used (even for $\xi$ on the whole stratum), as the results are cached and reused, see \autoref{sec:caching}.
\end{example}

\begin{rem}\label{xipowzero}
Note that by Sauvaget \cite[Prop. 3.3]{SauvagetMinimal}, for a connected holomorphic stratum $X$ of genus $g$, $\int_X\xi^k=0$ for $k\geq 2g$. In this case, nothing is computed and \verb|X.ZERO| is returned immediately.
\end{rem}

\subsection{Interface}
\label{subsec:mult:int}

While the operations \verb|*| and \verb|^| should work intuitively, the syntax of the underlying methods is more involved. In particular, when enhanced profiles, \verb|AdditiveGenerator|s and \verb|ELGTautClass|es appear is a subtle issue that we try to explain in the following summary.

In the following, \verb|ambient| is always an \emph{enhanced profile} corresponding to the ambient stratum in which the multiplication is taking place (see \autoref{subsec:intersection} for details). By default, it always corresponds to the full enveloping stratum, i.e. the empty enhanced profile \verb|((), 0)|. 

The following are methods of \verb|GeneralisedStratum|:
\begin{itemize}
\item\verb|intersection| calculates the product of two tautological classes inside \verb|ambient|.
\begin{description}
\item[Arguments] \verb|ELGTautClass|, \verb|ELGTautClass|, \verb|ambient|
\item[Returns] \verb|ELGTautClass|
\end{description}
\item\verb|intersection_AG| calculates the product of two additive generators inside \verb|ambient|
\begin{description}
\item[Arguments] \verb|AdditiveGenerator|, \verb|AdditiveGenerator|, \verb|ambient|
\item[Returns] \verb|ELGTautClass|
\end{description}
\item\verb|normal_bundle| calculates the normal bundle of a graph inside \verb|ambient|.
\begin{description}
\item[Arguments] enhanced profile, \verb|ambient|
\item[Returns] \verb|ELGTautClass|
\end{description}
\item\verb|cnb| calculates the common normal bundle of two graphs inside \verb|ambient|.
\begin{description}
\item[Arguments] Enhanced profile, enhanced profile, \verb|ambient|
\item[Returns] \verb|ELGTautClass|
\end{description}
\item\verb|gen_pullback| calculates the generalised pullback of an additive generator to a graph inside \verb|ambient|.
\begin{description}
\item[Arguments] \verb|AdditiveGenerator|, enhanced profile, \verb|ambient|
\item[Returns] \verb|ELGTautClass|
\end{description}
\item\verb|gen_pullback_taut| calculates the generalised pullback of a tautological class to a graph inside \verb|ambient|.
\begin{description}
\item[Arguments] \verb|ELGTautClass|, enhanced profile, \verb|ambient|
\item[Returns] \verb|ELTTautClass|
\end{description}
\end{itemize}

Moreover, the graphs needed are calculated by a series of methods that are straight-forwardly implemented using the observations of \autoref{subsec:degenerations}. In contrast to the methods described above, the arguments are always enhanced profiles and they return \verb|list|s of enhanced profiles:
\begin{itemize}
\item\verb|common_degenerations| finds the common degenerations of two graphs;
\item\verb|codim_one_degenerations| finds degenerations of an enhanced profile associated to a graph with one more level;
\item\verb|codim_one_common_undegenerations| finds degenerations of \verb|ambient| with one more level that are additionally undegenerations of two given enhanced profiles;
\item\verb|minimal_common_undegeneration| finds the graph of minimal dimension that is undegeneration of two supplied graphs.
\end{itemize}
Finally, an expression of the form \verb|A*B| is evaluated according to the following rules:
\begin{itemize}
\item if \verb|A| and \verb|B| are both \verb|ELGTautClass|es, \verb|intersection| is called (with empty \verb|ambient|);
\item if one of the two is an \verb|AdditiveGenerator| and the other an \verb|ELGTautClass|, they are multiplied using the \verb|as_taut| method;
\item otherwise, if one of the two is an \verb|ELGTautClass|, scalar multiplication with the other is attempted. This can be used, e.g., for multiplication with symbolic variables, cf. \autoref{subsec:tautclasses}.
\end{itemize}
Note that two \verb|AdditiveGenerator|s may only be multiplied (using \verb|*|) if the underlying graphs are the same; otherwise \verb|intersection_AG| should be used.
Moreover, multiplying with \verb|X.ONE| is the identity and with \verb|0| or \verb|X.ZERO| gives \verb|X.ZERO|.

Using \verb|^|, we may calculate powers of a \verb|ELGTautClass| (with empty \verb|ambient|). The method \verb|pow| of \verb|GeneralisedStratum| allows the specification of an \verb|ambient|.

\subsection{Intersections}
\label{subsec:intersection}

Mathematically, the situation is the following: we implement the general version of the excess intersection formula \cite[eq. (61)]{strataEC}. 

Let $\Lambda_1$ and $\Lambda_2$ be degenerations of a $k$-level graph $\Gamma$ and denote the associated strata by $D_{\Lambda_1}$, $D_{\Lambda_2}$ and $D_\Gamma$, as in the following diagram:
\[	\begin{tikzcd}
	D_\Pi  \arrow{r}{\frakj_{\Pi,\Lambda_2}} 
	\arrow{d}{\frakj_{\Pi,\Lambda_1}}  & D_{\Lambda_2}  \arrow{d}{\frakj_{\Lambda_2,\Gamma}}\\
	D_{\Lambda_1} \arrow{r}{\frakj_{\Lambda_1,\Gamma}} & D_{\Gamma}
	\end{tikzcd}
\]
Furthermore, we define $\nu^\Pi_{(\Lambda_1\cap \Lambda_2)/\Gamma}$
to be the product of the pullback to $D_\Pi$ of the normal bundles
$\cN_{\Gamma',\Gamma}$ where~$\Gamma'$ ranges
over all $k+1$-level graphs~$\Gamma'$ that are a degeneration of $\Gamma$ and
that are moreover common to $\Lambda_1$ and $\Lambda_2$, see \autoref{subsec:NB}.
For any $\alpha \in \CH^\bullet(D_{\Lambda_2})$ we can then express its push-forward 
pulled back to $\Lambda_1$ as
\[
\frakj_{\Lambda_1,\Gamma}^* \frakj_{\Lambda_2,\Gamma*} \alpha \= \sum_{\Pi}
\frakj_{\Pi, \Lambda_1,*} \Bigl(\nu^\Pi_{(\Lambda_1\cap \Lambda_2)/\Gamma} \cdot \frakj_{\Pi, \Lambda_2}^* 
\alpha \Bigr)\,,
\]
where the sum ranges over all $(\Lambda_1,\Lambda_2)$-graphs~$\Pi$.
This corresponds to the product of $\alpha$ with the class of $\Lambda_2$ in $\CH(D_\Gamma)$.
A level graph $\Pi$ is {\em a $({\Lambda_1,\Lambda_2})$-graph} if there
are undegeneration morphisms $\rho_i \colon \Pi \to \Lambda_i$, i.e.\ edge
contraction morphisms with the property that there are subsets~$I$ and~$J$
of levels such that $\delta_I(\Pi) = \Lambda_1$ and $\delta_J(\Pi) = \Lambda_2$.
For details, see \cite[\S 8.1]{strataEC}.

\begin{rem}
We briefly summarise the relationship to \autoref{subsec:mult:int}. The correspondence between the mathematical and the \verb|diffstrata| objects is as follows:
\begin{itemize}
\item the enveloping (projective) stratum is a \verb|GeneralisedStratum| \verb|X|;
\item the graphs $\Lambda_1, \Lambda_2$ correspond to enhanced profiles \verb|ep_1| and \verb|ep_2| in \verb|X|;
\item the graph $\Gamma$ corresponds to the enhanced profile \verb|ambient|;
\item the class $\alpha$ corresponds to an \verb|AdditiveGenerator| \verb|A| with underlying graph \verb|ep_1|.
\end{itemize}
Then we may perform the following calculations in \verb|diffstrata|:
\begin{itemize}
\item \verb|X.common_undegenerations(ep_1, ep_2)| gives the set of $(\Lambda_1,\Lambda_2)$-graphs;
\item each $\rho_i$ is given by \verb|X.explicit_leg_maps|;
\item $\frakj_{\Pi, \Lambda_2}^* \alpha$ is given by \verb|X.gen_pullback(A, pi_i, ambient)| for \verb|pi_i| an enhanced profile in \verb|X.common_undegenerations(ep_1, ep_2)|;
\item the individual normal bundles $\cN_{\Gamma',\Gamma}$ can be calculated using \verb|normal_bundle|;
\item \verb|X.cnb(ep_1, ep_2, ambient)| corresponds to the product of normal bundles appearing in $\nu_{(\Lambda_1\cap \Lambda_2)/\Gamma}$ (not yet pulled back to $\Lambda_1\cap \Lambda_2$, see \autoref{subsec:NB});
\item $\frakj_{\Lambda_1,\Gamma}^* \frakj_{\Lambda_2,\Gamma*} \alpha$ is given by setting \verb|AG_2=X.additive_generator(ep_2)| and then performing
\verb|X.intersection_AG(A, AG_2, ambient)|.
\end{itemize}
Of course, the product of two \verb|AdditiveGenerator|s on the same graph is simply the product of their $\psi$-polynomials (sum of the \verb|leg_dict|s) and the product of two tautological classes is the sum of the products of their \verb|AdditiveGenerator|s.
\end{rem}
Therefore, the excess intersection formula is given by a pullback as the base case (no common intersection) and recursive multiplication with normal bundles.
Note that the codimension of \verb|ambient| increases at each step and therefore the normal bundle becomes trivial after finitely many iterations.

\subsection{Normal bundles}
\label{subsec:NB}

The key ingredient for the multiplication is the calculation of normal bundles. 

The first Chern class of the normal bundle of a divisor is computed in \cite[Thm 7.1]{strataEC}. In \verb|diffstrata|, this can be computed using \verb|normal_bundle|:
\begin{lstlisting}
sage: X=Stratum((2,))
sage: X.normal_bundle(((0,),0)) == X.taut_from_graph((0,),0)^2
True
sage: X.normal_bundle(((1,),0)) == X.taut_from_graph((1,),0)^2
True
\end{lstlisting}
However, for the excess intersection formula, we need
to compute $\nu^\Pi_{(\Lambda_1\cap \Lambda_2)/\Gamma}$. The situation is summarised in the following diagram of the involved graphs:
\[	\begin{tikzcd}
	\Lambda_1\cap\Lambda_2  \arrow[rrr] \arrow[ddd]  \arrow[rd, dashed] 
	& && \Lambda_2  \arrow[ddd] \arrow[dll] \arrow[ddl]\\
	& \Gamma'' \arrow[rd]& \\
	&& \Gamma'\arrow[rd, "\mathrm{codim} 1" description] & \\
	\Lambda_1 \arrow[rrr] \arrow[uur] \arrow[urr]& && \Gamma
	\end{tikzcd}
\]
where the arrows represent undegeneration maps. The graph $\Gamma''$ is the \emph{minimal common undegeneration} of $\Lambda_1$ and $\Lambda_2$ (\verb|X.minimal_common_undegeneration|), corresponding to the intersection of the profiles, while $\Lambda_1\cap\Lambda_2$ corresponds to the union of the profiles (up to reducibility issues, cf. \autoref{subsec:degenerations}). The graph(s) $\Gamma'$ are codimension one degenerations of $\Gamma$ that are common undegenerations of $\Lambda_1$ and $\Lambda_2$ (\verb|X.codim_one_common_undegenerations|).

For the excess intersection formula, we need the product of the normal bundles $\cN_{\Gamma',\Gamma}=\cN_{D_{\Gamma'}/D_{\Gamma}}$ where $\Gamma'$ is a $k+1$-level degeneration of the $k$-level graph $\Gamma$. 
More precisely, the normal bundles $\cN_{D_{\Gamma'}/D_{\Gamma}}$ are pulled back to $D_{\Gamma''}$ and the product is computed in $\CH(D_{\Gamma''})$.
The normal bundle is computed in \cite[Prop. 7.5]{strataEC}:
\[
c_1(\cN_{\Gamma',\Gamma}) \= \frac{1}{\ell_{\delta_i(\Gamma')}} \bigl(-\xi_\Gamma'^{[i]}
- c_1(\cL_\Gamma'^{[i]}) + \xi_\Gamma'^{[i+1]} \bigr)\quad \text{in} \quad
\CH^1(D_{\Gamma'})\,,
\]
where 
\[
\cL_{\Gamma'}^{[i]} \=  \cO_{D_\Gamma'}
\Bigl(\sum_{\Gamma' \overset{[i]}{\rightsquigarrow} %
  \wh{\Delta} } \ell_{\delta_{i+1}(\wh{\Delta})}D_{\wh{\Delta}} \Bigr)\,,
\]
where the sum runs over all graphs $\wh{\Delta} \in \LG_{k+2}(\ol{B})$
that yield divisors in~$D_\Gamma'$ by splitting the $i$-th level. These are then pulled back to $D_{\Gamma''}$, see \autoref{subsec:pullback}, and then they are multiplied (in $\CH(D_{\Gamma''})$, i.e. using \verb|ambient| $\Gamma''$).

\begin{rem}
\label{nb:terminates}
Observe that this recursive procedure terminates: indeed, the product in $\CH(D_{\Gamma})$ has been transformed to a product in $\CH(D_{\Gamma''})$ which is of strictly smaller dimension. Moreover, if the dimension is small enough, transversality is ensured by dimension reasons.
\end{rem}

In \verb|diffstrata|, $\xi_\Gamma^{[i]}$ is given by \verb|X.xi_at_level(i, ep)|, where \verb|ep| is the enhanced profile corresponding to $\Gamma$ in \verb|X|, see \autoref{xiatlevel}. The easiest way to produce $\cL_{\Gamma}^{[i]}$ is by generating all the BICs inside the \verb|GeneralisedStratum| at level $i$ and glue these into $\Gamma$ (see \autoref{clutching} for details about gluing BICs into a level). The class of $\cL^{[i]}_\Gamma$ is computed by the method \verb|calL|.

Recall that $\ell_{\delta_i(\Gamma)}$ is the $\lcm$ of the prongs of the BIC $\delta_i(\Gamma)$. In \verb|diffstrata| this is stored in the attribute \verb|ell| of \verb|EmbeddedLevelGraph| (if it is a BIC).

We summarise the normal bundle algorithm:
\begin{algorithm}[Normal Bundle]\mbox{}
\begin{description}
\item[Step 1] Compute the \verb|minimal_common_undegeneration|, \verb|min_com| (corres\-ponding to $\Gamma''$).
\item[Step 2] If \verb|min_com| is \verb|ambient| the intersection is transversal, we return \verb|1|.
\item[Step 3] Loop through \verb|ep| in \verb|codim_one_common_undegenerations| (corresponding to $\Gamma'$).
\item[Step 4] Calculate the level $i$ where \verb|ep| and \verb|min_com| differ.
\item[Step 5] Calculate the normal bundle as in \cite[Eq. (58)]{strataEC} with
  the function \verb|xi_at_level| for $\xi^{[i]}$ and \verb|glue_AG_at_level| for $\cL$ (see \autoref{clutching} for details).
\item[Step 6] Pull this normal bundle back to \verb|min_com|, cf. \autoref{subsec:pullback}.
\item[Step 7] Return the product of the normal bundles (one for each \verb|ep|) inside \verb|min_com|.
\end{description}
\end{algorithm}

\begin{rem}
Note that in the case of a transversal intersection, the ``dummy'' value \verb|1| is returned:
\begin{lstlisting}
sage: X.cnb(((1,),0),((0,),0))
1
\end{lstlisting}
The reason for this inconsistency is that the normal bundle should correspond to the class \verb|ONE|, but inside an \verb|ambient| this would be \verb|ambient|, which is not what we want. Therefore, this case must be handled separately!
\end{rem}

\begin{example}
\label{ex:cnb}
We continue in the setting of \autoref{topxieval} in the boundary of $\Omega\cM_2(1,1)$. Recall that \verb|Y.bics[3]| was the compact-type graph with top-level an $\cM_2(2)$. We may intersect this graph with the banana graph and the other compact-type graph in \verb|Y|:
\begin{lstlisting}
sage: Y.bics[3]
EmbeddedLevelGraph(LG=LevelGraph([2, 0],[[1], [2, 3, 4]],[(1, 4)],{1: 2, 2: 1, 3: 1, 4: -4},[0, -1],True),dmp={2: (0, 0), 3: (0, 1)},dlevels={0: 0, -1: -1})
sage: Y.bics[1]
EmbeddedLevelGraph(LG=LevelGraph([1, 0],[[1, 2], [3, 4, 5, 6]],[(1, 5), (2, 6)],{1: 0, 2: 0, 3: 1, 4: 1, 5: -2, 6: -2},[0, -1],True),dmp={3: (0, 0), 4: (0, 1)},dlevels={0: 0, -1: -1})
sage: Y.bics[2]
EmbeddedLevelGraph(LG=LevelGraph([1, 1],[[1], [2, 3, 4]],[(1, 4)],{1: 0, 2: 1, 3: 1, 4: -2},[0, -1],True),dmp={2: (0, 0), 3: (0, 1)},dlevels={0: 0, -1: -1})
\end{lstlisting}
Calculating the common normal bundle of these intersections we obtain:
\begin{lstlisting}
sage: print(Y.cnb(((1,3),0),((2,3),0)))
Tautological class on Stratum: (1, 1)
with residue conditions: []

-1 * Psi class 1 with exponent 1 on level 0 * Graph ((3,), 0) +
\end{lstlisting}
As expected, this is the normal of the BIC \verb|3|, as the other two BICs intersect transversally:
\begin{lstlisting}
sage: print(Y.normal_bundle(((3,),0)))
Tautological class on Stratum: (1, 1)
with residue conditions: []

-1 * Psi class 1 with exponent 1 on level 0 * Graph ((3,), 0) +
\end{lstlisting}
Calculating, e.g.
\begin{lstlisting}
sage: print(Y.cnb(((1,3),0),((1,),0)))
Tautological class on Stratum: (1, 1)
with residue conditions: []

-1/2 * Psi class 2 with exponent 1 on level 0 * Graph ((1,), 0) +
-1/2 * Psi class 1 with exponent 1 on level 0 * Graph ((1,), 0) +
-1/2 * Psi class 5 with exponent 1 on level 1 * Graph ((1,), 0) +
-1/2 * Psi class 6 with exponent 1 on level 1 * Graph ((1,), 0) +
\end{lstlisting}
gives the normal bundle of the banana graph.
\end{example}

\subsection{Pulling back classes}
\label{subsec:pullback}

To calculate the pullback of an \verb|AdditiveGenerator|, we consider first the base case and then the generalised case. 

Let \verb|ep| be the enhanced profile of a graph $\Lambda$ in \verb|X|, \verb|A| an \verb|AdditiveGenerator| on \verb|ep| corresponding to a class $\alpha$ on $D_\Lambda$ and \verb|ep_deg| the enhanced profile of a degeneration $\Pi$ of $\Lambda$, i.e. we obtain $\Lambda$ by contracting some of the level-crossings of $\Pi$. Then there are finitely many contraction morphisms $\rho\colon\Pi\to\Lambda$ and each of these gives a well-defined pullback map of $\alpha$. The pullback is the weighted sum of these and is given by \verb|A.pull_back(ep_deg)|.

\begin{example}
Consider the minimal stratum in genus $2$, $\Omega\cM_2(2)$, with the notation of \autoref{H2}, see also \autoref{degenerationgraphH2}. Let \verb|A| denote the $\psi$-class on the top-level of the compact-type divisor:
\begin{lstlisting}
sage: X=Stratum((2,))
sage: X.bics
[EmbeddedLevelGraph(LG=LevelGraph([1, 0],[[1, 2], [3, 4, 5]],[(1, 4), (2, 5)],{1: 0, 2: 0, 3: 2, 4: -2, 5: -2},[0, -1],True),dmp={3: (0, 0)},dlevels={0: 0, -1: -1}),
 EmbeddedLevelGraph(LG=LevelGraph([1, 1],[[1], [2, 3]],[(1, 3)],{1: 0, 2: 2, 3: -2},[0, -1],True),dmp={2: (0, 0)},dlevels={0: 0, -1: -1})]
sage: A=X.additive_generator(((1,),0), {1:1})
sage: print(A)
Psi class 1 with exponent 1 on level 0 * Graph ((1,), 0)
\end{lstlisting}
We may pull this class back to the intersection with the banana graph \verb|((1,0), 0)|:
\begin{lstlisting}
sage: print(A.pull_back(((1,0),0)))
Tautological class on Stratum: (2,)
with residue conditions: []

1 * Psi class 1 with exponent 1 on level 0 * Graph ((1, 0), 0) +
\end{lstlisting}
Considering instead the class \verb|B| of the $\psi$-class on bottom-level of the same divisor:
\begin{lstlisting}
sage: B=X.additive_generator(((1,),0), {3:1}); print(B)
Psi class 3 with exponent 1 on level 1 * Graph ((1,), 0)
sage: print(B.pull_back(((1,0),0)))
Tautological class on Stratum: (2,)
with residue conditions: []
\end{lstlisting}
It vanishes for dimension reasons when pulled back to the intersection.

Finally, consider the class \verb|S| of one of the top half-legs of the banana graph. Pulling it back to the intersection will also vanish for dimension reasons, but pulling it back to the banana graph itself illustrates the weighted sum: there are two graph morphisms (switching the edges).
\begin{lstlisting}
sage: S=X.additive_generator(((0,),0), {1:1}); print(S)
Psi class 1 with exponent 1 on level 0 * Graph ((0,), 0)
sage: print(S.pull_back(((0,),0)))
Tautological class on Stratum: (2,)
with residue conditions: []

1/2 * Psi class 2 with exponent 1 on level 0 * Graph ((0,), 0) +
1/2 * Psi class 1 with exponent 1 on level 0 * Graph ((0,), 0) +
\end{lstlisting}
Notice that in all of the above cases an \verb|ELGTautClass| is returned.
\end{example}

The implementation of \verb|pull_back| is simply a matter of transforming the \verb|leg_dict| of \verb|A| via the map $\rho$ to \verb|ep_deg| and dividing by the number of these maps. The maps are given by \verb|explicit_leg_maps|, cf. \autoref{subsec:degenerations}.

However, for the excess intersection formula, we require a more general notion of pullback. As above, let \verb|ep| be the enhanced profile of a graph $\Lambda$ in \verb|X|, \verb|A| an \verb|AdditiveGenerator| on \verb|ep| corresponding to a class $\alpha$ on $D_\Lambda$, but now we do not require the ``target'' graph $\Lambda'$ to be a degeneration of $\Lambda$. Instead, we will pull $\alpha$ back to the intersection $\Lambda\cap\Lambda'$:
\[	\begin{tikzcd}
	\Lambda\cap\Lambda'\arrow[rr]\arrow[dd]\arrow[rd]
	&& \Lambda   \arrow[dd] \arrow[dl]&&\\
	& \Gamma'' \arrow[rd]& \\
	\Lambda' \arrow[rr] \arrow[ur]&& \Gamma
	\end{tikzcd}
\]
Note that this must include the normal bundle contribution of the minimal common undegeneration $\Gamma''$ of $\Lambda$ and $\Lambda'$ (in $D_\Gamma$).

More precisely, the algorithm to compute the pullback of $\alpha$ to $\Lambda\cap\Lambda'$ is:
\begin{algorithm}[Pullback]\mbox{}
\begin{description}
\item[Step 1] Compute the common normal bundle of $\alpha$ and $\Lambda'$ in $D_\Gamma$ (\autoref{subsec:NB}).
\item[Step 2] If it is \verb|1|, i.e. the intersection is transversal, perform the pullback to each graph of $\Lambda\cap\Lambda'$ as described above.
\item[Step 3] Otherwise, multiply (in $\CH(\Lambda\cap\Lambda')$) the pullback to each graph of $\Lambda\cap\Lambda'$ with the normal bundle, with ambient $\Gamma''$.
\end{description}
\end{algorithm}

\begin{rem}
Observe that this recursive procedure terminates: indeed, as for the normal bundle calculation, the involved intersections are always performed in an \verb|ambient| stratum of strictly lower dimension, cf. \autoref{nb:terminates}.
\end{rem}

\begin{example}
We continue in the setting of \autoref{topxieval} and \autoref{ex:cnb} in the boundary of $\Omega\cM_2(1,1)$. Recall that \verb|Y.bics[3]| was the compact-type graph with top-level an $\cM_2(2)$. The normal bundle of this graph is simply a $\psi$-class on top-level:
\begin{lstlisting}
sage: N = Y.normal_bundle(((3,),0))
sage: print(N)
Tautological class on Stratum: (1, 1)
with residue conditions: []

-1 * Psi class 1 with exponent 1 on level 0 * Graph ((3,), 0) +
\end{lstlisting}
Here, the V-graph is BIC \verb|0|. It has empty intersection with the BIC \verb|3|. We may still perform the pullback and obtain the \verb|ZERO| class:
\begin{lstlisting}
sage: Y.bics[0]
EmbeddedLevelGraph(LG=LevelGraph([1, 1, 0],[[1], [2], [3, 4, 5, 6]],[(2, 5), (1, 6)],{1: 0, 2: 0, 3: 1, 4: 1, 5: -2, 6: -2},[0, 0, -1],True),dmp={3: (0, 0), 4: (0, 1)},dlevels={0: 0, -1: -1})
sage: print(Y.gen_pullback_taut(N, ((0,),0)))
Tautological class on Stratum: (1, 1)
with residue conditions: []
\end{lstlisting}
Note that we must use \verb|gen_pullback_taut| if we want to pull back an \verb|ELGTautClass| instead of an \verb|AdditiveGenerator|!

Consider now the normal bundle \verb|N2| of the BIC \verb|2| (the other compact type graph). Pulling this back to the BIC \verb|3|, we obtain the normal bundle on the intersection:
\begin{lstlisting}
sage: N2=Y.normal_bundle(((2,),0))
sage: print(Y.gen_pullback_taut(N2, ((3,),0)))
Tautological class on Stratum: (1, 1)
with residue conditions: []

-1 * Psi class 1 with exponent 1 on level 0 * Graph ((2, 3), 0) +
-1 * Psi class 3 with exponent 1 on level 1 * Graph ((2, 3), 0) +
\end{lstlisting}
However, this is \verb|N2| pulled back to \verb|(2, 3)| \emph{inside} BIC \verb|2|:
\begin{lstlisting}
sage: print(Y.gen_pullback_taut(N2, ((2,3),0)))
Tautological class on Stratum: (1, 1)
with residue conditions: []

2 * Psi class 1 with exponent 1 on level 0 * Psi class 3 with exponent 1 on level 1 * Graph ((2, 3), 0) +

sage: print(Y.gen_pullback_taut(N2, ((2,3),0), ((2,),0)))
Tautological class on Stratum: (1, 1)
with residue conditions: []

-1 * Psi class 1 with exponent 1 on level 0 * Graph ((2, 3), 0) +
-1 * Psi class 3 with exponent 1 on level 1 * Graph ((2, 3), 0) +
\end{lstlisting}
The generalised pullback in \verb|Y| is this pullback multiplied with \verb|N2| inside the BIC \verb|2|:
\begin{lstlisting}
sage: print(Y.intersection(Y.gen_pullback_taut(N2, ((2,3),0), ((2,),0)),N2, ((2,),0)))
Tautological class on Stratum: (1, 1)
with residue conditions: []

2 * Psi class 1 with exponent 1 on level 0 * Psi class 3 with exponent 1 on level 1 * Graph ((2, 3), 0) +
\end{lstlisting}
\end{example}

Using \verb|gen_pullback| and \verb|cnb| it is not difficult to implement the  multiplication of arbitrary tautological classes of a \verb|GeneralisedStratum| using the excess intersection formula \cite[Eq. (61)]{strataEC}.

%% file: sec_gluing.tex
\section{Clutching, Splitting and Gluing}
\label{clutching}

Following the philosophy of trading geometric for combinatorial complexity, we wish to transform an expression of classes inside a stratum $X$ into expressions on the levels of the graphs in $X$. Mathematically, this requires a well-behaved level-projection map, cf. \cite[\S 4.2]{strataEC}. In practice, this requires us, given a graph $\Gamma$ and a level $L$ of $\Gamma$, to glue a graph $\Gamma'$ of the generalised stratum $L$ into $\Gamma$ to obtain a graph $\Gamma''$ in $X$.

In light of \autoref{leveltypes} this allows to reduce any calculation to the levels of two and three level graphs together with the combinatorial data of the degeneration graph of $X$.

For our purposes, it will suffice to consider the situation that $\Gamma'$ is a BIC, resulting in a graph in $X$ with one more level than $\Gamma$. However, the implementation of this is rather delicate. The main difficulty stems from a level occuring in several degenerations of the same graph but with different automorphism groups acting on the legs. These automorphisms will, in general, not respect the numbering of the BICs, so extra care must be taken when calculating the profile of $\Gamma''$ from that of $\Gamma$ and $\Gamma'$.

\begin{example}\label{reducible_aut}
Consider the $V$-graph in the boundary of $\Omega\cM_3(4)$, cf. \autoref{fig:vgraph}. 
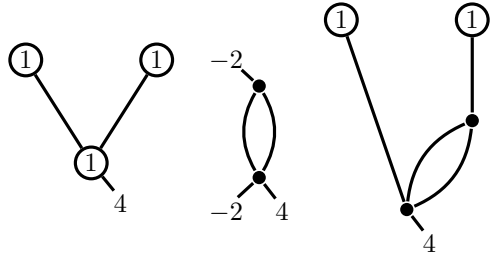
\begin{figure}
\[\begin{tikzpicture}[
		baseline={([yshift=-.5ex]current bounding box.center)},
		scale=2,very thick,
		bend angle=30,
		every loop/.style={very thick},
     		comp/.style={circle,fill,black,,inner sep=0pt,minimum size=5pt},
		order bottom left/.style={pos=.05,left,font=\tiny},
		order top left/.style={pos=.85,left,font=\tiny},
		order bottom right/.style={pos=.05,right,font=\tiny},
		order top right/.style={pos=.85,right,font=\tiny},
		circled number/.style={circle, draw, inner sep=0pt, minimum size=12pt},
		bottom right with distance/.style={below right,text height=10pt},
		bottom left with distance/.style={below left,text height=10pt}]
\begin{scope}[node distance=.5cm]
\node[circled number] (L) [] {$1$}; 
\node [] (T) [right=of L] {};
\node[circled number] (R) [right=of T] {$1$};
\end{scope}
\node[circled number] (B) [below=of T] {$1$}
	edge (L)
	edge (R);
\node [bottom right with distance] (B-4) at (B.south east) {$4$};
\path (B) edge [shorten >=4pt] (B-4.center);
\end{tikzpicture}
\quad
\begin{tikzpicture}[
		baseline={([yshift=-.5ex]current bounding box.center)},
		scale=2,very thick,
		bend angle=30,
		every loop/.style={very thick},
     		comp/.style={circle,fill,black,,inner sep=0pt,minimum size=5pt},
		prong left/.style={pos=.5,left,font=\small},
		prong right/.style={pos=.5,right,font=\small},
		circled number/.style={circle, draw, inner sep=0pt, minimum size=12pt},
		bottom right with distance/.style={below right,text height=10pt},
		bottom left with distance/.style={below left,text height=10pt},
		above right with distance/.style={above right,text height=10pt},
		above left with distance/.style={above left,text height=10pt}]
\node[comp] (T) {};
\node[comp] (B) [below=of T] {}
	edge [bend right] (T)
	edge [bend left] (T);
\node [bottom right with distance] (B-4) at (B.south east) {$4$};
\path (B) edge [shorten >=4pt] (B-4.center);
\node [bottom left with distance] (B-2) at (B.south west) {$-2$};
\path (B) edge [shorten >=6pt] (B-2.center);
\node [above left with distance] (T-22) at (T.north west) {$-2$};
\path (T) edge [shorten >=8pt] (T-22.center);
\end{tikzpicture}
\quad
\begin{tikzpicture}[
		baseline={([yshift=-.5ex]current bounding box.center)},
		scale=2,very thick,
		bend angle=30,
		every loop/.style={very thick},
     		comp/.style={circle,fill,black,,inner sep=0pt,minimum size=5pt},
		order bottom left/.style={pos=.05,left,font=\tiny},
		order top left/.style={pos=.85,left,font=\tiny},
		order bottom right/.style={pos=.05,right,font=\tiny},
		order top right/.style={pos=.85,right,font=\tiny},
		circled number/.style={circle, draw, inner sep=0pt, minimum size=12pt},
		bottom right with distance/.style={below right,text height=10pt},
		bottom left with distance/.style={below left,text height=10pt}]
\begin{scope}[node distance=.5cm]
\node[circled number] (L) [] {$1$}; 
\node [] (T) [right=of L] {};
\node[circled number] (R) [right=of T] {$1$};
\end{scope}
\node[comp] (MR) [below=of R] {}
	edge (R);
\node (M) [below=of T] {};
\node[comp] (B) [below=of M] {}
	edge (L)
	edge [bend left] (MR)
	edge [bend right] (MR);
\node [bottom right with distance] (B-4) at (B.south east) {$4$};
\path (B) edge [shorten >=4pt] (B-4.center);
\end{tikzpicture}
\]
\caption{From left to right: the V-graph in $\Omega\cM_3(4)$; the BICs \texttt{0} and \texttt{2} in the bottom level of V, distinguished only by the labelings of the marked points on top and bottom level; the unique graph in the profile \texttt{(0, 5)}.}
\label{fig:vgraph}
\end{figure}
\begin{figure}
\[
\begin{tikzpicture}[
		baseline={([yshift=-.5ex]current bounding box.center)},
		scale=2,very thick,
		bend angle=30,
		every loop/.style={very thick},
     		comp/.style={circle,fill,black,,inner sep=0pt,minimum size=5pt},
		order bottom left/.style={pos=.05,left,font=\tiny},
		order top left/.style={pos=.85,left,font=\tiny},
		order bottom right/.style={pos=.05,right,font=\tiny},
		order top right/.style={pos=.85,right,font=\tiny},
		circled number/.style={circle, draw, inner sep=0pt, minimum size=12pt},
		bottom right with distance/.style={below right,text height=10pt},
		bottom left with distance/.style={below left,text height=10pt}]
\begin{scope}[node distance=.5cm]
\node[circled number] (L) [] {$1$}; 
\node (T) [right=of L] {};
\node (R) [right=of T] {};
\end{scope}
\node[circled number] (TR) [above=of R] {$1$};
\node[circled number] (B) [below=of T] {$1$}
	edge (L)
	edge (TR);
\node [bottom right with distance] (B-4) at (B.south east) {$4$};
\path (B) edge [shorten >=4pt] (B-4.center);
\end{tikzpicture}
\quad
\begin{tikzpicture}[
		baseline={([yshift=-.5ex]current bounding box.center)},
		scale=2,very thick,
		bend angle=30,
		every loop/.style={very thick},
     		comp/.style={circle,fill,black,,inner sep=0pt,minimum size=5pt},
		order bottom left/.style={pos=.05,left,font=\tiny},
		order top left/.style={pos=.85,left,font=\tiny},
		order bottom right/.style={pos=.05,right,font=\tiny},
		order top right/.style={pos=.85,right,font=\tiny},
		circled number/.style={circle, draw, inner sep=0pt, minimum size=12pt},
		bottom right with distance/.style={below right,text height=10pt},
		bottom left with distance/.style={below left,text height=10pt}]
\begin{scope}[node distance=.5cm]
\node (L) [] {}; 
\node [] (T) [right=of L] {};
\node[circled number] (R) [right=of T] {$1$};
\end{scope}
\node[circled number] (TL) [above=of L] {$1$};
\node[comp] (MR) [below=of R] {}
	edge (R);
\node (M) [below=of T] {};
\node[comp] (B) [below=of M] {}
	edge (TL)
	edge [bend left] (MR)
	edge [bend right] (MR);
\node [bottom right with distance] (B-4) at (B.south east) {$4$};
\path (B) edge [shorten >=4pt] (B-4.center);
\end{tikzpicture}
\quad
\begin{tikzpicture}[
		baseline={([yshift=-.5ex]current bounding box.center)},
		scale=2,very thick,
		bend angle=30,
		every loop/.style={very thick},
     		comp/.style={circle,fill,black,,inner sep=0pt,minimum size=5pt},
		order bottom left/.style={pos=.05,left,font=\tiny},
		order top left/.style={pos=.85,left,font=\tiny},
		order bottom right/.style={pos=.05,right,font=\tiny},
		order top right/.style={pos=.85,right,font=\tiny},
		circled number/.style={circle, draw, inner sep=0pt, minimum size=12pt},
		bottom right with distance/.style={below right,text height=10pt},
		bottom left with distance/.style={below left,text height=10pt}]
\begin{scope}[node distance=.5cm]
\node[circled number] (L) [] {$1$}; 
\node [] (T) [right=of L] {};
\node (R) [right=of T] {};
\end{scope}
\node[circled number] (TR) [above=of R] {$1$};
\node[comp] (MR) [below=of R] {}
	edge (TR);
\node (M) [below=of T] {};
\node[comp] (B) [below=of M] {}
	edge (L)
	edge [bend left] (MR)
	edge [bend right] (MR);
\node [bottom right with distance] (B-4) at (B.south east) {$4$};
\path (B) edge [shorten >=4pt] (B-4.center);
\end{tikzpicture}
\]
\caption{From left to right: the long V \texttt{(1, 0)}; the two graphs in \texttt{(1, 0, 5)} (distinguished by one long edge versus two long edges).}
\label{fig:longvgraph}
\end{figure}
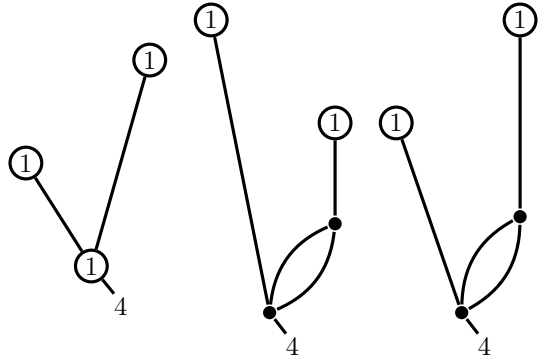
Here it is BIC \verb|0| (cf. \autoref{numbering})
\begin{lstlisting}
sage: X=Stratum((4,))
sage: V=X.bics[0]   # index might change!
sage: V.explain()
LevelGraph embedded into stratum Stratum: (4,)
with residue conditions: []
 with:
On level 0:
* A vertex (number 0) of genus 1
* A vertex (number 1) of genus 1
On level 1:
* A vertex (number 2) of genus 1
The marked points are on level 1.
More precisely, we have:
* Marked point (0, 0) of order 4 on vertex 2 on level 1
Finally, we have 2 edges. More precisely:
* one edge between vertex 0 (on level 0) and vertex 2 (on level 1) with prong 1.
* one edge between vertex 1 (on level 0) and vertex 2 (on level 1) with prong 1.
\end{lstlisting}
In the following discussion, the underlying \verb|LevelGraph| and the leg numbering will be important. Note that the bottom level of \verb|V| has three BICs, each in the shape of a banana, distinguished only by the location of the marked points: for BIC \verb|1|, both are on top, while for \verb|0| and \verb|2| one is on top and one on bottom.
\begin{lstlisting}
sage: V
EmbeddedLevelGraph(LG=LevelGraph([1, 1, 1],[[1], [2], [3, 4, 5]],[(1, 4), (2, 5)],{1: 0, 2: 0, 3: 4, 4: -2, 5: -2},[0, 0, -1],True),dmp={3: (0, 0)},dlevels={0: 0, -1: -1})
sage: V.bot.bics
[EmbeddedLevelGraph(LG=LevelGraph([0, 0],[[1, 2, 3], [4, 5, 6, 7]],[(2, 6), (3, 7)],{1: -2, 2: 0, 3: 0, 4: 4, 5: -2, 6: -2, 7: -2},[0, -1],True),dmp={1: (0, 1), 4: (0, 0), 5: (0, 2)},dlevels={0: 0, -1: -1}),
 EmbeddedLevelGraph(LG=LevelGraph([0, 0],[[1, 2, 3, 4], [5, 6, 7]],[(3, 6), (4, 7)],{1: -2, 2: -2, 3: 1, 4: 1, 5: 4, 6: -3, 7: -3},[0, -1],True),dmp={1: (0, 1), 2: (0, 2), 5: (0, 0)},dlevels={0: 0, -1: -1}),
 EmbeddedLevelGraph(LG=LevelGraph([0, 0],[[1, 2, 3], [4, 5, 6, 7]],[(2, 6), (3, 7)],{1: -2, 2: 0, 3: 0, 4: 4, 5: -2, 6: -2, 7: -2},[0, -1],True),dmp={1: (0, 2), 4: (0, 0), 5: (0, 1)},dlevels={0: 0, -1: -1})]
\end{lstlisting}
Gluing either of \verb|0| or \verb|2| into bottom level results in the same graph (cf. \autoref{fig:vgraph}), as they are exchanged by an automorphism of \verb|V|. This corresponds to both being mapped to the same index by \verb|bot_to_bic| and the corresponding profile \verb|(0, 5)| being irreducible:
\begin{lstlisting}
sage: X.DG.bot_to_bic(0)
{0: 5, 1: 7, 2: 5}
sage: X.lookup((0,5))
[EmbeddedLevelGraph(LG=LevelGraph([1, 0, 1, 0],[[1], [2, 3, 4], [5], [6, 7, 8, 9]],[(1, 4), (2, 7), (3, 8), (5, 9)],{1: 0, 2: 0, 3: 0, 4: -2, 5: 0, 6: 4, 7: -2, 8: -2, 9: -2},[0, -1, 0, -2],True),dmp={6: (0, 0)},dlevels={0: 0, -1: -1, -2: -2})]
\end{lstlisting}
However, the edges are distinguishable in this graph, so e.g. $\psi$-classes might behave differently!

If we do a top-level degeneration of \verb|V| first, resulting in the ``long'' V graph (cf. \autoref{fig:longvgraph}), however, the situation changes. The profile of the long V can be found by inspecting \verb|top_to_bic|.
\begin{lstlisting}
sage: X.DG.top_to_bic(0)
{0: 1, 1: 1}
sage: long_V=X.lookup_graph((1,0)); long_V
EmbeddedLevelGraph(LG=LevelGraph([1, 1, 1],[[1], [2], [3, 4, 5]],[(1, 4), (2, 5)],{1: 0, 2: 0, 3: 4, 4: -2, 5: -2},[0, -1, -2],True),dmp={3: (0, 0)},dlevels={0: 0, -1: -1, -2: -2})
\end{lstlisting}
Indeed, even though the bottom level has not changed, gluing in BIC \verb|0| or \verb|1| now results in different graphs: the profile \verb|(1, 0, 5)| is reducible (cf. \autoref{fig:vgraph}).
\begin{lstlisting}
sage: X.lookup((1,0,5))
[EmbeddedLevelGraph(LG=LevelGraph([1, 1, 0, 0],[[1], [2], [3, 4, 5], [6, 7, 8, 9]],[(1, 5), (3, 7), (4, 8), (2, 9)],{1: 0, 2: 0, 3: 0, 4: 0, 5: -2, 6: 4, 7: -2, 8: -2, 9: -2},[0, -1, -2, -3],True),dmp={6: (0, 0)},dlevels={0: 0, -1: -1, -2: -2, -3: -3}),
 EmbeddedLevelGraph(LG=LevelGraph([1, 1, 0, 0],[[1], [2], [3, 4, 5], [6, 7, 8, 9]],[(1, 5), (3, 7), (4, 8), (2, 9)],{1: 0, 2: 0, 3: 0, 4: 0, 5: -2, 6: 4, 7: -2, 8: -2, 9: -2},[-1, 0, -2, -3],True),dmp={6: (0, 0)},dlevels={0: 0, -1: -1, -2: -2, -3: -3})]
\end{lstlisting}
Therefore, when gluing in a graph of the ``reference'' level (here \verb|X.bics[0].bot|, the bottom level of the bottom BIC), we must keep track of the isomorphism used to identify \verb|delta| of the graph with the ``reference'' BIC.
\end{example}

\subsection{Splitting Dictionaries}

To resolve the ambiguities described above, we always have to  work with a fixed reference stratum. In accordance with \autoref{leveltypes} this will be either a top or bottom level of a BIC or the middle level of a three-level graph. In the first two cases, this level can be found directly from the profile \verb|p|: it is either \verb|X.bics[p[0]].top| or \verb|X.bics[p[-1]].bot|.

However, the length-two profile around an intermediate level could be reducible. The method \verb|three_level_profile_for_level| retrieves the enhanced profile of the three-level graph around a given level of the graph associated to an enhanced profile.

To extract all the data we need from a graph to glue a BIC into one of its levels, we use \emph{splitting dictionaries}. The easiest way to create a splitting dictionary is via \verb|splitting_info_at_level|. We illustrate this for the graph \verb|(5,0)| in \verb|X=Stratum((4,))|, in the situation of \autoref{fig:vgraph}, see \autoref{reducible_aut}:
\begin{lstlisting}
sage: d, leg_dict, L = X.splitting_info_at_level(((0,5),0), 1)
sage: d
{'X': GeneralisedStratum(sig_list=[Signature((4,))],res_cond=[]),
 'top': EmbeddedLevelGraph(LG=LevelGraph([1, 1],[[1], [5]],[],{1: 0, 5: 0},[0, 0],True),dmp={5: (0, 0), 1: (1, 0)},dlevels={0: 0}),
 'bottom': EmbeddedLevelGraph(LG=LevelGraph([0],[[6, 7, 8, 9]],[],{6: 4, 7: -2, 8: -2, 9: -2},[-2],True),dmp={6: (0, 0), 9: (0, 3), 8: (0, 1), 7: (0, 2)},dlevels={-2: -2}),
 'middle': LevelStratum(sig_list=[Signature((0, 0, -2))],res_cond=[[(0, 2)]],leg_dict={2: (0, 0), 3: (0, 1), 4: (0, 2)}),
 'emb_dict_top': {},
 'emb_dict_mid': {},
 'emb_dict_bot': {(0, 0): (0, 0)},
 'clutch_dict': {(1, 0): (0, 2)},
 'clutch_dict_lower': {(0, 1): (0, 2), (0, 0): (0, 1)},
 'clutch_dict_long': {(0, 0): (0, 3)}}
sage: leg_dict
{4: (0, 2), 2: (0, 1), 3: (0, 0)}
sage: L
LevelStratum(sig_list=[Signature((0, 0, -2))],res_cond=[[(0, 2)]],leg_dict={2: (0, 0), 3: (0, 1), 4: (0, 2)})
\end{lstlisting}
As this illustrates, \verb|splitting_info_at_level| actually returns a \verb|tuple| consisting of a splitting dictionary, a dictionary \verb|leg_dict| of legs and a \verb|LevelStratum|. The \verb|LevelStratum| is the standardised level in the sense of \autoref{leveltypes} and will in general differ from the one returned by \verb|level(l)| of the \verb|EmbeddedLevelGraph| associated to the enhanced profile. Similarly, \verb|leg_dict| should be used instead of \verb|dmp| of \verb|level(l)|.

The splitting dictionary uses keywords to encode the relevant data for performing a clutch. More precisely, for an enhanced profile and a level \verb|l|, we obtain a \verb|dict| containing information about the graphs and strata:
\begin{description}
\item[\texttt{X}] the enveloping \verb|GeneralisedStratum|;
\item[\texttt{top}] an \verb|EmbeddedLevelGraph| inside the \verb|top| component of the top BIC of the three-level graph around the level \verb|l|;
\item[\texttt{bottom}] an \verb|EmbeddedLevelGraph| inside the \verb|bot| component of the bottom BIC of the three-level graph around the level \verb|l|;
\item[\texttt{middle}] the \verb|LevelStratum| corresponding to the middle level of the three-level graph around the level \verb|l|;
\end{description}
as well as dictionaries describing the clutching. Note that the ``removed'' edges are now marked points of the strata, so the \verb|dict|s are all in terms of stratum points (cf. \autoref{pointtypes}). More precisely:
\begin{description}
\item[\texttt{clutch\char`_dict}] a \verb|dict| mapping points of \verb|top| stratum to points of \verb|middle| stratum;
\item[\texttt{clutch\char`_dict\char`_lower}] a \verb|dict| mapping points of \verb|middle| stratum to points of \verb|bottom| stratum;
\item[\texttt{clutch\char`_dict\char`_long}] a \verb|dict| mapping points of \verb|top| stratum to points of \verb|bottom| stratum.
\end{description}
Finally, there are three \verb|emb_dict|s, relating the marked points of \verb|X| to the marked points in \verb|top|, \verb|middle| and \verb|bottom|.

Using these, we can simply replace \verb|top|, \verb|middle| or \verb|bottom| by degenerations (inside these strata) and use the information stored in the clutching dictionary to ``safely'' clutch these together.

\begin{example}\label{ex:gluing}
Consider the following situation in the boundary of $\Omega\cM_3(2,1,1)$: We consider the intersection \verb|G| of the asymmetric ``double v''-graph and a ``v''-graph depicted in \autoref{fig:Vs211}. Inspecting \verb|bics|, we find that these are numbered \verb|11| and \verb|23| (cf. \autoref{numbering}). Consequently, \verb|G| has profile \verb|(11, 23)|.
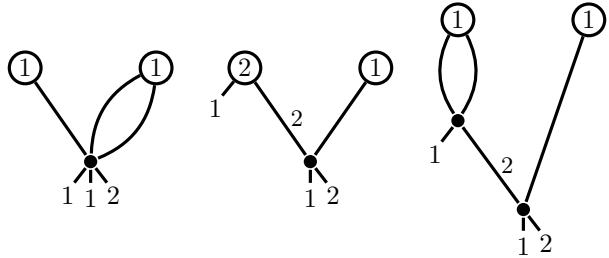
\begin{figure}
\[\begin{tikzpicture}[
		baseline={([yshift=-.5ex]current bounding box.center)},
		scale=2,very thick,
		bend angle=30,
		every loop/.style={very thick},
     		comp/.style={circle,fill,black,,inner sep=0pt,minimum size=5pt},
		order bottom left/.style={pos=.05,left,font=\tiny},
		order top left/.style={pos=.85,left,font=\tiny},
		order bottom right/.style={pos=.05,right,font=\tiny},
		order top right/.style={pos=.85,right,font=\tiny},
		circled number/.style={circle, draw, inner sep=0pt, minimum size=12pt},
		bottom right with distance/.style={below right,text height=10pt},
		bottom left with distance/.style={below left,text height=10pt}]
\begin{scope}[node distance=.5cm]
\node[circled number] (L) [] {$1$}; 
\node [] (T) [right=of L] {};
\node[circled number] (R) [right=of T] {$1$};
\end{scope}
\node[comp] (B) [below=of T] {}
	edge (L)
	edge [bend left] (R)
	edge [bend right] (R);
\node [bottom right with distance] (B-2) at (B.south east) {$2$};
\path (B) edge [shorten >=4pt] (B-2.center);
\node [below,text height=10pt] (B-11) at (B.south) {$1$};
\path (B) edge [shorten >=4pt] (B-11.center);
\node [bottom left with distance] (B-12) at (B.south west) {$1$};
\path (B) edge [shorten >=4pt] (B-12.center);
\end{tikzpicture}
\quad
\begin{tikzpicture}[
		baseline={([yshift=-.5ex]current bounding box.center)},
		scale=2,very thick,
		bend angle=30,
		every loop/.style={very thick},
     		comp/.style={circle,fill,black,,inner sep=0pt,minimum size=5pt},
		order bottom left/.style={pos=.05,left,font=\tiny},
		order top left/.style={pos=.85,left,font=\tiny},
		order bottom right/.style={pos=.05,right,font=\tiny},
		order top right/.style={pos=.85,right,font=\tiny},
		prong left/.style={pos=.5,left,font=\small},
		prong right/.style={pos=.5,right,font=\small},
		circled number/.style={circle, draw, inner sep=0pt, minimum size=12pt},
		bottom right with distance/.style={below right,text height=10pt},
		bottom left with distance/.style={below left,text height=10pt}]
\begin{scope}[node distance=.5cm]
\node[circled number] (L) [] {$2$}; 
\node [] (T) [right=of L] {};
\node[circled number] (R) [right=of T] {$1$};
\end{scope}
\node[comp] (B) [below=of T] {}
	edge node [prong right] {$2$} (L)
	edge (R);
\node [bottom right with distance] (B-2) at (B.south east) {$2$};
\path (B) edge [shorten >=4pt] (B-2.center);
\node [below,text height=10pt] (B-11) at (B.south) {$1$};
\path (B) edge [shorten >=4pt] (B-11.center);
\node [bottom left with distance] (L-12) at (L.south west) {$1$};
\path (L) edge [shorten >=4pt] (L-12.center);
\end{tikzpicture}
\quad
\begin{tikzpicture}[
		baseline={([yshift=-.5ex]current bounding box.center)},
		scale=2,very thick,
		bend angle=30,
		every loop/.style={very thick},
     		comp/.style={circle,fill,black,,inner sep=0pt,minimum size=5pt},
		order bottom left/.style={pos=.05,left,font=\tiny},
		order top left/.style={pos=.85,left,font=\tiny},
		order bottom right/.style={pos=.05,right,font=\tiny},
		order top right/.style={pos=.85,right,font=\tiny},
		prong left/.style={pos=.5,left,font=\small},
		prong right/.style={pos=.5,right,font=\small},
		circled number/.style={circle, draw, inner sep=0pt, minimum size=12pt},
		bottom right with distance/.style={below right,text height=10pt},
		bottom left with distance/.style={below left,text height=10pt}]
\begin{scope}[node distance=.5cm]
\node[circled number] (L1) [] {$1$}; 
\node [] (T1) [right=of L] {};
\node[circled number] (R) [right=of T] {$1$};
\end{scope}
\node[comp] (L2) [below=of L1] {}
	edge [bend left] (L1)
	edge [bend right] (L1);
\node [] (T2) [below=of T1] {};
\node[comp] (B) [below=of T2] {}
	edge (R)
	edge node [prong right] {$2$} (L2);
\node [bottom right with distance] (B-2) at (B.south east) {$2$};
\path (B) edge [shorten >=4pt] (B-2.center);
\node [below,text height=10pt] (B-11) at (B.south) {$1$};
\path (B) edge [shorten >=4pt] (B-11.center);
\node [bottom left with distance] (L2-12) at (L2.south west) {$1$};
\path (L2) edge [shorten >=4pt] (L2-12.center);
\end{tikzpicture}
\]
\caption{From left to right: the ``double v''-graph, \texttt{11}; the ``v''-graph, \texttt{23}; their intersection \texttt{(11, 23)} all in the boundary of $\Omega\cM_3(2,1,1)$. Note that we omit the prongs whenever they are $1$.}
\label{fig:Vs211}
\end{figure}
\begin{lstlisting}
sage: X=Stratum((2,1,1))
sage: G=X.lookup_graph((11, 23))
sage: X.lookup_graph((11,23)).explain()
LevelGraph embedded into stratum Stratum: (2, 1, 1)
with residue conditions: []
 with:
On level 0:
* A vertex (number 0) of genus 1
* A vertex (number 2) of genus 1
On level 1:
* A vertex (number 1) of genus 0
On level 2:
* A vertex (number 3) of genus 0
The marked points are on levels 1 and 2.
More precisely, we have:
* Marked point (0, 1) of order 1 on vertex 1 on level 1
* Marked point (0, 0) of order 2 on vertex 3 on level 2
* Marked point (0, 2) of order 1 on vertex 3 on level 2
Finally, we have 4 edges. More precisely:
* 2 edges between vertex 0 (on level 0) and vertex 1 (on level 1) with prongs 1 and 1.
* one edge between vertex 1 (on level 1) and vertex 3 (on level 2) with prong 2.
* one edge between vertex 2 (on level 0) and vertex 3 (on level 2) with prong 1.
\end{lstlisting}
However, comparing the top level of \verb|G| to the \verb|top| of the top BIC, \verb|11|, we find that the components have been switched:
\begin{lstlisting}
sage: G.level(0)
LevelStratum(sig_list=[Signature((0, 0)), Signature((0,))],res_cond=[],leg_dict={1: (0, 0), 2: (0, 1), 7: (1, 0)})
sage: X.bics[11].top
LevelStratum(sig_list=[Signature((0,)), Signature((0, 0))],res_cond=[],leg_dict={1: (0, 0), 2: (1, 0), 3: (1, 1)})
\end{lstlisting}
The top-level degenerations of \verb|11| are given by \verb|top_to_bic| (cf. \autoref{subsec:graphsandprofiles}). However, to successfully degenerate via clutching, we have to use the \verb|bics| inside \verb|top| of BIC \verb|11|. This stratum is part of the information given by \verb|splitting_info_at_level|:
\begin{lstlisting}
sage: X.DG.top_to_bic(11)
{0: 16, 1: 12, 2: 12, 3: 31}
sage: d, leg_dict, L = X.splitting_info_at_level(((11,23),0), 0)
sage: L
LevelStratum(sig_list=[Signature((0,)), Signature((0, 0))],res_cond=[],leg_dict={1: (0, 0), 2: (1, 0), 3: (1, 1)})
\end{lstlisting}
Note that it corresponds to \verb|X.bics[11].top|, not \verb|G.level(0)|!

To insert a BIC into the top level, we replace the appropriate entry of the clutching dictionary and feed this to clutch (using Python's \verb|**| operator):
\begin{lstlisting}
sage: d
{'X': GeneralisedStratum(sig_list=[Signature((2, 1, 1))],res_cond=[]),
 'top': EmbeddedLevelGraph(LG=LevelGraph([1, 1],[[1, 2], [7]],[],{1: 0, 2: 0, 7: 0},[0, 0],True),dmp={7: (0, 0), 2: (1, 0), 1: (1, 1)},dlevels={0: 0}),
 'bottom': EmbeddedLevelGraph(LG=LevelGraph([0, 0],[[3, 4, 5, 6], [8, 9, 10, 11]],[(4, 10)],{3: 1, 4: 1, 5: -2, 6: -2, 8: 2, 9: 1, 10: -3, 11: -2},[-1, -2],True),dmp={3: (0, 1), 8: (0, 0), 9: (0, 2), 11: (0, 3), 6: (0, 4), 5: (0, 5)},dlevels={-1: -1, -2: -2}),
 'clutch_dict': {(1, 1): (0, 5), (1, 0): (0, 4), (0, 0): (0, 3)},
 'emb_dict_top': {},
 'emb_dict_bot': {(0, 0): (0, 0), (0, 1): (0, 1), (0, 2): (0, 2)}}
sage: d['top']=L.bics[1]
sage: H=clutch(**d)
\end{lstlisting}
Note that using \verb|G.level(0).bics[1]| would result in an error here, as swapping of the components invalidates the gluing data.

In accordance with \verb|top_to_bic|, we find that the profile of \verb|H| is in fact given by \verb|(12, 11, 23)|. However, this profile is not reducible:
\begin{lstlisting}
sage: H.is_isomorphic(X.lookup_graph((12,11,23), 0))
False
sage: H.is_isomorphic(X.lookup_graph((12,11,23), 1))
True
\end{lstlisting}
Indeed, we obtain the other component by clutching in instead the other BIC of \verb|L| that is mapped to \verb|12| by \verb|top_to_bic|:
\begin{lstlisting}
sage: d['top']=L.bics[2]
sage: HH=clutch(**d)
sage: HH.is_isomorphic(X.lookup_graph((12,11,23),0))
True
\end{lstlisting}
Continuing with this graph, we see that we run into a similar problem when trying to degenerate level \verb|1|: the reference level around level \verb|1| is the middle level of the three-level graph around this level. We find it using \verb|three_level_profile_for_level| (indeed, the profile is irreducible and we need the ``non-standard'' component):
\begin{lstlisting}
sage: X.three_level_profile_for_level(((12,11,23),0),1)
((12, 11), 1)
sage: X.lookup_graph((12,11),1).level(1)
LevelStratum(sig_list=[Signature((0,)), Signature((0, 0, -2))],res_cond=[[(1, 2)]],leg_dict={1: (0, 0), 3: (1, 0), 4: (1, 1), 5: (1, 2)})
\end{lstlisting}
Now, it turns out that this reference level is in fact the level \verb|1| of \verb|HH| but \emph{not} of the reference graph associated to the enhanced profile \verb|((12, 11, 23), 0)| (the fixed representative of the isomorphism class of graphs):
\begin{lstlisting}
sage: HH.level(1)
LevelStratum(sig_list=[Signature((0,)), Signature((0, 0, -2))],res_cond=[[(1, 2)]],leg_dict={1: (0, 0), 3: (1, 0), 4: (1, 1), 5: (1, 2)})
sage: X.lookup_graph((12,11,23),0).level(1)
LevelStratum(sig_list=[Signature((0, 0, -2)), Signature((0,))],res_cond=[[(0, 2)]],leg_dict={2: (0, 0), 3: (0, 1), 4: (0, 2), 5: (1, 0)})
\end{lstlisting}
Again, the components have been switched. However, \verb|splitting_info_at_level| gives the reference level (the middle level of \verb|three_level_profile_for_level|) and this is the one we should use for clutching:
\begin{lstlisting}
sage: d, leg_dict, L = X.splitting_info_at_level(((12,11,23),0), 1)
sage: L
LevelStratum(sig_list=[Signature((0,)), Signature((0, 0, -2))],res_cond=[[(1, 2)]],leg_dict={1: (0, 0), 3: (1, 0), 4: (1, 1), 5: (1, 2)})
\end{lstlisting}
In particular, this agrees with \verb|middle_to_bic|:
\begin{lstlisting}
sage: X.DG.middle_to_bic(((12,11),1))
{0: 16, 1: 31}
sage: d
{'X': GeneralisedStratum(sig_list=[Signature((2, 1, 1))],res_cond=[]),
 'top': EmbeddedLevelGraph(LG=LevelGraph([1],[[1]],[],{1: 0},[0],True),dmp={1: (0, 0)},dlevels={0: 0}),
 'bottom': EmbeddedLevelGraph(LG=LevelGraph([0, 0],[[6, 7, 8, 9], [10, 11, 12, 13]],[(7, 12)],{6: 1, 7: 1, 8: -2, 9: -2, 10: 2, 11: 1, 12: -3, 13: -2},[-2, -3],True),dmp={6: (0, 1), 10: (0, 0), 11: (0, 2), 13: (0, 3), 9: (0, 4), 8: (0, 5)},dlevels={-2: -2, -3: -3}),
 'middle': LevelStratum(sig_list=[Signature((0,)), Signature((0, 0, -2))],res_cond=[[(1, 2)]],leg_dict={1: (0, 0), 3: (1, 0), 4: (1, 1), 5: (1, 2)}),
 'emb_dict_top': {},
 'emb_dict_mid': {},
 'emb_dict_bot': {(0, 0): (0, 0), (0, 1): (0, 1), (0, 2): (0, 2)},
 'clutch_dict': {(0, 0): (1, 2)},
 'clutch_dict_lower': {(1, 1): (0, 5), (1, 0): (0, 4), (0, 0): (0, 3)},
 'clutch_dict_long': {}}
sage: d['middle'] = L.bics[0]
sage: HHH=clutch(**d)
sage: HHH.is_isomorphic(X.lookup_graph((12,16,11,23)))
True
sage: d['middle'] = L.bics[1]
sage: HHH=clutch(**d)
sage: HHH.is_isomorphic(X.lookup_graph((12,31,11,23)))
True
sage: HHH.is_isomorphic(X.lookup_graph((12,16,11,23)))
False
\end{lstlisting}
Note that in the above example, the top strata differed in an obvious manner and clutching immediately raised an error. However, more subtle renumberings of BICs can occur:
\begin{lstlisting}
sage: X.lookup_graph((15,33)).level(0)
LevelStratum(sig_list=[Signature((0, 2))],res_cond=[],leg_dict={1: (0, 0), 2: (0, 1)})
sage: X.bics[15].top
LevelStratum(sig_list=[Signature((2, 0))],res_cond=[],leg_dict={1: (0, 0), 2: (0, 1)})
sage: X.lookup_graph((15,33)).level(0).bics == X.bics[15].top.bics
False
sage: X.lookup_graph((15,33)).level(0).bics
[EmbeddedLevelGraph(LG=LevelGraph([1, 1],[[1], [2, 3, 4]],[(1, 4)],{1: 0, 2: 2, 3: 0, 4: -2},[0, -1],True),dmp={2: (0, 1), 3: (0, 0)},dlevels={0: 0, -1: -1}),
 EmbeddedLevelGraph(LG=LevelGraph([1, 0],[[1, 2, 3], [4, 5, 6]],[(2, 5), (3, 6)],{1: 0, 2: 0, 3: 0, 4: 2, 5: -2, 6: -2},[0, -1],True),dmp={1: (0, 0), 4: (0, 1)},dlevels={0: 0, -1: -1}),
 EmbeddedLevelGraph(LG=LevelGraph([1, 1],[[1, 2], [3, 4]],[(2, 4)],{1: 0, 2: 0, 3: 2, 4: -2},[0, -1],True),dmp={1: (0, 0), 3: (0, 1)},dlevels={0: 0, -1: -1}),
 EmbeddedLevelGraph(LG=LevelGraph([2, 0],[[1], [2, 3, 4]],[(1, 4)],{1: 2, 2: 2, 3: 0, 4: -4},[0, -1],True),dmp={2: (0, 1), 3: (0, 0)},dlevels={0: 0, -1: -1}),
 EmbeddedLevelGraph(LG=LevelGraph([1, 0],[[1, 2], [3, 4, 5, 6]],[(1, 5), (2, 6)],{1: 0, 2: 0, 3: 2, 4: 0, 5: -2, 6: -2},[0, -1],True),dmp={3: (0, 1), 4: (0, 0)},dlevels={0: 0, -1: -1})]
sage: X.bics[15].top.bics
[EmbeddedLevelGraph(LG=LevelGraph([1, 0],[[1, 2, 3], [4, 5, 6]],[(2, 5), (3, 6)],{1: 0, 2: 0, 3: 0, 4: 2, 5: -2, 6: -2},[0, -1],True),dmp={1: (0, 1), 4: (0, 0)},dlevels={0: 0, -1: -1}),
 EmbeddedLevelGraph(LG=LevelGraph([1, 1],[[1], [2, 3, 4]],[(1, 4)],{1: 0, 2: 2, 3: 0, 4: -2},[0, -1],True),dmp={2: (0, 0), 3: (0, 1)},dlevels={0: 0, -1: -1}),
 EmbeddedLevelGraph(LG=LevelGraph([1, 0],[[1, 2], [3, 4, 5, 6]],[(1, 5), (2, 6)],{1: 0, 2: 0, 3: 2, 4: 0, 5: -2, 6: -2},[0, -1],True),dmp={3: (0, 0), 4: (0, 1)},dlevels={0: 0, -1: -1}),
 EmbeddedLevelGraph(LG=LevelGraph([2, 0],[[1], [2, 3, 4]],[(1, 4)],{1: 2, 2: 2, 3: 0, 4: -4},[0, -1],True),dmp={2: (0, 0), 3: (0, 1)},dlevels={0: 0, -1: -1}),
 EmbeddedLevelGraph(LG=LevelGraph([1, 1],[[1, 2], [3, 4]],[(2, 4)],{1: 0, 2: 0, 3: 2, 4: -2},[0, -1],True),dmp={1: (0, 1), 3: (0, 0)},dlevels={0: 0, -1: -1})]
\end{lstlisting}
In this case, the graphs clutch without error and the mistake is much harder to detect!
\begin{lstlisting}
sage: d, leg_dict, L = X.splitting_info_at_level(((15,33),0), 0)
sage: d
{'X': GeneralisedStratum(sig_list=[Signature((2, 1, 1))],res_cond=[]),
 'top': EmbeddedLevelGraph(LG=LevelGraph([2],[[1, 2]],[],{1: 0, 2: 2},[0],True),dmp={1: (0, 1), 2: (0, 0)},dlevels={0: 0}),
 'bottom': EmbeddedLevelGraph(LG=LevelGraph([0, 0],[[3, 4, 5, 6], [7, 8, 9]],[(5, 8)],{3: 1, 4: 1, 5: 0, 6: -4, 7: 2, 8: -2, 9: -2},[-1, -2],True),dmp={3: (0, 1), 4: (0, 2), 7: (0, 0), 9: (0, 4), 6: (0, 3)},dlevels={-1: -1, -2: -2}),
 'clutch_dict': {(0, 0): (0, 3), (0, 1): (0, 4)},
 'emb_dict_top': {},
 'emb_dict_bot': {(0, 0): (0, 0), (0, 1): (0, 1), (0, 2): (0, 2)}}
sage: L
LevelStratum(sig_list=[Signature((2, 0))],res_cond=[],leg_dict={1: (0, 0), 2: (0, 1)})
sage: X.DG.top_to_bic(15)
{0: 32, 1: 12, 2: 31, 3: 26, 4: 31}
sage: d['top'] = X.lookup_graph((15,33)).level(0).bics[0]
sage: clutch(**d).is_isomorphic(X.lookup_graph((32,15,33)))
False
sage: d['top'] = L.bics[0]
sage: clutch(**d).is_isomorphic(X.lookup_graph((32,15,33)))
True
\end{lstlisting}
\end{example}

\subsection{Pullback of Classes on a Level}

Since for our goals it is only necessary to pull back the classes $\xi$ and $\cL$ from a level, we restrict this discussion to the case of codimension-one classes. There are thus two cases to distinguish. Let \verb|G| be a graph with enhanced profile \verb|ep| and \verb|L| be the standardised level stratum at level \verb|l| of \verb|G|. Then a codimension-one class in \verb|L| is either
\begin{itemize}
\item a $\psi$-class on \verb|L| or
\item a BIC in \verb|L|.
\end{itemize}
In the first case, the pullback will consist of $\psi$-classes on the graph \verb|G|, in the second case, the pullback class will be a one-level degeneration of \verb|G|. The pullback of $\xi$ from level \verb|l| is accomplished by \verb|xi_at_level|.

\begin{rem}\label{xiatlevel}
The method \verb|xi_at_level(i,ep)| corresponds to the class $\xi^{[i]}_{\Gamma}$ where $\Gamma$ is the level graph associated to the enhanced profile \verb|ep|, see \cite[Prop. 4.7]{strataEC} for details. As the class $\xi^{[i]}_{\Gamma}$ is an element of $\CH(D_\Gamma)$, care must be taken when multiplying, cf. \autoref{topxieval}.
\end{rem}

For example, for any Stratum \verb|X|, multiplying any graph with $\xi$ yields an equivalent class to pulling back $\xi$ from top level. We can check this for one-dimensional graphs by evaluating:
\begin{lstlisting}
sage: all((X.xi*X.taut_from_graph(*ep)).evaluate() == X.xi_at_level(0, ep).evaluate() for ep in X.enhanced_profiles_of_length(X.dim()-1))
True
\end{lstlisting}

\begin{rem}
Note that it is essential for \verb|L| to be the standardised level in the sense of \autoref{leveltypes} and obtained via \verb|splitting_info_at_level| (not \verb|G.level(l)|!), to obtain the correct class.
\end{rem}

Using \verb|splitting_info_at_level|, the pullback of a $\psi$-class is straight-forward, as we have access to \verb|L.smooth_LG.dmp| and \verb|leg_dict| to find the correct leg number on \verb|G|.

Pulling back the class of a BIC is a slightly more involved, as we need to determine the enhanced profile of the one-level degeneration of \verb|G|:
\begin{algorithm}[Gluing in a BIC]\mbox{}
\begin{description}
\item[Step 1] Determine the new profile. This is given by the degeneration graph via \verb|X.DG.top_to_bic|, \verb|X.DG.bot_to_bic| or \verb|X.DG.middle_to_bic|, depending on the location of \verb|L|.
\item[Step 2] Determine the graph. For this, we replace \verb|L| by the BIC in the splitting dictionary as illustrated in \autoref{ex:gluing} and use \verb|clutch| to build the \verb|EmbeddedLevelGraph| in \verb|X|.
\item[Step 3] We find the enhanced profile by locating the isomorphism class of the clutched graph inside the new profile.
\end{description}
\end{algorithm}
Moreover, we need to weigh the pulled back class with the contribution from
comparison of multiplicities in the level projections as given in
\cite[Prop.~4.7]{strataEC}. 
The correction factor is the product of the edge contribution and the automorphism contribution. The edge contribution is the quotient of the \verb|ell| of the BIC of \verb|X| that extended the profile and the \verb|ell| of the BIC of \verb|L| that was inserted. The automorphism factor is the quotient of the number of automorphisms of the glued graph and the product of the number of automorphisms of \verb|G| and the BIC of~\verb|L|.

\subsubsection*{Leg choice}
As the class of $\xi$ is always implemented via Sauvaget's relation, it requires a choice of leg. By default, the one giving the shortest expression is chosen (i.e. the one appearing on bottom level for the fewest BICs in \verb|L|). Using the optional argument \verb|leg|, we may specify a \emph{leg of} \verb|G| that is to be used. \verb|xi_at_level| will raise a \verb|ValueError| if the leg is not found on level \verb|l|.

\subsubsection*{Empty Profile}
We may apply \verb|xi_at_level| to the empty profile \verb|((), 0)| and obtain $\xi$ on \verb|X|. Note, however, that while \verb|xi_with_leg| requires a stratum point (cf. \autoref{pointtypes}), the optional \verb|leg| argument of \verb|xi_at_level| requires a leg of \verb|X.smooth_LG|:
\begin{lstlisting}
sage: X=Stratum((2,-2))
sage: X.smooth_LG
EmbeddedLevelGraph(LG=LevelGraph([1],[[1, 2]],[],{1: 2, 2: -2},[0],True),dmp={1: (0, 0), 2: (0, 1)},dlevels={0: 0})
sage: X.xi == X.xi_at_level(0, ((),0))
True
sage: X.xi_with_leg((0,0)) == X.xi_at_level(0, ((),0), leg=1)
True
sage: X.xi_with_leg((0,1)) == X.xi_at_level(0, ((),0), leg=2)
True
\end{lstlisting}

\subsection{Splitting Graphs}

Let \verb|G| be the \verb|EmbeddedLevelGraph| associated to the enhanced profile \verb|enh_profile|, written as \verb|(p, i)| inside the \verb|GeneralisedStratum|~\verb|X|. To construct the clutching dictionaries used above, we must, given a level \verb|l|, realise the subgraph of \verb|G| above level \verb|l| inside the \verb|top| level of the BIC \verb|p[l-1]|. Denote by \verb|L| the standardised level \verb|l| of \verb|G|.

Once this has been accomplished, the splitting essentially reduces to the case of a BIC (if \verb|l| is \verb|0|, i.e. top level, or \verb|len(p)|, i.e. bottom level) or that of a three-level graph split around the middle level.

The method \verb|splitting_info_at_level| performs this distinction and serves as a wrapper for the low-level splitting functions, returning all the information needed for clutching.

\subsubsection*{Extracting subgraphs}
The extraction of the subgraph is accomplished by the method \verb|sub_graph_from_level| and yields a ``true'' subgraph, i.e. the names of the vertices and legs are the same as in \verb|G|.

To get the embedding into the appropriate level, we have to fix an undegeneration to the appropriate BIC via \verb|explicit_leg_maps|, cf. \autoref{subsec:degenerations}.
Note that the marked points of \verb|L| correspond to the marked points of \verb|X| above (resp. below) level \verb|l| on \verb|G| and to the top (resp. bottom) legs of the edges that cross level \verb|l|.

Using \verb|sub_graph_from_level|, it suffices to retrieve the splitting dictionary of the appropriate BIC or three-level graph. The splitting dictionary of the level graph~\verb|G| is then obtained by replacing the appropriate strata by the subgraphs given by \verb|sub_graph_from_level|.

\subsubsection*{Splitting BICs}

Splitting BICs is implemented, for a graph with exactly two levels, in the method \verb|split| of \verb|EmbeddedLevelGraph|.

For a BIC, splitting happens in several steps:
\begin{algorithm}[Splitting a BIC]\mbox{}
\begin{description}
\item[Step 1] Extract \verb|top| and \verb|bot| (via \verb|level|).
\item[Step 2] Construct \verb|emb_top| and \verb|emb_bot| by combining \verb|dmp| with the \verb|leg_dicts| of \verb|top| and \verb|bot|.
\item[Step 3] Save the gluing information from the cut edges. Because this is a BIC, all edges are cut in this process.
\end{description}
\end{algorithm}

This yields the splitting dictionary.

\subsubsection*{Splitting three-level graphs}

A three-level graph is determined by its enhanced profile. Consequently, while the splitting of BICs is a method of \verb|EmbeddedLevelGraph|, the corresponding method for three-level graphs, \verb|doublesplit|, is a method of \verb|GeneralisedStratum|.

This is important, because when splitting around a level, everything should be embedded into \verb|top| and \verb|bot| of the two BICs surrounding the level, the only strata we can control. In particular, we need to split into \verb|top| and \verb|bot| of the BICs \verb|p[0]| and \verb|p[1]| and not into \verb|level(0)| and \verb|level(2)| of the three-level graph, as these could differ by a non-trivial automorphism! See \autoref{ex:gluing}.

We therefore have to work with \verb|explicit_leg_maps| to fix undegeneration maps to \verb|p[0]| and \verb|p[1]|. This gives a map from points in the three-level graphs to points in the BIC and we can compose this with the embedding of \verb|top| and \verb|bot| to construct the embedding maps for the splitting dictionary.

As this is a three-level graph, again all edges are cut in this process. In contrast to the BIC case, we need to distinguish edges from top to middle, middle to bottom and top to bottom (long edges) here.

All this is stored in the splitting dictionary and returned by \verb|doublesplit|.

\subsection{Clutching}

The \verb|clutch| method of the \verb|stratatautring| module takes a splitting dictionary and produces from it an \verb|EmbeddedLevelGraph|.

Note that the case clutching two graphs (e.g. degenerating top or bottom level) and clutching three graphs (e.g. degenerating an ``interior'' level) must be distinguished; one cannot perform clutch twice or recursively, as the intermediate graph is not embedded into the same stratum.

Note that clutching is performed on graphs, so the top, middle and bottom components are converted to their respective \verb|smooth_LG| if one of them is stratum. The clutch is then performed in several steps:
\begin{algorithm}[Clutching]\mbox{}
\begin{description}
\item[Step 1] Unite the vertices, renumbering the levels appropriately.
\item[Step 2] Unite the legs, renumbering appropriately. Here we have to keep track of the renumbering for the new \verb|dmp| and to apply the clutching information. We also insert the renumbered (old) edges and use the clutching dictionaries to store the legs that are to be identified.
\item[Step 3] We create the new edges from the information gathered in Step 2.
\item[Step 4] We use this data to create an \verb|EmbeddedLevelGraph| that we return.
\end{description}
\end{algorithm}
Using Python's \verb|**| operator, we can feed a splitting dictionary directly into \verb|clutch|:
\begin{lstlisting}
sage: X=GeneralisedStratum([Signature((1,1))])
sage: assert clutch(**X.bics[1].split()).is_isomorphic(X.bics[1])
sage: assert all(clutch(**B.split()).is_isomorphic(B) for B in X.bics)
\end{lstlisting}
Of course, because of the renumbering, we cannot assume the graphs to be the same after splitting and clutching. They are, however, isomorphic.
The same works for three-level graphs and \verb|doublesplit|, also for a more complicated stratum:
\begin{lstlisting}
sage: X=GeneralisedStratum([Signature((1,1))])
sage: assert all(X.lookup_graph(*ep).is_isomorphic(clutch(**X.doublesplit(ep))) for ep in X.enhanced_profiles_of_length(2))
sage: X=GeneralisedStratum([Signature((2,2,-2))])
sage: assert all(X.lookup_graph(*ep).is_isomorphic(clutch(**X.doublesplit(ep))) for ep in X.enhanced_profiles_of_length(2))
\end{lstlisting}
In particular, there are three-level graphs with long edges in both these strata.

%% file: sec_caching.tex
\section{Caching}
\label{sec:caching}

The boundary strata grow in size very quickly. Even for holomorphic strata in genus $3$, class calculations would not be possible without extensive caching.

The downside of this is that \verb|diffstrata| has quite an extensive memory footprint; there is surely still room for much optimisation.

\subsection{Graphs}

Working with explicit \verb|LevelGraph|s is painfully slow. The main achievement of 
\autoref{subsec:graphsandprofiles} was to associate to each \verb|EmbeddedLevelGraph| an enhanced profile, i.e. a \verb|tuple| (of \verb|tuple|s) of integers for each isomorphism class. Ideally, we always refer to a graph by its enhanced profile.

\begin{rem}
Note that it is important to use \verb|tuple|s and not \verb|list|s for (enhanced) profiles, as \verb|tuple|s are immutable and may thus be used as arguments of cached functions and as keys of dictionaries.
\end{rem}

However, for certain degeneration questions as well as level extraction, we do need concrete representations of the graphs as well as explicit isomorphisms and leg maps. For this, we store, for each enhanced profile, a \emph{reference graph}. More precisely, for a \verb|GeneralisedStratum| \verb|X| and a profile \verb|p|, \verb|X.lookup(p)| will, on first call, generate all graphs in the profile \verb|p| (cf. \autoref{subsec:graphsandprofiles}) and stores this \verb|list| in the dictionary \verb|X._lookup| with key \verb|p|. On subsequent calls, this dictionary is used, so that the clutching is only done once and the \emph{same} \verb|EmbeddedLevelGraph| is returned on every lookup. Moreover, for profiles of length $1$, the entries of \verb|bics| are used:
\begin{lstlisting}
sage: G=X.bics[0]
sage: G is X.lookup_graph((0,),0)
True
sage: H=X.lookup_graph((1,0),0)
sage: H is X.lookup_graph((1,0),0)
True
\end{lstlisting}
Moreover, each \verb|EmbeddedLevelGraph| \verb|G| has a \verb|list| \verb|G._levels| that is filled with the extracted levels on first call of \verb|G.level| (and these are subsequently reused). Therefore, also the individual levels of the graph associated to an enhanced profile will not be regenerated. Also, other intrinsic values, such as the number of automorphisms, are stored on first computation and then retrieved.

Finally, \verb|AdditiveGenerator|s are considered immutable and are even hashable, i.e. using them, e.g., as keys in dictionaries works:
\begin{lstlisting}
sage: a=X.additive_generator(((0,),0))
sage: {a : 1}
{AdditiveGenerator(X=GeneralisedStratum(sig_list=[Signature((2,))],res_cond=[]),enh_profile=((0,), 0),leg_dict={}): 1}
\end{lstlisting}
By contrast, \verb|ELGTautClass|es are mutable (their \verb|psi_list|s can and will be changed, for example by \verb|reduce|) and therefore may \emph{not} be used as keys.

This allows any method whose arguments consist \emph{only} of enhanced profiles and \verb|AdditiveGenerator|s to be cached (using \verb|sage|'s \verb|@cached_method| decorator). This is a key reason for splitting all methods in \autoref{sec:mult} into operations involving only these objects (instead of working only with \verb|ELGTautClass|es of \verb|EmbeddedLevelGraph|s).

Also, note that \verb|AdditiveGenerator|s should always be created and used via \verb|X.additive_generator| as this stores them in the \verb|_AGs| dictionary of \verb|X| and allows them to be reused (instead of being created newly on each call): 
\begin{lstlisting}
sage: a=X.additive_generator(((0,),0))
sage: a is X.additive_generator(((0,),0))
True
sage: a is AdditiveGenerator(X, ((0,),0))
False
\end{lstlisting}

This allows computations in strata of genus $3$ and $4$ in feasible time. However, the memory footprint is considerable: already in genus $3$, the larger strata use about $20$GB, while in genus $4$ already more than a TB is required.

\subsection{Files and Values}
\label{subsec:files}

As described in \autoref{subsec:eval}, \verb|diffstrata| uses the package \verb|admcycles| to evaluate top-degree \verb|ELGTautClass|es. As the computations of \verb|admcycles| are also very involved, we cache every use and, in fact, (attempt to) write any computed value into a local file that is automatically (attempted to be) reused.

More precisely, given the signature \verb|sig| of a stratum (note that \verb|admcycles| works only with connected strata without residue conditions, cf. \autoref{subsec:eval}) and
a $\psi$-polynomial \verb|psis| (as a \verb|dict| mapping points of the stratum to exponents), the method \verb|adm_evaluate| uses \verb|adm_key| to compute a key consisting of the signature and \verb|psis| transformed into a \verb|tuple|. To avoid needless recomputations, the signature is sorted and \verb|psis| renumbered accordingly:
\begin{lstlisting}
sage: from admcycles.diffstrata.levelstratum import adm_key
sage: adm_key((2,-2), {1: 2, 2: 1})
((-2, 2), ((1, 1), (2, 2)))
\end{lstlisting}
Next, we check if the cache dictionary exists and if not, we attempt to load it from the file \verb|adm_evals.sobj| (see below for details). If the key exists in the cache, we return its value, otherwise we use \verb|admcycles| to compute the value, store it in the cache and write this back into the file.

Similarly, whenever we evaluate a top-power of $\xi$, we save this to file, as the Euler characteristic can be computed purely from the degeneration graph and this information (cf. \autoref{subsec:ec}). 

More precisely, for a \verb|GeneralisedStratum| \verb|X|, an enhanced profile and a level~$l$, the method \verb|X.top_xi_at_level| computes the evaluation of the top-power of $\xi^{[l]}_\Gamma$ on the graph $\Gamma$ associated to the enhanced profile, i.e. \verb|X.xi_at_level_pow| as in \autoref{topxieval}. These numbers are stored in a cache dictionary that is synchronised with the file \verb|top_xis.sobj| as described above.

The $\xi$-cache uses \verb|LevelStratum|'s method \verb|dict_key| to compute a key, again with the aim of performing as few evaluations as necessary: 
the list of components is sorted, as is the signature of each component; then the residue conditions are renumbered appropriately and also sorted. Finally, everything is converted to a nested \verb|tuple|.
\begin{example}
We illustrate the \verb|dict_key|s in the stratum $\Omega\cM_{2}(1,1)$. For the V-graph, the top stratum is disconnected and produces:
\begin{lstlisting}
sage: print(VT)
Product of Strata:
Signature((0,))
Signature((0,))
with residue conditions:
dimension: 3
leg dictionary: {1: (0, 0), 2: (1, 0)}
leg orbits: [[(1, 0), (0, 0)]]

sage: VT.dict_key()
(((0,), (0,)), ())
\end{lstlisting}
The bottom stratum has residue conditions:
\begin{lstlisting}
sage: print(VB)
Stratum: Signature((1, 1, -2, -2))
with residue conditions: [(0, 3)] [(0, 2)]
dimension: 0
leg dictionary: {3: (0, 0), 4: (0, 1), 5: (0, 2), 6: (0, 3)}
leg orbits: [[(0, 0)], [(0, 1)], [(0, 3), (0, 2)]]

sage: VB.dict_key()
(((-2, -2, 1, 1),), (((0, 0),), ((0, 1),)))
\end{lstlisting}
By contrast, the bottom level of the banana graph produces:
\begin{lstlisting}
sage: print(BB)
Stratum: Signature((1, 1, -2, -2))
with residue conditions: [(0, 2), (0, 3)]
dimension: 1
leg dictionary: {3: (0, 0), 4: (0, 1), 5: (0, 2), 6: (0, 3)}
leg orbits: [[(0, 0)], [(0, 1)], [(0, 3), (0, 2)]]

sage: BB.dict_key()
(((-2, -2, 1, 1),), (((0, 0), (0, 1)),))
\end{lstlisting}
\end{example}

\subsubsection*{Importing and Exporting Values}

As described above, any value computed with \verb|adm_evaluate| is cached and synchronised with the file \verb|adm_evals.sobj| (actually the \verb|global| variable \verb|FILENAME|, which is set to \verb|adm_evals.sobj| by default). This is accomplished by the method \verb|load_adm_evals|.

The \verb|diffstrata| method \verb|import_adm_evals| takes a filename as an argument and attempts to update the cache dictionary from \verb|load_adm_evals| with the dictionary read from this file. The result is written immediately to \verb|FILENAME|.

Similarly, the $\xi$-cache is synchronised with \verb|XI_FILENAME|, which is \verb|top_xis.sobj| by default, the loading is handled by \verb|load_xis| and these can be imported using \verb|import_top_xis|.

For file handling, \verb|sage|'s methods \verb|save| and \verb|load| are used, i.e. the files are stored in the \emph{current working directory}.
However, \verb|admcycles| versions \emph{after} \verb|v1.1| include the modules \verb|adm_eval_cache| and \verb|xi_cache| that automatically initialise the cache with a large number of values for low-genus strata.
The caching is then controlled by the more versatile \verb|cache| module and the files are stored in \verb|DOT_SAGE|.

\begin{example}\label{ex:export}
Assume we have a computed only the top $\xi$ of the minimal stratum in genus $3$, our $\xi$-cache would look like this and we may export it:
\begin{lstlisting}
sage: my_xi
{(((4,),), ()): 305/580608}
sage: save(my_xi, 'my_xi.sobj')
\end{lstlisting}
We may import this anywhere (and check):
\begin{lstlisting}
sage: import_top_xis('my_xi.sobj')
sage: load_xis()
{(((4,),), ()): 305/580608}
\end{lstlisting}
Timing this calculation confirms that the cached value is being used:
\begin{lstlisting}
sage: X=Stratum((4,))
sage: %
CPU times: user 14.9 ms, sys: 970 µs, total: 15.9 ms
Wall time: 15.5 ms
305/580608
\end{lstlisting}
\end{example}

\subsubsection*{Printing Values}

As the generation of keys for strata and $\psi$-polynomials described above hinder the readability of the content of the cache, \verb|diffstrata| includes the methods \verb|print_top_xis| and \verb|print_adm_evals| to print these values in a more human-readable format.

\begin{example}
Continuing \autoref{ex:export} from above, assume we are importing a second file and want to confirm that our previous value was not deleted:
\begin{lstlisting}
sage: import_top_xis('my_other_xi.sobj')
sage: print_top_xis()
Stratum            | Residue Conditions           | xi^dim
----------------------------------------------------------------
(-2, -2, 1, 1)     | [[(0, 0)], [(0, 1)]]         | 1
(4,)               | ()                           | 305/580608
\end{lstlisting}
We can also pass the \verb|dict| to be printed as an argument, e.g. to inspect a file before loading:
\begin{lstlisting}
sage: print_top_xis(load('my_xi.sobj'))
Stratum            | Residue Conditions           | xi^dim
----------------------------------------------------------------
(4,)               | ()                           | 305/580608
sage: print_top_xis(load('my_other_xi.sobj'))
Stratum            | Residue Conditions           | xi^dim
----------------------------------------------------------------
(-2, -2, 1, 1)     | [[(0, 0)], [(0, 1)]]         | 1
\end{lstlisting}
Additionally, one may want to filter the content of the cache:
\begin{lstlisting}
sage: xi_cache=load_xis()
sage: val_one = {k : v for k, v in xi_cache.items() if v == 1}
sage: print_top_xis(val_one)
Stratum            | Residue Conditions           | xi^dim
----------------------------------------------------------------
(-2, -2, 1, 1)     | [[(0, 0)], [(0, 1)]]         | 1
\end{lstlisting}
To facilitate this, the method \verb|list_top_xis| picks apart the key of the $\xi$-cache:
\begin{lstlisting}
sage: for sigs, res_conds, value in list_top_xis():
....:     print('%
....:
((-2, -2, 1, 1),) (((0, 0),), ((0, 1),)) 1
((4,),) () 305/580608
\end{lstlisting}
\end{example}

%% file: sec_tests.tex
\section{Tests and Computations}
\label{sec:tests}

We use this section to briefly illustrate how the described methods of \verb|diffstrata| may be used to implement some of the key results and crosschecks of formulas in \cite{strataEC} 

\subsection{Euler Characteristics}
\label{subsec:ec}

We illustrate the methods described above to implement \cite[Thm. 1.3]{strataEC} for computing the Euler characteristics of strata. Recall that the (orbifold) Euler characteristic of the moduli space $\omoduli[g,n](\mu)$
is the dimension-weighted sum over all level graphs~$\Gamma \in \LG_L(B)$
without horizontal nodes
\[
\chi(B) \= (-1)^{d} \sum_{L=0}^d \sum_{\Gamma \in \LG_L(B)}  \ell_\Gamma
N_\Gamma^\top \int_{D_\Gamma} \prod_{i=0}^{L}  (\xi_\Gamma^{[i]})^{d_{\Gamma}^{[i]}}
\]
of the product of the top power of the of the first Chern class~$\xi_\Gamma^{[i]}$
of the tautological bundle at each level, where $d_{\Gamma}^{[i]}$ is the
dimension of the projectivized moduli space at level~$i$ and where
$d = \dim(B) = N-1$. The equivalence of the above formula and the one
stated in \cite[Eq.~(2)]{strataEC} follows from \cite[Lemma~9.12]{strataEC}.
\par
We may implement this using \verb|diffstrata| for a stratum \verb|X| as follows:
\begin{algorithm}[Euler Characteristic]\mbox{}
\begin{description}
\item[Step 1] Loop over \verb|L| from \verb|0| to \verb|X.dim()|.
\item[Step 2] Loop over \verb|ep| in \verb|X.enhanced_profiles_of_length(L)|.
\item[Step 3] Writing the profile \verb|ep=(p, enh)|, the number $l_\Gamma$ is the product over \verb|X.bics[p[i]].ell| for all \verb|i| and $N_\Gamma^\top$ is \verb|X.bics[p[0]].top.dim() + 1|.
\item[Step 4] Calculate the product using \verb|top_xi_at_level| (this uses $d_{\Gamma}^{[i]}$ and returns a number). It's also cached, cf. \autoref{sec:caching}.
\item[Step 5] Sum all these \emph{numbers} together.
\end{description}
\end{algorithm}

\begin{example}
The above algorithm is implemented by \verb|euler_characteristic|. We run this with an empty cache:
\begin{lstlisting}
sage: print_adm_evals()
Stratum            | Psis                         | eval
----------------------------------------------------------------
sage: print_top_xis()
Stratum            | Residue Conditions           | xi^dim
----------------------------------------------------------------
sage: X=Stratum((4,))
sage: %
CPU times: user 7.31 s, sys: 108 ms, total: 7.41 s
Wall time: 7.43 s
-55/504
\end{lstlisting}
Re-inspecting the cache, we see that it has been filled:
\begin{lstlisting}
sage: print_top_xis()
Stratum            | Residue Conditions           | xi^dim
----------------------------------------------------------------
(-4, -2, 4)        | [(0, 0), (0, 1)]             | 1
(-4, 0, 2)         | [(0, 0)]                     | 1
(-4, 1, 1)         | [(0, 0)]                     | 1
(-4, 4)            | [(0, 0)]                     | -15/8
(-3, -3, 4)        | [(0, 0), (0, 1)]             | 1
(-2, -2, -2, 4)    | [(0, 0), (0, 1), (0, 2)]     | -4
(-2, -2, -2, 4)    | [[(0, 0), (0, 2)], [(0, 1)]] | 1
(-2, -2, 0, 2)     | [(0, 0), (0, 1)]             | -2
(-2, -2, 1, 1)     | [[(0, 0)], [(0, 1)]]         | 1
(-2, -2, 1, 1)     | [(0, 0), (0, 1)]             | -1
(-2, -2, 2)        | [(0, 0), (0, 1)]             | 1
(-2, -2, 4)        | [[(0, 0)], [(0, 1)]]         | -11/12
(-2, -2, 4)        | [(0, 0), (0, 1)]             | 13/8
(-2, 0, 0)         | [(0, 0)]                     | 1
(-2, 0, 0, 0)      | [(0, 0)]                     | -1
(-2, 0, 2)         | [(0, 0)]                     | 1/8
(-2, 1, 1)         | [(0, 0)]                     | 0
(-2, 2)            | [(0, 0)]                     | -1/8
(-2, 4)            | [(0, 0)]                     | -23/1152
(0,)               | ()                           | 1/24
[(0,), (-2, 0, 0)] | [(1, 0)]                     | -1/24
[(0,), (0,)]       | ()                           | -1/576
[(0,), (0, 0)]     | ()                           | 0
(0, 0)             | ()                           | 0
(0, 0, 0)          | ()                           | 0
(0, 2)             | ()                           | 0
(1, 1)             | ()                           | 0
(2,)               | ()                           | -1/640
(4,)               | ()                           | 305/580608
sage: print_adm_evals()
Stratum            | Psis                         | eval
----------------------------------------------------------------
(0,)               | {1: 1}                       | 1/24
(-2, 2)            | {1: 1}                       | 1/8
(-2, 0, 0, 0)      | {1: 1}                       | 1
(-2, -2, -2, 4)    | {1: 1}                       | 1
(-2, -2, 1, 1)     | {1: 1}                       | 1
(-4, 4)            | {1: 1}                       | 5/8
(-2, -2, 0, 2)     | {1: 1}                       | 1
(0, 0)             | {1: 1, 2: 1}                 | 1/24
(-2, -2, 4)        | {1: 1, 2: 1}                 | 11/12
(-2, 0, 2)         | {1: 1, 2: 1}                 | 1/4
(-2, 1, 1)         | {1: 1, 2: 1}                 | 1/6
(0, 0, 0)          | {1: 1, 2: 1, 3: 1}           | 1/12
(-2, 4)            | {1: 1, 2: 2}                 | 73/1152
(0, 0, 0)          | {1: 1, 2: 2}                 | 1/12
(1, 1)             | {1: 1, 2: 3}                 | 1/720
(-2, -2, 4)        | {1: 1, 3: 1}                 | 11/12
(0, 0)             | {1: 2}                       | 1/24
(-2, -2, 4)        | {1: 2}                       | 19/24
(-2, 0, 2)         | {1: 2}                       | 1/8
(-2, 1, 1)         | {1: 2}                       | 1/24
(-2, 4)            | {1: 2, 2: 1}                 | 97/1152
(1, 1)             | {1: 2, 2: 2}                 | 1/720
(2,)               | {1: 3}                       | 1/1920
(0, 0, 0)          | {1: 3}                       | 1/24
(-2, 4)            | {1: 3}                       | 43/1152
(0, 2)             | {1: 4}                       | 11/1920
(1, 1)             | {1: 4}                       | 1/720
(4,)               | {1: 5}                       | 13/580608
(-4, 4)            | {2: 1}                       | 5/8
(-2, 2)            | {2: 1}                       | 1/8
(-2, 0, 0, 0)      | {2: 1}                       | 1
(-2, 1, 1)         | {2: 1, 3: 1}                 | 1/6
(-2, 0, 2)         | {2: 2}                       | 1/4
(-2, 1, 1)         | {2: 2}                       | 1/6
(-2, 4)            | {2: 3}                       | 19/1152
(-2, -2, 0, 2)     | {3: 1}                       | 1
(-2, -2, 1, 1)     | {3: 1}                       | 1
(-2, -2, 4)        | {3: 2}                       | 7/24
(-2, -2, -2, 4)    | {4: 1}                       | 1
\end{lstlisting}
Of course this affects any future calculations:
\begin{lstlisting}
sage: %
CPU times: user 2.52 ms, sys: 5.39 ms, total: 7.91 ms
Wall time: 22.3 ms
-55/504
\end{lstlisting}
In fact, the cache was already used in the first calculation, as all levels appear already in the two and three-level graphs (see \autoref{leveltypes}).
\end{example}

\begin{example}
Note that calling \verb|euler_characteristic| is simply a frontend for calling the method \verb|euler_char_immediate_evaluation|. Calling this directly, we may use the \verb|quiet=False| option to give extensive output:
\begin{lstlisting}
sage: X=Stratum((2,))
sage: X.euler_char_immediate_evaluation(quiet=False)
Generating enhanced profiles of length 0...
Building all graphs in () (1/1)...
Going through 1 profiles of length 0...
1 / 1, ((), 0): Calculating xi at level 0 (cache) -1/640 Done.
Generating enhanced profiles of length 1...
Building all graphs in (0,) (1/2)...
Building all graphs in (1,) (2/2)...
Going through 2 profiles of length 1...
1 / 2, ((0,), 0): Calculating xi at level 0 (cache) 1/24 level 1 (cache) -1/8 Done.
2 / 2, ((1,), 0): Calculating xi at level 0 (cache) 0 Product 0. Done.
Generating enhanced profiles of length 2...
Building all graphs in (0, 1) (1/1)...
Going through 1 profiles of length 2...
1 / 1, ((0, 1), 0): Calculating xi at level 0 (cache) 1/24 level 1 (cache) 1 level 2 (cache) 1 Done.
Generating enhanced profiles of length 3...
Going through 0 profiles of length 3...
-1/40
\end{lstlisting}
\end{example}

\begin{example}
Alternatively, we can calculate the Chern character of the logarithmic cotangent bundle using \cite[Thm. 1.2]{strataEC} and use Newton's identity to calculate the Chern polynomial. Of course, this is a longer calculation and there is less caching (\verb|top_xi_at_level| is not called), but it may be used to check the consistency of the formulas:
\begin{lstlisting}
sage: X=Stratum((2,))
sage: X.top_chern_class().evaluate()
1/40
sage: X.euler_char()
-1/40
\end{lstlisting}
Note that \verb|euler_char| is simply a frontend for various methods computing the Chern polynomial. More details may be found in the docstrings of the various methods.
\end{example}

\subsection{Crosschecks}

The module \verb|tests| includes some more cross-checks and example computations using \verb|diffstrata|. For example, \verb|leg_tests| tests on each one-dimensional graph $\Gamma$ of a stratum if the evaluation of the $\xi_\Gamma^{[i]}$ at that level is the same (for every leg!) as the product of $\Gamma$ with $\xi$. Note that these expressions can be evaluated and the numbers compared.

The method \verb|commutativity_check| runs an extensive commutativity check on a stratum, i.e. multiplying products of BICs in various orders to give top-level classes and check that these evaluate to the same number, testing the normal bundle and intersection formulas along the way.

Finally, the class \verb|BananaSuite| tests the strata $\Omega\cM_1(k,1,-k-1)$ (cf. \cite[\S 10.3]{strataEC}) implementing the $D$-notation introduced there and including a method to test \cite[Prop. 10.2]{strataEC}.

\begin{example}
We illustrate the tests on some small strata:
\begin{lstlisting}
sage: leg_test((4,))
Graph ((3, 6, 7, 2), 0): xi evaluated: 1/48 (dim of Level 0: 1)
level: 0, leg: 1, xi ev: 1/48
Graph ((3, 6, 7, 2), 1): xi evaluated: 1/24 (dim of Level 0: 1)
level: 0, leg: 1, xi ev: 1/24
Graph ((3, 6, 5, 2), 0): xi evaluated: 1/48 (dim of Level 0: 1)
level: 0, leg: 1, xi ev: 1/48
Graph ((3, 6, 5, 4), 0): xi evaluated: 1/48 (dim of Level 0: 1)
level: 0, leg: 1, xi ev: 1/48
sage: commutativity_check((2,))
Starting IPs
(0, 0, 0)
0 0
Starting IPs
(0, 0, 1)
0 0
Starting IPs
(0, 1, 0)
0 1
Starting IPs
(0, 1, 1)
0 1
Starting IPs
(1, 0, 0)
1 0
Starting IPs
(1, 0, 1)
1 0
Starting IPs
(1, 1, 0)
1 1
Starting IPs
(1, 1, 1)
1 1
sage: B=BananaSuite(2)
sage: B.check()
D(1,1)^2 = -1, RHS = -1
D(1,2)^2 = -1, RHS = -1
D(5,1)^2 = -3/2, RHS = -3/2
True
\end{lstlisting}
\end{example}